\documentclass[10pt,a4paper]{article}
\usepackage{amssymb,latexsym,amsmath,amsfonts,amsthm}
\usepackage{graphicx}

\usepackage{epsfig}
\usepackage{comment}
\usepackage{subfigure}
\usepackage{graphicx,ebezier}
\usepackage{overpic}

\DeclareMathOperator{\diag}{diag}
\DeclareMathOperator{\Res}{Res}

\newcommand{\lam}{\lambda}

\newcommand{\bewijs}{\textsc{proof}}
\newcommand{\bol}{\hfill\square\\}
\newcommand{\til}{\tilde}

\newcommand{\vecm}{\mathbf{m}}
\newcommand{\vecn}{\mathbf{n}}
\newcommand{\vece}{\mathbf{e}}
\newcommand{\vecA}{\mathbf{A}}
\newcommand{\crit}{\textrm{crit}}

\newcommand{\Pee}{{\cal P}}
\renewcommand{\Re}{\mathrm{Re}\,}
\renewcommand{\Im}{\mathrm{Im}\,}
\newcommand{\ds}{\displaystyle}

\newtheorem{theorem}{Theorem}[section]
\newtheorem{lemma}[theorem]{Lemma}
\newtheorem{proposition}[theorem]{Proposition}

\newtheorem{corollary}[theorem]{Corollary}

\theoremstyle{definition}
\newtheorem{definition}[theorem]{Definition}

\theoremstyle{remark}

\newtheorem{remark}[theorem]{Remark}

\newtheorem{convention}[theorem]{Convention}

\numberwithin{equation}{section}

\title{A phase transition for non-intersecting Brownian motions, and the Painlev\'e~II equation}

\author{Steven Delvaux\footnotemark[1],\quad Arno B.~J.~Kuijlaars\footnotemark[1]}
\date{}

\begin{document}

\maketitle
\renewcommand{\thefootnote}{\fnsymbol{footnote}}
\footnotetext[1]{Department of Mathematics, Katholieke Universiteit Leuven,
Celestijnenlaan 200B, B-3001 Leuven, Belgium. email:
\{steven.delvaux,arno.kuijlaars\}\symbol{'100}wis.kuleuven.be. \\
The first author is a Postdoctoral Fellow of the Fund for Scientific Research -
Flanders (Belgium). His work is supported by the Onderzoeksfonds
K.U.Leuven/Research Fund
K.U.Leuven. \\
The work of the second author is supported by FWO-Flanders project
G.0455.04, by K.U. Leuven research grants OT/04/21 and OT/08/33, by the Belgian
Interuniversity Attraction Pole P06/02, by the European Science Foundation
Program MISGAM, and by a grant from the Ministry of Education and Science of
Spain, project code MTM2005- 08648-C02-01.}

\begin{abstract}
We consider $n$ non-intersecting Brownian motions with two fixed starting
positions and two fixed ending positions in the large $n$ limit. We show that
in case of \lq large separation\rq\ between the endpoints, the particles are
asymptotically distributed in two separate groups, with no interaction between
them, as one would intuitively expect. We give a rigorous proof using the
Riemann-Hilbert formalism. In the case of \lq critical separation\rq\ between
the endpoints we are led to a model Riemann-Hilbert problem associated to the
Hastings-McLeod solution of the Painlev\'e~II~equation. We show that the
Painlev\'e~II equation also appears in the large $n$ asymptotics of the
recurrence coefficients of the multiple Hermite polynomials that are associated
with the Riemann-Hilbert problem.


\end{abstract}

\section{Introduction and statement of results}
\label{sectionintro}

\subsection{Non-intersecting Brownian motion}
\label{subsectionintrointro}

This paper deals with non-intersecting Brownian motions with two
fixed starting positions $a_1$, $a_2$ and two fixed ending
positions $b_1$, $b_2$. Let us first describe the general framework
in which this paper fits.


Let $p,q\in\mathbb N $. Consider sequences of real numbers
$\{a_k\}_{k=1}^p$, $\{b_l\}_{l=1}^q$ and sequences of positive
integers $\{n_k\}_{k=1}^p$, $\{m_l\}_{l=1}^q$ satisfying
\begin{itemize}
  \item $a_1>\cdots>a_p$,
  \item $b_1>\cdots>b_q$,
  \item $\sum_{k=1}^p n_k = \sum_{l=1}^q m_l =: n$.
\end{itemize}
Consider $n$ one-dimensional Brownian motions (actually
Brownian bridges) such that $n_k$ of them start from the
point $a_k$ at time $t=0$, for $k=1,\ldots,p$, and $m_l$ of them
arrive in the point $b_l$ at time $t=1$, for $l=1,\ldots,q$, and
such that the $n$ particles are conditioned not to collide with
each other in the time interval $t \in (0,1)$.

We are interested in the limiting behavior where the number of
Brownian particles $n$ tends to infinity in such a way that each
of the fractions $n_k/n$, $k= 1, \ldots, p$ and $m_l/n$,
$l=1,\ldots,q$ has a limit in $(0,1)$. To obtain non-trivial
limiting behavior we assume that the transition probability density
of the Brownian motions scales as
\begin{equation}\label{transitionprob}
    P_N(t,x,y) = \frac{\sqrt{N}}{\sqrt{2\pi t}} e^{-\frac{N}{2t}(x-y)^2}
\end{equation}
where $N$ is a parameter (the inverse of the overall variance)
that increases as $n$ increases. In a non-critical regime we will simply
take $N = n$.

\begin{figure}[tbp]
\begin{center}
\includegraphics[scale=0.6]{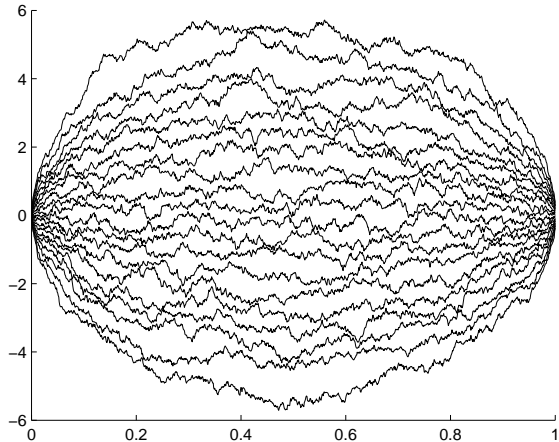}
\end{center}
\caption{Non-intersecting Brownian motions in case of one starting and ending
point. Here the horizontal axis denotes the time, $t\in[0,1]$, and for each
fixed $t$ the positions of the $n$ non-intersecting Brownian motions at time
$t$ are denoted on the vertical line through $t$.} \label{figwigner}
\end{figure}

Under these assumptions, it is expected that the particles will asymptotically
fill a bounded region in the time-space plane ($tx$-plane). It is also expected
that for every $t\in(0,1)$, the particles at time $t$ are asymptotically
distributed according to some well-defined limiting distribution. See Figures
\ref{figwigner} and \ref{figbrownianmotions} for possible behaviors when
$p=q=1$ and $p=q=2$, respectively.

It is natural to ask for an explicit description of this limiting distribution.
Such a result is known in the classical case where $p=q=1$, which is connected
to Dyson's Brownian motion \cite{Dyson}. Here it is known that up to a suitable
scaling and translation, the Brownian particles at time $t\in(0,1)$ have the
same joint distribution as the eigenvalues of the Gaussian unitary ensemble
(GUE) of size $n\times n$, and the limiting distribution as $n\to\infty$ is
given by Wigner's semicircle law. More precisely, if all Brownian motions start
from $a:=a_1$ at time $t=0$ and end up at  $b:=b_1$ at time $t=1$, and if $N
=n/T$, then the non-intersecting Brownian particles at time $t \in (0,1)$ are
asymptotically distributed on the interval $[\alpha(t),\beta(t)]$ with
endpoints
\begin{align}
    \label{defalpha0} \alpha(t) & = (1-t)a+tb-\sqrt{4 T t(1-t)},\\
    \label{defbeta0} \beta(t) & =  (1-t)a+tb+\sqrt{4 T t(1-t)},
\end{align}
and with limiting density of the particles given by the semicircle
law on that interval
\begin{equation}\label{wignerdensity}
  \frac{1}{2\pi T t(1-t)}\sqrt{(\beta(t)-x)(x-\alpha(t))},\qquad
   x\in[\alpha(t),\beta(t)].
\end{equation}
By varying $t \in [0,1]$, the endpoints \eqref{defalpha0}, \eqref{defbeta0}
parameterize an ellipse in the time-space plane, see
Figure~\ref{figwigner}.

Explicit descriptions of the limiting distribution of the non-intersecting
Brownian motions are also known when $p=1$, $q>1$. Also in this case there is
an underlying random matrix model \cite{ABK,BK2,BK3,KKPS,O,Pastur}.

The limiting distribution is also known when $p=q=2$, $m_1=n_1=m_2=n_2=n/2$, $N
=n$, and $(a_1-a_2)(b_1-b_2)<2$. In this case the limiting distribution as
$n\to\infty$ of the Brownian particles at time $t\in(0,1)$ is obtained from
a certain algebraic curve of degree four \cite{DKV}. In contrast to the
previous cases, however, it is not known if the Brownian particles for finite
$n$ can be described as the eigenvalues of a random matrix ensemble. See
Figure~\ref{figsmallsep} for an illustration of this case.

\subsection{Separation of the endpoints}

The main results of this paper will concern non-intersecting Brownian motions
with two starting points $a_1>a_2$ at time $t=0$, two ending points $b_1>b_2$
at time $t=1$, and in addition
\begin{align}
\label{balanced} n_1 = m_1,\qquad n_2=m_2,\quad (\textrm{for all }n)
\end{align}
and if we put
\begin{align} \label{defp1p2}
    p_1 = \frac{n_1}{n}, \qquad p_2 = \frac{n_2}{n}
\end{align}
which are varying with $n$, then we assume that
\begin{align}
\label{defp1} \frac{n_1}{n} & =
    p_1 = p_1^*+O(1/n),\quad n\to\infty,  \\
\label{defp2} \frac{n_2}{n}  &=
    p_2 = p_2^*+O(1/n),\quad n\to\infty,
\end{align}
for certain limiting values $p_1^*,p_2^* \in (0,1)$. We also assume that $N$
increases with $n$ such that
\begin{align}
\label{defnN}
    T = n/N  > 0
    \end{align}
is fixed.

\begin{figure}[tbp]
\begin{center}
\subfigure{\label{figlargesep}}\includegraphics[scale=0.3]{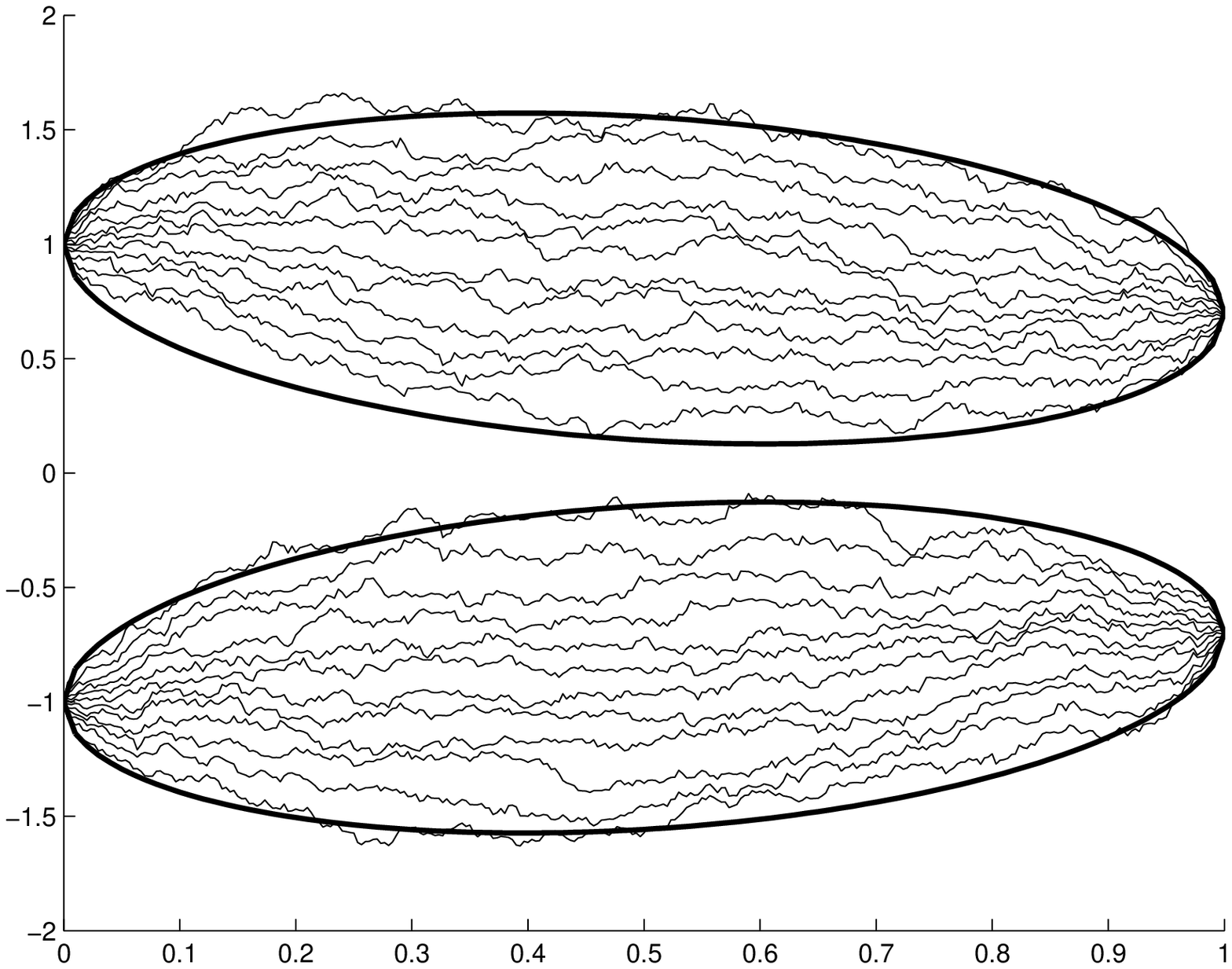}\hspace{5mm}
\subfigure{\label{figsmallsep}}\includegraphics[scale=0.3]{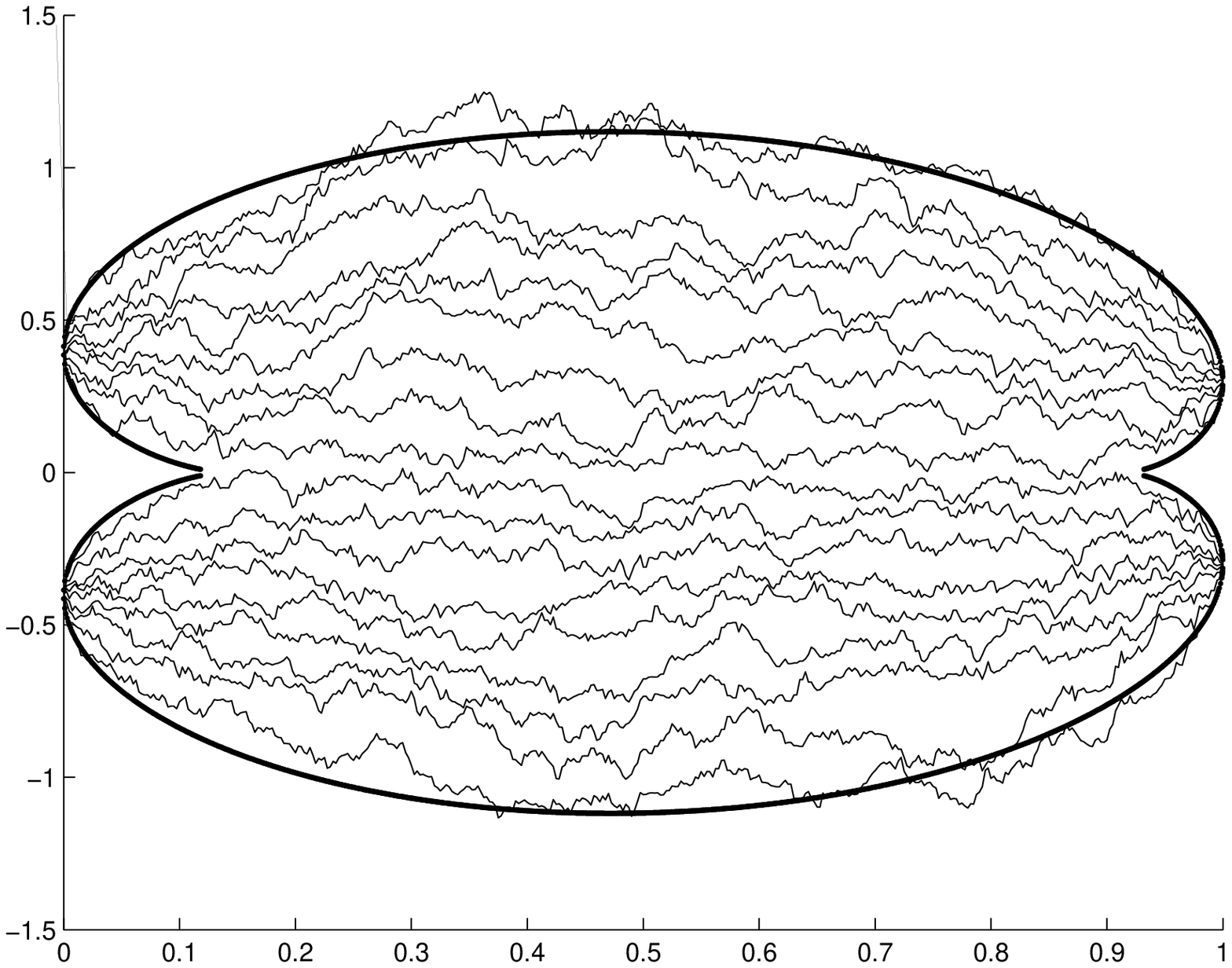}
\subfigure{\label{figcriticalsep}}\includegraphics[scale=0.3]{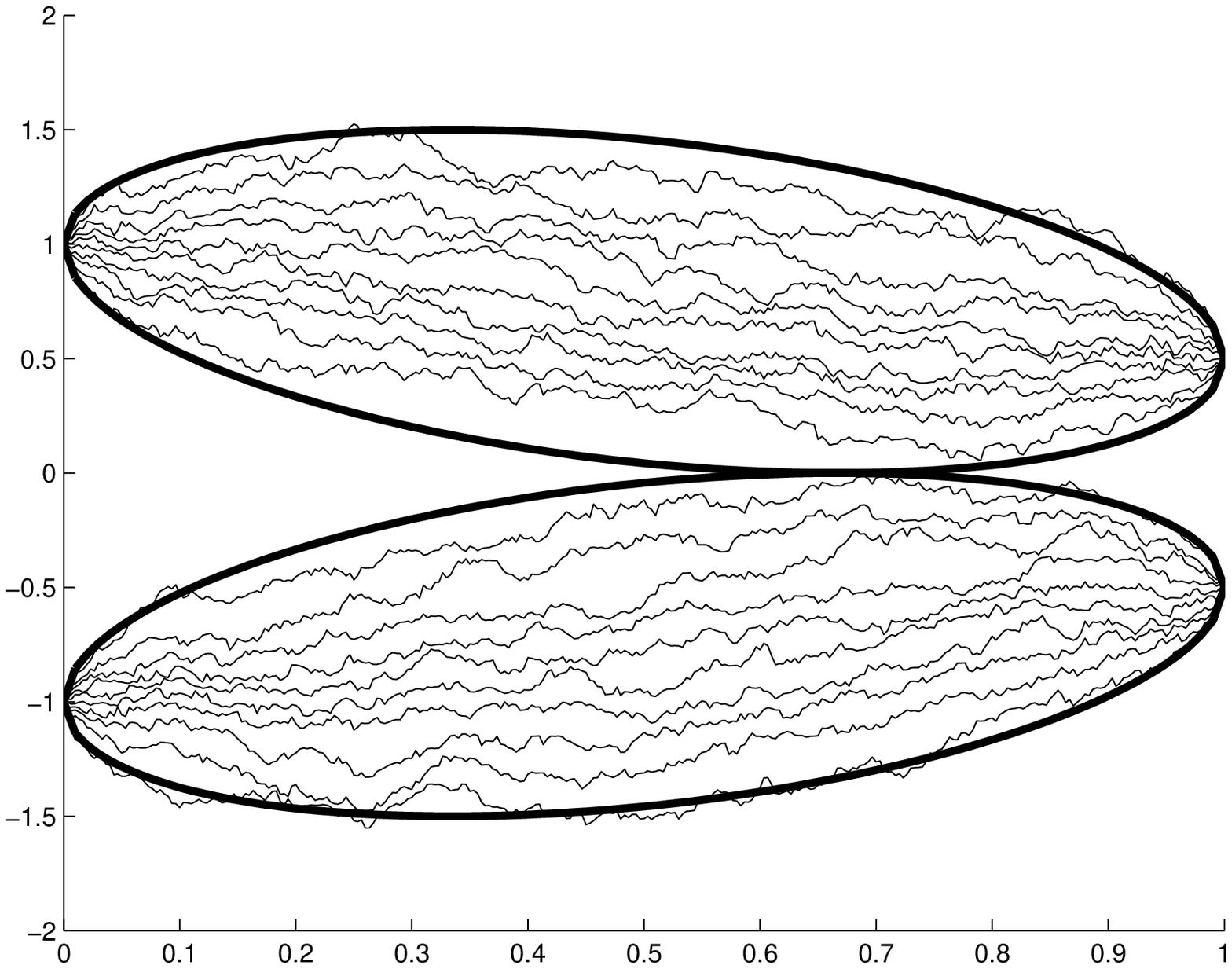}
\end{center}
\caption{Non-intersecting Brownian motions with two starting and two ending
positions in case of (a) large, (b) small, and (c) critical separation between
the endpoints. Here the horizontal axis denotes the time, $t\in[0,1]$, and for
each fixed $t$ the positions of the $n$ non-intersecting Brownian motions at
time $t$ are denoted on the vertical line through $t$. Note that for
$n\to\infty$ the positions of the Brownian motions fill a prescribed region in
the time-space plane, which is bounded by the boldface lines in the figures.
Here we have chosen $N=n=20$ and $p_1=p_2=1/2$ in each of the figures, and (a)
$a_1=-a_2=1$, $b_1=-b_2=0.7$, (b) $a_1=-a_2=0.4$, $b_1=-b_2=0.3$, and (c)
$a_1=-a_2=1$, $b_1=-b_2=1/2$, in the cases of large, small and critical
separation, respectively.} \label{figbrownianmotions}
\end{figure}

If the assumptions \eqref{balanced}--\eqref{defnN} hold, and if the separation
between the starting points $a_1$ and $a_2$, and the ending points $b_1$ and
$b_2$ is large enough, then in analogy with
\eqref{defalpha0}--\eqref{wignerdensity}, one would expect the Brownian
particles to be asymptotically distributed on two disjoint ellipses in the
$tx$-plane, whose intersections with the vertical line through $t$ are given by
the two intervals $[\alpha_1^*,\beta_1^*]$ and $[\alpha_2^*,\beta_2^*]$ with
\begin{align} \label{defalpha}
    \alpha_j^* & = \alpha_j^*(t) =  (1-t)a_j+tb_j-\sqrt{4p_j^* T t(1-t)},\\
    \label{defbeta}
    \beta_j^* & =  \beta_j^*(t) = (1-t)a_j+tb_j+\sqrt{4p_j^* T t(1-t)},
\end{align}
for $j=1,2$,
and with limiting densities on these two intervals given by the semicircle laws
\begin{align}
    \label{wignerdensities}
    \frac{1}{2\pi T t(1-t)}\sqrt{(\beta_j^*-x)(x-\alpha_j^*)},\quad
    x\in[\alpha_j^*,\beta_j^*],
        \quad j=1,2,
\end{align}
for each $t\in(0,1)$. This situation is illustrated in Figure
\ref{figlargesep}. Note that
\[ \frac{1}{2\pi T t(1-t)} \int_{\alpha_j^*}^{\beta_j^*}  \sqrt{(\beta_j^*-x)(x-\alpha_j^*)} dx = p_j^*. \]
We derive the precise condition for this two-ellipse scenario to happen. It is
clear that a necessary condition is the disjointness of the two ellipses.
\begin{lemma} \label{lemmadisjointellipses} (Disjointness of the two ellipses)
The two ellipses parameterized by \eqref{defalpha}--\eqref{defbeta} are
disjoint if and only if
\begin{equation}\label{largeseparationeq}
    (a_1-a_2)(b_1-b_2) > T \left(\sqrt{p_1^*}+\sqrt{p_2^*}\right)^2.
\end{equation}
\end{lemma}
\bewijs.
The two ellipses are disjoint if and only if
$\alpha_1^*(t)>\beta_2^*(t)$ for all $t\in (0,1)$. From
\eqref{defalpha}--\eqref{defbeta} this leads to the condition
\begin{align} \label{conditiondisjointellipses}
(1-t)(a_1-a_2)+t(b_1-b_2)  > \sqrt{4Tt(1-t)}\left(\sqrt{p_1^*}+\sqrt{p_2^*}\right),
\quad t \in (0,1),
\end{align}
which after putting $u = \sqrt{\frac{t}{1-t}}$ is equivalent to
\begin{align*}
        u^2(b_1-b_2) -2 u \sqrt{T} \left(\sqrt{p_1^*}+\sqrt{p_2^*}\right) +(a_1-a_2) >0, \qquad u \in
        (0,\infty).
\end{align*}
The left-hand side is a quadratic expression in
$u$, whose discriminant is negative if and only if \eqref{largeseparationeq}
holds. The lemma then easily follows since $b_1 > b_2$ and $a_1 > a_2$. $\bol$

 We will call \eqref{largeseparationeq} the case
of large separation.

\begin{definition} (Large, critical and small separation)
For each $n$, consider $n$ non-intersecting Brownian motions with two starting
points $a_1>a_2$ at time $t=0$ and two ending points $b_1>b_2$ at time $t=1$,
and assume that the hypotheses \eqref{balanced}--\eqref{defnN} hold. If
\eqref{largeseparationeq} holds, we say that we are in a situation of
\emph{large separation} of the endpoints. If instead
\begin{equation}\label{smallseparationeq}
(a_1-a_2)(b_1-b_2) < T (\sqrt{p_1^*}+\sqrt{p_2^*})^2,
\end{equation}
we are in a situation of \emph{small separation}, and if
\begin{equation}\label{criticalseparationeq}
(a_1-a_2)(b_1-b_2) = T(\sqrt{p_1^*}+\sqrt{p_2^*})^2,
\end{equation}
we are in a situation of \emph{critical separation} between the endpoints. In
the latter case, we define the \emph{critical time} $t_{\crit}\in(0,1)$ as
\begin{equation}\label{deftcrit0}
t_{\crit} = \frac{a_1-a_2}{(a_1-a_2)+(b_1-b_2)}.
\end{equation}
\end{definition}
\medskip

The non-intersecting Brownian motions corresponding to each of these three
cases are illustrated in Figure \ref{figlargesep}--\ref{figcriticalsep}. In the
case of critical separation, the time $t_{\crit}$ is precisely the time where
the two ellipses with endpoints parameterized by
\eqref{defalpha}--\eqref{defbeta} are tangent, cf.\ Figure
\ref{figcriticalsep}.

We may equivalently view the critical behavior in terms of the temperature
parameter $T =  n/N$.
With $a_j$, $b_j$, and $p_j^*$ fixed for $j=1,2$, there is a critical temperature
\[ T_{\crit} = \frac{(a_1-a_2)(b_1-b_2)}{(\sqrt{p_1^*} + \sqrt{p_2^*})^2} \]
such that $T < T_{\crit}$ (low temperature), $T > T_{\crit}$ (high
temperature), and $T = T_{\crit}$ (critical temperature) correspond to large,
small, and critical separation, respectively.

\subsection{Large separation: decoupling of the Brownian motions}

The case of small separation of the endpoints was considered in \cite{DKV} for
the special case when $p_1=p_2=1/2$. In the present paper, we will focus
instead on the cases of large and critical separation. As is  to be expected,
in these regimes the Brownian motions asymptotically decouple into two separate
groups, as is evidenced by Figures \ref{figlargesep} and \ref{figcriticalsep}.
This is our first main theorem.

\begin{theorem}\label{st2ellipses} (Decoupling of Brownian motions)
Consider $n$ non-intersecting Brownian motions with two starting
points $a_1>a_2$ at time $t=0$ and two ending points $b_1>b_2$ at
time $t=1$. Assume that the hypotheses
\eqref{balanced}--\eqref{defnN} hold, and assume that either
\begin{itemize}
    \item there is large separation of the endpoints
    \eqref{largeseparationeq}, or
    \item there is critical separation
    of the endpoints \eqref{criticalseparationeq} and $t\neq t_{\crit}$,
    where $t_{\crit}$ is given by \eqref{deftcrit0}.
\end{itemize}
Then as $n  \to\infty$ the Brownian
particles at time $t\in(0,1)$ are asymptotically supported on the
two disjoint intervals $[\alpha_1^*(t),\beta_1^*(t)]$ and
$[\alpha_2^*(t),\beta_2^*(t)]$ given by \eqref{defalpha}--\eqref{defbeta},
with limiting densities given by the semicircle laws
\eqref{wignerdensities}.
\end{theorem}

In the case of critical separation, we strongly expect that the
conclusion of the theorem should also be valid when $t=t_{\crit}$.
However, we will not consider the critical time in this paper.

Our proof of Theorem \ref{st2ellipses} follows from a steepest descent analysis
of the underlying Riemann-Hilbert problem that will be described in Section
\ref{subsectionintroRH}. The Riemann-Hilbert problem is of size $4 \times 4$
and it was also used in \cite{DKV} to analyze the case of small separation. In
the (apparently) simpler case of large separation one expects that for large
$n$, the $4 \times 4$ RH problem asymptotically decouples into two smaller RH
problems of size $2 \times 2$. This is indeed the case, but in order to show
this, we need a preliminary transformation where we introduce auxiliary curves
in the complex plane and subsequently perform a Gaussian elimination step in
the jump matrix of the RH problem, serving to annihilate some undesired entries
of this matrix. This Gaussian elimination step is similar to the so-called
global opening of the lens discussed in \cite{ABK,AKV,DKV,LW}.

It follows from our analysis that the interaction between the two
groups of Brownian particles decays exponentially with $n$
in case of large separation, and polynomially (like a power
$n^{-1/3}$) in case of critical separation at a non-critical time.
In the limit when $n\to\infty$ the particles will then indeed be
distributed inside two disjoint ellipses in the time-space plane.

\subsection{Riemann-Hilbert problem}
\label{subsectionintroRH}

The non-intersecting Brownian motions described in the previous subsections are
related to the following Riemann-Hilbert problem (RH problem) which we already
alluded to above. The RH problem was introduced in \cite{DK2} as a
generalization of the RH problem for orthogonal polynomials in \cite{FIK}, see
also \cite{VAGK}.  In accordance with Section \ref{subsectionintrointro} we
will state the RH problem for general numbers of starting and ending positions
$p,q$ of the Brownian motions, although in our applications we will eventually
take $p=q=2$.

Define weight functions
\begin{align}
\label{gaussianw1} w_{1,k}(x) & = e^{-\frac{N}{2t}(x^2-2a_kx)},\quad k=1,\ldots,p, \\
\label{gaussianw2} w_{2,l}(x) & = e^{-\frac{N}{2(1-t)}(x^2-2b_lx)},\quad
l=1,\ldots,q.
\end{align}

The RH problem consists in finding a matrix-valued function
$Y(z)=$\\$Y_{n_1,\ldots,n_p;m_1,\ldots,m_q}(z)$ of size $p+q$ by $p+q$ such
that
\begin{itemize}
\item[(1)] $Y(z)$ is analytic
in $\mathbb C \setminus\mathbb R $;
 \item[(2)] For $x\in\mathbb R $, it holds that
\begin{equation}\label{defjumpmatrix0}
Y_{+}(x) = Y_{-}(x)
\begin{pmatrix} I_p & W(x)\\
0 & I_q
\end{pmatrix},
\end{equation}
where $I_k$ denotes the identity matrix of size $k$; where $W(x)$ denotes the
rank-one matrix (outer product of two vectors)
\begin{equation}\label{defWblock0}
W(x) = \begin{pmatrix} w_{1,1}(x) \\ \vdots \\ w_{1,p}(x)
\end{pmatrix}\begin{pmatrix} w_{2,1}(x) & \ldots & w_{2,q}(x)
\end{pmatrix},
\end{equation}
and where the notation $Y_+(x), Y_-(x)$ denotes the limit of $Y(z)$ with $z$
approaching $x\in\mathbb R $ from the upper or lower half plane in $\mathbb C $,
respectively;
  \item[(3)] As $z\to\infty$, we have that
\begin{equation}\label{asymptoticconditionY0} Y(z) =
    (I_{p+q}+O(1/z))\diag(z^{n_1},\ldots,z^{n_p},z^{-m_1},\ldots,z^{-m_q}).
\end{equation}
\end{itemize}

The RH problem has a unique solution \cite{DK2} that can be described in terms
of certain multiple orthogonal polynomials (actually multiple Hermite
polynomials); details will be given in Section \ref{subsectionintromop}.

Let us explain the connection between the non-intersecting Brownian motions and
the RH problem in the case $p=q=2$. It is well-known \cite{KMcG}, see also
\cite{Joh1,Kat1}, that the distribution of the non-intersecting Brownian
motions at time $t\in (0,1)$ describes a determinantal point process,
determined by an associated correlation kernel. According to \cite{DK2} the
correlation kernel $K(x,y)= K_{n_1,n_2;m_1,m_2}(x,y)$ can be expressed in terms
of the solution to the RH problem as
\begin{equation} \label{correlationkernel}
    K(x,y) = \frac{1}{2\pi i(x-y)}\begin{pmatrix} 0 & 0 & w_{2,1}(y) &
    w_{2,2}(y)\end{pmatrix} Y_{+}^{-1}(y)Y_{+}(x)\begin{pmatrix} w_{1,1}(x)\\
    w_{1,2}(x)\\ 0 \\ 0 \end{pmatrix}.
\end{equation}

By general properties of determinantal point processes,
Theorem~\ref{st2ellipses} then comes down to the statement that under the
conditions of Theorem~\ref{st2ellipses}, the limit of $\frac{1}{n}
K_{n_1,n_2;n_1,n_2}(x,x)$ as $n \to \infty$ exists and is equal to
\begin{equation} \label{correlationkernellimit}
    \lim_{n\to\infty}\frac{1}{n} K_{n_1,n_2;n_1,n_2}(x,x)
        = \frac{1}{2\pi Tt(1-t)} \sqrt{(\beta_j^*-x)(x-\alpha_j^*)},
            \qquad x \in [\alpha_j^*, \beta_j^*],
\end{equation}
for $j=1,2$. This is what we will show in Sections
\ref{sectionsteepestdescentlargesep} and \ref{sectionsteepestdescentcritical}.

Our method will also allow us to obtain the local scaling limits
of the correlation kernel that are common in random matrix theory
and related areas, namely the sine kernel in the bulk and the
Airy kernel at the endpoints of the intervals $[\alpha_1^*, \beta_1^*]$
and $[\alpha_2^*, \beta_2^*]$. We will not go into details about this
in this paper.

\subsection{Critical separation and the double scaling limit}
\label{subsectionintromodelRHPII}

Assume that the endpoints are such that
\begin{align}
    \label{deftau} (a_1-a_2)(b_1-b_2) = (\sqrt{p_1^*}+\sqrt{p_2^*})^2.
\end{align}
We put
\[ N = n/T \]
where $T$ can be interpreted as a temperature variable. For $T = 1$ (critical
temperature) we have by \eqref{criticalseparationeq}, \eqref{deftau} and
Theorem~\ref{st2ellipses} that in the large $n$ limit, the particles fill out
two ellipses, which are tangent to each other at the critical time $t_{\crit}$,
see Figure~\ref{figcriticalsep}. By varying $T$ around the critical value
$T_{\crit} = 1$ we move from a case of disjoint ellipses (for $T < 1$) to a
case of small separation (for $T > 1$), where the two-ellipses scenario is not
valid anymore. Hence we see a phase transition in the case of critical
separation, which is clearly seen at the critical time $t_{\crit}$. At a
non-critical time $t \neq t_{\crit}$ the phase transition is less obvious, but
there is also a nontrivial transitional effect. Indeed, the distribution of
particles at time $t \neq t_{\crit}$ in the case of small separation differs
from the distribution of two semicircle laws of two disjoint intervals, which
it is for large separation. Hence the endpoints of the intervals do not depend
analytically on the starting and ending points, which indicates the phase
transition. It is a surprising outcome of our analysis that the phase
transition for the case where $t\neq t_{\crit}$ can be described by the
Painlev\'e II equation.

In the case of critical separation we investigate the behavior of the Brownian
particles in a double scaling limit where $n, N \to \infty$, and simultaneously
$T \to T_{\crit} = 1$. More precisely, we consider the endpoints $a_1$, $a_2$,
$b_1$, $b_2$ fixed so that \eqref{deftau} holds for all $n$.
The temperature $T = T_n = n/N$ is varying with $n$ as follows:
\begin{equation}\label{doublescalinglimit}
    T_n = 1 + L n^{-2/3},
\end{equation}
where $L$ is an arbitrary real constant.

We will show that in the double scaling regime described above, the steepest
descent analysis of the Riemann-Hilbert problem leads in a natural way to a
model Riemann-Hilbert problem related to the Painlev\'e II equation.
More precisely, we will be led to the construction of a local parametrix that
can be mapped onto the model RH problem \cite{FN} satisfied by the
$\Psi$-functions (Lax pair) associated with the Hastings-McLeod solution of the
Painlev\'e II equation
\begin{equation}\label{defPII}
    q''(s) = sq(s)+2q^3(s).
\end{equation}
The Hasting-McLeod solution \cite{HML} is the special solution $q(s)$ of
\eqref{defPII} which is real for real $s$ and satisfies $q(s)\sim
\textrm{Ai}(s)$ as $s\to\infty$, where $\textrm{Ai}$ denotes the usual Airy
function. The precise form of the model RH problem will be described in Section
\ref{subsectionlocalparametrixx0star}.

The Hastings-McLeod solution of the Painlev\'e II equation also appears in
the famous Tracy-Widom distributions
\cite{TW1,TW2} for the largest eigenvalues of large
random matrices.
It also appears in the critical unitarily invariant matrix models,
where the parameters in the model
are fine-tuned so that the limiting mean eigenvalue density vanishes
quadratically at an interior point of its support \cite{BI2,Claeys1}. In this
case it leads to a new family of local scaling limits of the eigenvalue
correlation kernel that involve the $\Psi$-functions associated with $q(s)$.

In our situation the Painlev\'e II equation does not manifest
itself in the local scaling limits of the correlation kernel. The
construction of the local parametrix is done at a point $x_0$
strictly outside of the support and it does not influence the
local correlation functions for the positions of any of the
particles. The point $x_0$ does not seem to have any physical
meaning.

We emphasize that our asymptotic analysis will be only valid when
$t\neq t_{\crit}$. At the critical time $t=t_{\crit}$ where the
two ellipses are tangent, one is led to a considerably more
difficult, multi-critical situation. Here one expects the
appearance of a model RH problem related to some as yet unknown
fourth order ODE. As already mentioned, we will not attempt to
study this case in the present paper.

\subsection{Generalities on multiple orthogonal polynomials}
\label{subsectionintromop}

While the appearance of the Hastings-McLeod solution of the RH problem does not
affect any of the local scaling limits, it is felt by the recurrence
coefficients of the multiple Hermite polynomials. These polynomials appear in
the solution of the RH problem given in Section \ref{subsectionintroRH}.

To state the results, let us first recall some generalities on multiple
orthogonal polynomials in the sense of \cite{DK2}. In accordance with Sections
\ref{subsectionintrointro} and \ref{subsectionintroRH} we will again give the
definitions for general values of $p$ and $q$, although in our applications we
will eventually take $p=q=2$.

\begin{definition}\label{defmop} (Multiple orthogonal polynomials; cf.\ \cite{DK2})
Let $p,q\in\mathbb N $ be two positive integers. Let there be given
\begin{itemize}
\item A (finite) sequence of positive integers $n_1,n_2,\ldots,n_p\in\mathbb N $;
\item A sequence of weight functions $w_{1,1}(x),w_{1,2}(x),\ldots,w_{1,p}(x):\mathbb R \to\mathbb R $;
\item A sequence of positive integers $m_1,m_2,\ldots,m_q\in\mathbb N $;
\item A sequence of weight functions
$w_{2,1}(x),w_{2,2}(x),\ldots,w_{2,q}(x):\mathbb R \to\mathbb R $.
\end{itemize}
Put $\vecn := (n_1, \ldots, n_p)$, $|\vecn| := \sum_{k=1}^p n_k$ and similarly
for $\vecm$ and $|\vecm|$. Assume $|\vecn| = |\vecm| +1$. We say that a
sequence of polynomials $A_1(x),A_2(x),\ldots,A_p(x)$ is multiple
orthogonal with respect to the above data if (i) the polynomials $A_k(x)$ have
degrees bounded by $n_k-1$:
\begin{equation}\label{degreesAk} \deg A_k\leq n_k-1,\quad
k=1,\ldots,p, \end{equation} and (ii) the function
\begin{equation}
\label{defQ} Q(x) := \sum_{k=1}^p A_k(x)w_{1,k}(x) \end{equation} satisfies the
orthogonality relations \begin{equation}
\label{orthogonalityrelationsMOP}\int_{-\infty}^{\infty} Q(x)x^j w_{2,l}(x)\ dx
= 0,
\end{equation}
for $j=0,1,\ldots,m_l-1$ and $l=1,\ldots,q$.
\end{definition}

Note that \eqref{orthogonalityrelationsMOP} states that $Q(x)$ has $m_l$
vanishing moments with respect to the weight $w_{2,l}(x)$,
$l=1,\ldots,q$.

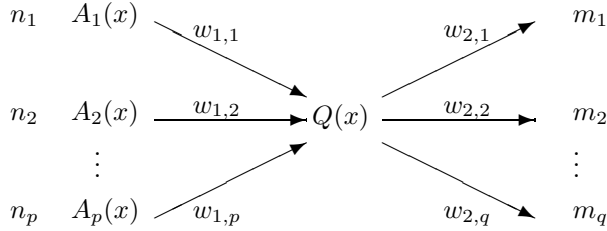
\begin{figure}[t]
\begin{center}
   \setlength{\unitlength}{1truemm}
   \begin{picture}(100,48.5)(0,2.5)

       \put(50,28){\line(-2,1){20}}
       \put(50,22){\line(-2,-1){20}}
       \put(50,25){\line(-1,0){20}}
       \put(50.7,24.5){$Q(x)$}
       \put(60,28){\line(2,1){20}}
       \put(60,22){\line(2,-1){20}}
       \put(60,25){\line(1,0){20}}

       \put(35,36){$w_{1,1}$}
       \put(35,26){$w_{1,2}$}
       \put(35,12){$w_{1,p}$}
       \put(68,36){$w_{2,1}$}
       \put(68,26){$w_{2,2}$}
       \put(68,12){$w_{2,q}$}

       \put(11,38){$n_1$}
       \put(11,25){$n_2$}
       \put(11,12){$n_p$}
       \put(19,38){$A_1(x)$}
       \put(19,25){$A_2(x)$}
       \put(22,17.5){$\vdots$}
       \put(19,12){$A_p(x)$}
       \put(85,38){$m_1$}
       \put(85,25){$m_2$}
       \put(85,17.5){$\vdots$}
       \put(85,12){$m_q$}

       \put(50,28){\thicklines\vector(2,-1){.0001}}
       \put(50,22){\thicklines\vector(2,1){.0001}}
       \put(50,25){\thicklines\vector(1,0){.0001}}
       \put(80,38){\thicklines\vector(2,1){.0001}}
       \put(80,12){\thicklines\vector(2,-1){.0001}}
       \put(80,25){\thicklines\vector(1,0){.0001}}
   \end{picture}
   \caption{The figure shows a schematic illustration of the multiple orthogonal polynomials $A_k(x)$ in Definition~\ref{defmop}.}
   \label{figdefmop}
\end{center}
\end{figure}

A schematic illustration of Definition~\ref{defmop} is shown in
Figure~\ref{figdefmop}. Let us comment on this figure. The left part of the
figure shows the polynomials $A_k(x)$ and their corresponding number of free
coefficients $n_k$, $k=1,\ldots,p$.
The middle part of the figure shows how the polynomials $A_k(x)$ should be
assembled into the function $Q(x)$ defined in \eqref{defQ}. Finally, the right
part of the figure schematically shows the orthogonality relations of $Q(x)$
with respect to the different weights $w_{2,l}(x)$, and it shows next to each
weight $w_{2,l}(x)$ also the number of vanishing moments $m_l$ of $Q(x)$ with
respect to this weight.

We will refer to the polynomials $A_k(x)$ in Definition~\ref{defmop} as
\emph{multiple orthogonal polynomials} (MOP). Note that these polynomials were
called multiple orthogonal polynomials of mixed type in \cite{DK2} and mixed
MOPs in \cite{AvMV}. We will also find it convenient to use the vectorial
notation $\vecA(x) = (A_1(x),\ldots,A_p(x))$. To stress the dependence on the
multi-indices $\vecn$, $\vecm$ we will sometimes write $\vecA(x) =
\vecA_{\vecn,\vecm}(x)$ and similarly $A_j(x) = (A_j)_{\vecn,\vecm}(x)$.

The coefficients of the multiple orthogonal polynomials in Definition
\ref{defmop} can be found from a homogeneous linear system with $|\vecn|$
unknowns (polynomial coefficients) and $|\vecm|$ equations (orthogonality
conditions). The restriction $|\vecn|=|\vecm|+1$ in Definition \ref{defmop}
guarantees that this system has a nontrivial solution. In fact, the solution
space to this linear system will be at least one-dimensional. This corresponds
to the fact that the MOP are only determined up to some multiplicative factor.
If the MOP are unique up to a multiplicative factor then the pair of
indices $\vecn,\vecm$ is called normal \cite{DK2}.

The fact that the multiple orthogonal polynomials are only determined up to
some multiplicative factor allows for different choices of normalization.

\begin{definition}\label{defnormalizationtypes} (Normalization types of MOP; cf.\  \cite{DK2})
Assume the data in Definition~\ref{defmop} and assume that $\vecn,\vecm$ is a
normal pair of indices. Then the MOP in Definition \ref{defmop} are said to
satisfy \begin{itemize}\item the \emph{normalization of type I with respect to
the $l$th index}, $l\in\{1,\ldots,q\}$, if the $m_l$th moment of $Q(x)$ with
respect to $w_{2,l}(x)$ is equal to one, i.e., if
\begin{equation}\label{type1norm}
\int_{-\infty}^{\infty} Q(x)x^{m_l}w_{2,l}(x)\ dx = 1.
\end{equation}
\item the \emph{normalization of type II with respect to the $k$th
index}, $k\in\{1,\ldots,p\}$, if the leading coefficient of $A_k(x)$ is equal
to one, i.e., if
\begin{equation}\label{type2norm}
A_k(x) = x^{n_k-1}+O(x^{n_k-2}).
\end{equation}
\end{itemize}
\end{definition}

The vectors of MOP $\vecA_{\vecn,\vecm}(x)$ corresponding to the above
normalizations will be denoted as $\vecA_{\vecn,\vecm}^{(I,l)}(x)$ and
$\vecA_{\vecn,\vecm}^{(II,k)}(x)$, respectively.

The above normalizations might not always be possible. The type $(I,l)$
normalization is not possible in those cases where the integral on the left
side of \eqref{type1norm} is equal to zero. Similarly, the type
$(II,k)$ normalization is not possible in those cases where the $k$th
polynomial $A_k(x)$ has degree strictly smaller than $n_k-1$.

As mentioned above, the MOP appear in the solution of the RH problem in
Section~\ref{subsectionintroRH}. Let us describe this is somewhat more detail.

Recall that in Definition \ref{defmop} we needed the condition
$|\vecn|=|\vecm|+1$ to ensure the existence of the multiple orthogonal
polynomials. But let us now assume that $|\vecn|=|\vecm|$. In this case, the
definition of MOP makes no sense. Indeed, the coefficients of the MOP would
then solve a homogeneous linear system with as many equations as unknowns,
which has in general only the trivial solution.

Therefore, to apply the definition of MOP in a meaningful way for a pair of
multi-indices satisfying $|\vecn|=|\vecm|$, we should first adapt the
multi-indices. There are essentially $p+q$ ways to proceed:
\begin{enumerate}
\item One can increase one of the components $n_k$, i.e., one can work with the pair of multi-indices
$\vecn+\vece_k,\vecm$, for some $k\in\{1,2,\ldots,p\}$.
\item One can decrease one of the components $m_l$, i.e., one can work with the pair of multi-indices
$\vecn,\vecm-\vece_l$, for some $l\in\{1,2,\ldots,q\}$.
\end{enumerate}
Here $\vece_k$ denotes the vector which has all its entries equal to zero,
except for the $k$th entry which equals one. The length of $\vece_k$ should be
clear from the context.

Let us discuss the MOP corresponding to each of the $p+q$ pairs of
multi-indices above. In the case of a pair of multi-indices
$\vecn+\vece_k,\vecm$, we are dealing with multiple orthogonal polynomials
where the $k$th polynomial $A_k$, $k\in\{1,2,\ldots,p\}$ has increased
degree; it will be natural to normalize the resulting MOP such that this $k$th
polynomial $A_k$ is monic, i.e., to work with a normalization of type $(II,k)$.
This leads to the vector of MOP
\begin{equation}\label{MOPtype2}
\vecA^{(II,k)}_{\vecn+\vece_k,\vecm}(x) =
((A_1)^{(II,k)}_{\vecn+\vece_k,\vecm}(x),\ldots,(A_p)^{(II,k)}_{\vecn+\vece_k,\vecm}(x)).
\end{equation}

In the case of a pair of multi-indices $\vecn,\vecm-\vece_l$, we are dealing
with multiple orthogonal polynomials where the $l$th orthogonality condition,
$l\in\{1,2,\ldots,q\}$ has a decreased number of vanishing moments; it
will be natural to normalize the resulting MOP such that this omitted moment
equals one, i.e., to work with a normalization of type $(I,l)$. This leads to
the vector of MOP
\begin{equation}\label{MOPtype1}
\vecA^{(I,l)}_{\vecn,\vecm-\vece_l}(x) =
((A_1)^{(I,l)}_{\vecn,\vecm-\vece_l}(x),\ldots,(A_p)^{(I,l)}_{\vecn,\vecm-\vece_l}(x)).
\end{equation}

Of course we are assuming here that all type I and type II normalizations in
\eqref{MOPtype2} and \eqref{MOPtype1} exist. It turns out that the existence of
each of these $p+q$ vectors of MOP is equivalent to a single condition
\cite{DKV}, to which we will loosely refer here as the solvability
condition. This condition will always be satisfied in the case of Gaussian
weight functions \eqref{gaussianw1}--\eqref{gaussianw2}.

Now we can use the above MOP to solve the RH problem in Section
\ref{subsectionintroRH}. Indeed, the $p+q$ vectors of MOP in \eqref{MOPtype2},
\eqref{MOPtype1} are row vectors of length $p$, and they can therefore be
stacked into the first $p$ columns of a matrix of size $(p+q) \times (p+q)$.
Denote with $\widetilde{Y}(z)$ such a matrix. The entries in the remaining $q$
columns of $\widetilde{Y}(z)$ are defined as Cauchy transforms of the functions
$Q^{(II,k)}$ and $Q^{(I,l)}$ defined as in \eqref{defQ}. More precisely, the
last $q$ entries in the rows $k=1, \ldots, p$ of $\widetilde{Y}(z)$ are defined
by
\begin{equation}\label{cauchytransform1}
\frac{1}{2\pi i}\int_{-\infty}^{\infty} \frac{Q^{(II,k)}(x)w_{2,l}(x)}{x-z}\
dx,
    \qquad l = 1, \ldots, q,
\end{equation}
and in the rows $p+k, \ldots, p+q$ by
\begin{equation}\label{cauchytransform2}
\frac{1}{2\pi i}\int_{-\infty}^{\infty} \frac{Q^{(I,k)}(x)w_{2,l}(x)}{x-z}\ dx,
    \qquad l = 1, \ldots, q.
\end{equation}
The fact of the matter is the following.

\begin{theorem}\label{stsolutionRH} (Solution to the Riemann-Hilbert problem; cf.\ \cite{DK2})
Let $\vecn,\vecm$ with $|\vecn|=|\vecm|$ be such that the solvability condition
holds. Then there exists a unique solution $Y(z)=Y_{\vecn,\vecm}(z)$ of the RH
problem in Section \ref{subsectionintroRH}. This solution is given by
\begin{equation*} Y(z) = D\widetilde{Y}(z),\end{equation*} where $\widetilde{Y}(z)$ is the matrix
constructed from the MOP as described above, and where $D = \diag(I_p,-2\pi i
I_q)$.
\end{theorem}

For example, in the special case where $p=q=2$, the solution $Y(z)$ of the
Riemann-Hilbert problem in Section \ref{subsectionintroRH} is given by the
$4\times 4$ matrix
\begin{equation}\label{RHmatrix4x40} Y(z) =
D\times\left(\begin{array}{cccc}
(A_1)_{\vecn+\vece_1,\vecm}^{(II,1)} & (A_2)_{\vecn+\vece_1,\vecm}^{(II,1)} & * & * \\
(A_1)_{\vecn+\vece_2,\vecm}^{(II,2)} & (A_2)_{\vecn+\vece_2,\vecm}^{(II,2)} & * & * \\
(A_1)_{\vecn,\vecm-\vece_1}^{(I,1)} & (A_2)_{\vecn,\vecm-\vece_1}^{(I,1)} & * & * \\
(A_1)_{\vecn,\vecm-\vece_2}^{(I,2)} & (A_2)_{\vecn,\vecm-\vece_2}^{(I,2)} & * & * \\
\end{array}\right),
\end{equation}
where $D := \diag(1,1,-2\pi i,-2\pi i)$, and where the entries denoted with $*$
are certain Cauchy transforms as in \eqref{cauchytransform1} and
\eqref{cauchytransform2}. The latter entries will be irrelevant in what
follows.

\subsection{Recurrence relations for multiple Hermite polynomials}
\label{subsectionintorecHerm}

In analogy with the three-term recurrence relations for classical orthogonal
polynomials on the real line, one can show that the multiple orthogonal
polynomials in Section \ref{subsectionintromop} satisfy certain $p+q+1$ term
recurrence relations. We will state these relations in the case of
multiple Hermite polynomials, i.e., when the weight functions of the MOP
are given by the Gaussians \eqref{gaussianw1}--\eqref{gaussianw2}. For
simplicity, we assume throughout that $p=q=2$, as in \eqref{RHmatrix4x40}.

Define the next term in the asymptotic expansion of $Y(z)$ in
\eqref{asymptoticconditionY0} as
\begin{equation}\label{asymptoticconditionY1}
Y(z) = \left(I+\frac{Y_1}{z}+
O\left(\frac{1}{z^2}\right)\right)\diag(z^{n_1},z^{n_2},z^{-m_1},z^{-m_2}).
\end{equation}
The entries of the matrix $Y_1=(Y_1)_{\vecn,\vecm}$ in
\eqref{asymptoticconditionY1} will be denoted by $(c_{i,j})_{i,j=1}^{4}$.

\begin{proposition}\label{proprecurrence4x4} (Recurrence relations for multiple Hermite polynomials)
Assume $p=q=2$ and let the weight functions be defined by
\eqref{gaussianw1}--\eqref{gaussianw2}. Then the multiple Hermite polynomials
in \eqref{RHmatrix4x40} satisfy the $5$-term recurrence relations
\begin{multline}\label{examplerecurrenceintro1}
    (A_1)^{(II,1)}_{\vecn+\vece_1+\vece_1,\vecm+\vece_1} =
    \left(z-(1-t)a_1-tb_1+\frac{c_{1,2}c_{2,3}}{c_{1,3}}\right)
    (A_1)^{(II,1)}_{\vecn+\vece_1,\vecm}\\
    -c_{1,2}c_{2,1}(A_1)^{(II,1)}_{\vecn+\vece_2,\vecm}
    -c_{1,3}c_{3,1}(A_1)^{(II,1)}_{\vecn,\vecm-\vece_1}
    -c_{1,4}c_{4,1}(A_1)^{(II,1)}_{\vecn,\vecm-\vece_2},
\end{multline}
\begin{multline}\label{examplerecurrenceintro2}
    (A_1)^{(II,1)}_{\vecn+\vece_1+\vece_1,\vecm+\vece_2} = \left(z-(1-t)a_1-tb_2+\frac{c_{1,2}c_{2,4}}{c_{1,4}}\right)
    (A_1)^{(II,1)}_{\vecn+\vece_1,\vecm}\\
    -c_{1,2}c_{2,1}(A_1)^{(II,1)}_{\vecn+\vece_2,\vecm}
    -c_{1,3}c_{3,1}(A_1)^{(II,1)}_{\vecn,\vecm-\vece_1}
    -c_{1,4}c_{4,1}(A_1)^{(II,1)}_{\vecn,\vecm-\vece_2},
\end{multline}
\begin{multline}\label{examplerecurrenceintro3}
    (A_2)^{(II,2)}_{\vecn+\vece_2+\vece_2,\vecm+\vece_1} =
    \left(z-(1-t)a_2-tb_1+\frac{c_{2,1}c_{1,3}}{c_{2,3}}\right)(A_2)^{(II,2)}_{\vecn+\vece_2,\vecm}\\
    -c_{2,1}c_{1,2}(A_2)^{(II,2)}_{\vecn+\vece_1,\vecm}
    -c_{2,3}c_{3,2}(A_2)^{(II,2)}_{\vecn,\vecm-\vece_1}
    -c_{2,4}c_{4,2}(A_2)^{(II,2)}_{\vecn,\vecm-\vece_2},
\end{multline}
and
\begin{multline}\label{examplerecurrenceintro4}
    (A_2)^{(II,2)}_{\vecn+\vece_2+\vece_2,\vecm+\vece_2} =
    \left(z-(1-t)a_2-tb_2+\frac{c_{2,1}c_{1,4}}{c_{2,4}}\right)
    (A_2)^{(II,2)}_{\vecn+\vece_2,\vecm} \\
    -c_{2,1}c_{1,2}(A_2)^{(II,2)}_{\vecn+\vece_1,\vecm}
    -c_{2,3}c_{3,2}(A_2)^{(II,2)}_{\vecn,\vecm-\vece_1}
    -c_{2,4}c_{4,2}(A_2)^{(II,2)}_{\vecn,\vecm-\vece_2}.
\end{multline}
Here we denote with $c_{i,j}$ the entries of the matrix
$Y_1=(Y_1)_{\vecn,\vecm}$ in \eqref{asymptoticconditionY1}.
\end{proposition}

The proof of Proposition~\ref{proprecurrence4x4} will be given in a more
general setting in Section~\ref{subsectionrecmop}. The explicit form of the
first term in the right-hand side of each of
\eqref{examplerecurrenceintro1}--\eqref{examplerecurrenceintro4} is only valid
under the assumption of Gaussian weight functions
\eqref{gaussianw1}--\eqref{gaussianw2}; the explicit form of these terms will
be established in Section~\ref{subsubsectiondiagrec}.

Note that the recurrence relations
\eqref{examplerecurrenceintro1}--\eqref{examplerecurrenceintro4} contain
several recurrence coefficients of the form $c_{i,j}c_{j,i}$ with $i<j$; it
will be convenient to collect them in the $4$ by $4$ matrix
\begin{equation}\label{hadamardproductexampleintro}
\begin{pmatrix}
0 & c_{1,2}c_{2,1} & c_{1,3}c_{3,1} & c_{1,4}c_{4,1}\\
0 & 0 & c_{2,3}c_{3,2} & c_{2,4}c_{4,2}\\
0 & 0 & 0 & c_{3,4}c_{4,3}\\
0 & 0 & 0 & 0
\end{pmatrix}.
\end{equation}
It turns out that there exist certain connections between these recurrence
coefficients.

\begin{proposition}\label{propconnectionsrecurrencecoeff4x4} (Relations between recurrence coefficients)
Assume $p=q=2$ and let the weight functions be defined by
\eqref{gaussianw1}--\eqref{gaussianw2}. Then the $2$ by $2$ submatrix
\begin{equation}\label{hadamardproductexamplesubmatrixintro}
    \begin{pmatrix}c_{1,3}c_{3,1} & c_{1,4}c_{4,1}\\
    c_{2,3}c_{3,2} & c_{2,4}c_{4,2}\end{pmatrix}
\end{equation}
of \eqref{hadamardproductexampleintro} has row sums equal to $t(1-t)
\frac{n_k}{N}$, $k=1,2$, and column sums equal to $t(1-t) \frac{m_l}{N}$,
$l=1,2$. Next, assume that $n_1=m_1$ and $n_2=m_2$. Then all the recurrence
coefficients $c_{i,j} c_{j,i}$ in \eqref{hadamardproductexampleintro} can be expressed in terms
of $c_{1,2}c_{2,1}$ and $c_{1,4}c_{4,1}$ alone, by means of the following
relations:
\begin{align}
\label{relation1rowsums} c_{2,3}c_{3,2} &= c_{1,4}c_{4,1}\\
\label{relation2rowsums} c_{1,3}c_{3,1} &= t(1-t) \frac{n_1}{N} -c_{1,4}c_{4,1} \\
\label{relation3rowsums} c_{2,4}c_{4,2} &= t(1-t) \frac{n_2}{N} -c_{1,4}c_{4,1} \\
\label{fourthrelation} t^2 (b_1-b_2)^2 c_{3,4}c_{4,3} &= (1-t)^2(a_1-a_2)^2
c_{1,2}c_{2,1}.
\end{align}
\end{proposition}
Note that \eqref{relation1rowsums}--\eqref{relation3rowsums} follow immediately
from the stated row and column sum relations for the matrix
\eqref{hadamardproductexamplesubmatrixintro}. The latter will be established
for general values of $p$ and $q$ in Section \ref{subsectionscalarproducts};
see also Section \ref{subsectionspectralcurve} for a spectral curve
interpretation of these relations. On the other hand, the equation
\eqref{fourthrelation} will be established (in a slightly more general form) in
Section \ref{subsubsectiontwodeterminants}.

\subsection{Painlev\'e II asymptotics for recurrence coefficients}
\label{subsectionintroPIIasymptotics}

Finally we are in position to formulate the Painlev\'e II asymptotics of the
recurrence coefficients in
\eqref{examplerecurrenceintro1}--\eqref{examplerecurrenceintro4} under the
double scaling regime in Section \ref{subsectionintromodelRHPII}. This is our
second main result. We first consider the off-diagonal recurrence coefficients,
i.e., the recurrence coefficients of the form $c_{i,j}c_{j,i}$ with $i<j$.

\begin{theorem}\label{stp2asymptotics} (Asymptotics of off-diagonal recurrence coefficients)
Assume the double scaling regime \eqref{deftau}--\eqref{doublescalinglimit},
and let $t\in(0,1)$ be a non-critical time, i.e., $t\neq t_{\crit}$. Define the
constants
\begin{equation}\label{defK}
    K := \frac{(p_1^*p_2^*)^{1/6}}{(\sqrt{p_1^*}+\sqrt{p_2^*})^{4/3}} > 0,
\end{equation}
and
\begin{equation}\label{defs0}
    s := -(p_1^*p_2^*)^{1/6}(\sqrt{p_1^*}+\sqrt{p_2^*})^{2/3}L\in\mathbb R ,
\end{equation}
where $L$ is defined in \eqref{doublescalinglimit}. Then we have
\begin{align}
    \label{c12asymptotics}
    c_{1,2}c_{2,1} &= - K^2 t^2(b_1-b_2)^2q^2(s)n^{-2/3}+O(n^{-1}),\\ 
    \label{c14asymptotics}
    c_{1,4}c_{4,1} &=  K^2t(1-t)(a_1-a_2)(b_1-b_2)q^2(s)
    n^{-2/3}+O(n^{-1}),
\end{align}
as $n\to\infty$, where $q(s)$ denotes the Hastings-McLeod solution to the
Painlev\'e II equation. The asymptotic behavior of the other recurrence
coefficients $c_{i,j}c_{j,i}$ with $i<j$ is then determined by
\eqref{relation1rowsums}--\eqref{fourthrelation} in Proposition
\ref{propconnectionsrecurrencecoeff4x4}.
\end{theorem}

The key point of Theorem \ref{stp2asymptotics} is that it shows that the
Painlev\'e II equation shows up in the large $n$ behavior of the recurrence
coefficients in the case of critical separation at a non-critical time.

\begin{remark}\label{remarkasymptlargesep} (The case of large separation)
Using the results in this paper, one can prove a similar result for the case of
large separation of the endpoints \eqref{largeseparationeq}. In this case, it
can be shown that there exists a constant $c>0$ such that
\begin{align}
c_{1,2}c_{2,1} &= O(e^{-cn}), \\
c_{1,4}c_{4,1} &= O(e^{-cn}),
\end{align}
as $n\to\infty$.
\end{remark}

\begin{remark}\label{remarkasymptsmallsep} (The case of small separation)
Performing a similar analysis of the results in \cite{DKV}, one can show that
in case of small separation of the endpoints \eqref{smallseparationeq}, and
assuming the additional hypothesis $p_1=p_2 = 1/2$, the following expansions
hold:
\begin{align}
\label{typeaEvi1} c_{1,2}c_{2,1} &= -\frac{t^2}{16(a_1-a_2)^2}(4-(a_1-a_2)^2(b_1-b_2)^2)+O(n^{-1}),\\
\label{typeaEvi2} c_{1,4}c_{4,1} &=
\frac{t(1-t)}{8}(2-(a_1-a_2)(b_1-b_2))+O(n^{-1}),
\end{align}
as $n\to\infty$.
\end{remark}

We have a similar theorem for the diagonal recurrence coefficients in
\eqref{examplerecurrenceintro1}--\eqref{examplerecurrenceintro4}, i.e., for the
first terms in the right-hand side of each of these equations.

\begin{theorem}\label{stp2asymptoticsB} (Asymptotics of diagonal recurrence coefficients)
Under the same assumptions as in Theorem \ref{stp2asymptotics} we have that
\begin{align}
    \label{c12c23c13asymptotics}
    \frac{c_{1,2}c_{2,3}}{c_{1,3}} &= -K^2 q^2(s) t (b_1-b_2) \sqrt{\frac{(a_1-a_2)(b_1-b_2)}{p_1^*}}n^{-2/3} +O(n^{-1}), \\
    \label{c12c24c14asymptotics}
    \frac{c_{1,2}c_{2,4}}{c_{1,4}} &= -t  \sqrt{\frac{p_2^*(b_1-b_2)}{a_1-a_2}} +O(n^{-1/3}), \\
    \label{c21c13c23asymptotics}
    \frac{c_{2,1}c_{1,3}}{c_{2,3}} &=  t  \sqrt{\frac{p_1^*(b_1-b_2)}{a_1-a_2}} +O(n^{-1/3}), \\
    \label{c21c14c24asymptotics}
    \frac{c_{2,1}c_{1,4}}{c_{2,4}} &=
     K^2 q^2(s) t (b_1-b_2) \sqrt{\frac{(a_1-a_2)(b_1-b_2)}{p_2^*}}n^{-2/3} +O(n^{-1}),
\end{align}
as $n \to \infty$.
\end{theorem}

\subsection{Phase diagram}

\begin{figure}[tbp]
\centering \begin{overpic}[scale=0.5]%
{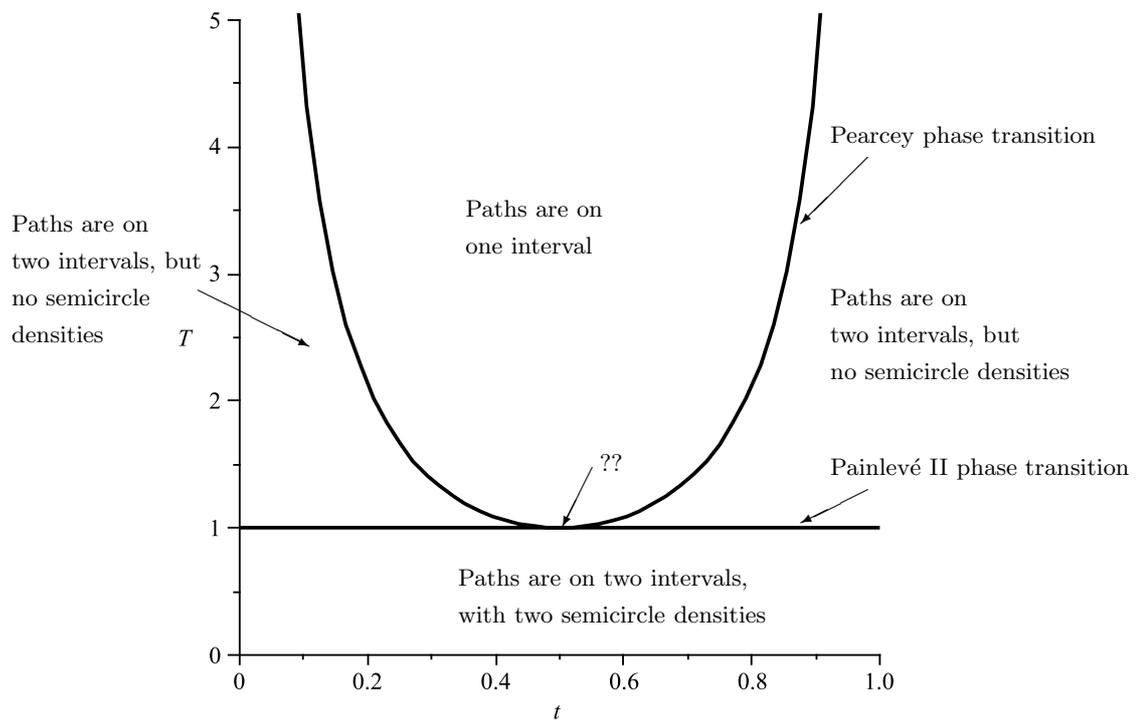}%
        \put(90,80){\small{Pearcey phase transition}}
        \put(95,78){\vector(-1,-1){9}}
        \put(90,35){\small{Painlev\'e II phase transition}}
        \put(95,33){\vector(-2,-1){9}}
        \put(59,35.5){\small{??}}
        \put(58,36){\vector(-1,-2){4}}
      \put(41,70){\small{Paths are on}}
      \put(41,65){\small{one interval}}
      \put(90,58){\small{Paths are on}}
      \put(90,53){\small{two intervals, but}}
      \put(90,48){\small{no semicircle densities}}
      \put(-20,68){\small{Paths are on}}
      \put(-20,63){\small{two intervals, but}}
      \put(-20,58){\small{no semicircle}}
      \put(-20,53){\small{densities}}
      \put(5,60){\vector(2,-1){15}}
      \put(40,20){\small{Paths are on two intervals,}}
      \put(40,15){\small{with two semicircle densities}}
\end{overpic}
\caption{The figure shows the phase diagram for non-intersecting Brownian
motions with two starting and ending points. We have chosen here
$p_1^*=p_2^*=1/2$, $a_1=-a_2=b_1=-b_2=1/\sqrt{2}$, hence $t_{\crit}=1/2$. The
boldface curves denote the places where a phase transition takes place. The
topmost curve has equation $T = (t^2-t+1/2)/(t(1-t))$, cf.\
\eqref{vglboundarycurve}.} \label{figphasediagram}
\end{figure}

The main results of this paper and their relation to the other results known in
the literature can be nicely summarized by means of a phase diagram. See Figure
\ref{figphasediagram}.

Let us comment on Figure \ref{figphasediagram}. Assume that the endpoints
$a_j$, $b_j$, $j=1,2$ are fixed and satisfy the critical separation
\eqref{deftau}. The horizontal axis in the figure denotes the time $t\in(0,1)$
and the vertical axis denotes the temperature $T>0$. The diagram is divided
into different regions according to the behavior of the limiting distribution
for $n\to\infty$ of the non-intersecting Brownian motions at time $t$ and
temperature $T$. The region where $T<1$ corresponds to the case of large
separation; according to Theorem~\ref{st2ellipses} the limiting distribution is
given here by two semicircles on the two disjoint intervals
$[\alpha_1^*(t),\beta_1^*(t)]$ and $[\alpha_2^*(t),\beta_2^*(t)]$. At
temperature $T=1$, these two intervals meet each other at a certain critical
time $t_{\crit}$. When $T$ further increases, the intersection region between
the two groups of Brownian motions starts to grow; this is indicated by the
boldface curve in the middle of the picture. In the region below this curve but
above $T=1$, the limiting distribution at time $t$ is still supported on two
disjoint intervals, but now with distribution described in terms of a certain
algebraic curve of degree 4 \cite{DKV}, rather than semicircle laws. In the
region above the curve, the limiting distribution is on one interval
\cite{DKV}.

In case where $p_1^*=p_2^*=1/2$, one can find an explicit description of the
boundary curve in Figure \ref{figphasediagram} from the results in
\cite{DKV}; it turns out that this curve is given by the equation
\begin{equation*}\label{vglboundarycurvebis}
T =
\left(\frac{a_1-a_2}{2}\right)^2\frac{1-t}{t}+\left(\frac{b_1-b_2}{2}\right)^2\frac{t}{1-t},
\end{equation*}
or equivalently
\begin{equation}\label{vglboundarycurve}
T = \frac{(a_1-a_2)^2(1-t)^2+(b_1-b_2)^2t^2}{4t(1-t)}.
\end{equation}
The curve plotted in Figure \ref{figphasediagram} corresponds to the choice of
endpoints $a_1=-a_2=b_1=-b_2=1/\sqrt{2}$ and hence $t_{\crit}=1/2$.

Figure \ref{figphasediagram} also displays the phase transitions between the
different regions. On the horizontal line $T=1$ the phase transition is
described in terms of the Painlev\'e II equation as shown in this paper. On the
curve \eqref{vglboundarycurve} one expects a description in terms of the
Pearcey kernels; although for the case of two starting and two ending points
this has not been strictly proven. For the case of non-intersecting Brownian
motion with one starting and two ending points, the Pearcey kernels were
obtained in \cite{BH1,BH2} (in the equivalent setting of Gaussian random
matrices with external source), see also \cite{AvM1,BK3,OR,TW3}. Finally, at
the place where the two curves in Figure \ref{figphasediagram} meet, one
expects a phase transition in terms of an as yet unknown family of kernels.
This is indicated by the question mark in the figure.

\subsection{Outline of the paper}

The remainder of this paper is organized as follows. In Sections
\ref{sectionsteepestdescentlargesep} and \ref{sectionsteepestdescentcritical}
we apply the Deift-Zhou steepest descent analysis to the Riemann-Hilbert
problem for multiple Hermite polynomials. We perform this analysis for the
case of large separation in Section~\ref{sectionsteepestdescentlargesep} and
for the case of critical separation at a non-critical time in
Section~\ref{sectionsteepestdescentcritical}, leading to the proof of
Theorem~\ref{st2ellipses}. In the critical case we are also led to a local
parametrix for the RH problem in terms of the Painlev\'e II equation. Next, we
investigate the large $n$ asymptotics of the recurrence coefficients of the
multiple Hermite polynomials under the double scaling regime. This is done in
Section~\ref{sectionPIIasymptotics}, leading to the proof of
Theorems~\ref{stp2asymptotics} and \ref{stp2asymptoticsB}. Finally, the
 Propositions~\ref{proprecurrence4x4} and
\ref{propconnectionsrecurrencecoeff4x4} are established in
Section~\ref{sectionlaxpair}.

\section{Steepest descent analysis in the case of large separation}
\label{sectionsteepestdescentlargesep}

In this section we analyse the non-intersecting Brownian motions in case of
large separation between the endpoints \eqref{largeseparationeq}. Using the
Deift/Zhou steepest descent analysis of the Riemann-Hilbert problem, we show
that the interaction between the two groups of Brownian particles decays
exponentially with $n$, thereby establishing Theorem~\ref{st2ellipses}.

\subsection{Starting RH problem}
Our starting point is the RH problem
\eqref{defjumpmatrix0}--\eqref{asymptoticconditionY0} with $p=q=2$ and in
addition $n_1 = m_1$ and $n_2 = m_2$. As in \eqref{defp1p2} we write $p_j =
n_j/n$. Without loss of generality we take $N =n$ (i.e., $T = 1$). For $T = 1$
the case of large separation corresponds to $(a_1-a_2) (b_1-b_2) >
(\sqrt{p_1^*} + \sqrt{p_2^*})^2$. Since $p_j \to p_j^*$ as $n \to \infty$, we
already assume that $n$ is so large that
\[ (a_1-a_2) (b_1-b_2) > (\sqrt{p_1} + \sqrt{p_2})^2. \]

Thus $Y$ satisfies the following RH problem.
\begin{itemize}
\item[(1)] $Y(z)$ is analytic in $\mathbb C \setminus \mathbb R $;
 \item[(2)] For $x \in \mathbb R$ we have
 \[ Y_+(x) = Y_-(x) \begin{pmatrix} I_2 & W(x) \\ 0 & I_2 \end{pmatrix} \]
 with the rank-one block $W$ given by
\begin{equation}\label{defWblockbis}
W = \begin{pmatrix} w_{1,1}\\ w_{1,2}\end{pmatrix}
\begin{pmatrix} w_{2,1}& w_{2,2}\end{pmatrix} =
\begin{pmatrix} w_{1,1}w_{2,1} & w_{1,1}w_{2,2}\\
w_{1,2}w_{2,1} & w_{1,2}w_{2,2}
\end{pmatrix}.
\end{equation}
  \item[(3)] As $z\to\infty$, we have that
\begin{equation*}Y(z) =
    (I_{4}+O(1/z))\diag(z^{n_1},z^{n_2},z^{-n_1},z^{-n_2}).
\end{equation*}
\end{itemize}
The entries of the rank-one block $W = W(x)$ in \eqref{defWblockbis} can be
written explicitly as
\begin{align}
    w_{1,k}(x)w_{2,l}(x)
    \label{productws} &= e^{-\frac{n}{2t(1-t)}(x^2-2((1-t)a_k+tb_l)x)},
\end{align}
for $k,l\in\{1,2\}$. Recall that $N = n$.

It will be convenient to write the diagonal entries of \eqref{defWblockbis} as
\begin{align}
\label{decoupledweight1} w_{1,1}(x)w_{2,1}(x) &= e^{-np_1V_1(x)} = e^{-n_1 V_1(x)},\\
\label{decoupledweight2} w_{1,2}(x)w_{2,2}(x) &= e^{-np_2V_2(x)} = e^{-n_2 V_2(x)},
\end{align}
with $p_j = n_j/n$, for $j=1,2$, and
\begin{align}
\label{defV1} V_1(x) & := \frac{1}{2p_1t(1-t)}(x^2-2\left((1-t)a_1+tb_1\right)x),  \\
\label{defV2} V_2(x) & := \frac{1}{2p_2t(1-t)}(x^2-2\left((1-t)a_2+tb_2\right)x).
\end{align}

Our goal is to show that in the large $n$ limit the $4\times 4$ matrix valued RH problem
for  $Y(z)$ essentially decouples into two smaller
$2\times 2$ problems with weight functions
\eqref{decoupledweight1} and \eqref{decoupledweight2}. To show
that this decoupling indeed occurs, we need to show that in some
sense, the off-diagonal entries in \eqref{defWblockbis} can be
neglected with respect to the diagonal entries of
\eqref{defWblockbis}. More precisely, one expects that
\begin{itemize}
\item around the interval $(\alpha_1^*,\beta_1^*)$, the $(1,1)$ entry of
\eqref{defWblockbis} is dominant (as $n\to\infty$) with respect to the other
three entries.
\item around the interval $(\alpha_2^*,\beta_2^*)$, the $(2,2)$ entry of
\eqref{defWblockbis} is dominant with respect to the other three entries.
\end{itemize}
Here $\alpha_j^*$, $\beta_j^*$ for $j=1,2$ are as in \eqref{defalpha}--\eqref{defbeta}.

Remarkably, these expectations are not confirmed by a straightforward
steepest descent analysis in which the first  transformation of the RH problem
is based on the two semicircle densities \eqref{wignerdensities} and the
corresponding $g$-functions. This approach turns out to be successful only for
$t$ near the critical time $t_{\crit}$ defined in \eqref{deftcrit0}. When $t$
is sufficiently close to $0$ however, one runs into difficulties since then
\begin{itemize}
\item the $(1,2)$ entry of
\eqref{defWblockbis} blows up (i.e., becomes exponentially large when
$n\to\infty$) somewhere in the interval $(\alpha_1^*,\beta_1^*)$, and
\item the $(2,1)$ entry of
\eqref{defWblockbis}  blows up somewhere in the interval $(\alpha_2^*,\beta_2^*)$.
\end{itemize}
Similar problems occur when $t$ is close to $1$, but then with the roles of the $(1,2)$
and $(2,1)$ entries of $W$  reversed.

In order to prevent the blow-up of undesired entries, we make a first preliminary
transformation that is described in the next subsection.
The transformation is different for the two cases $0 < t \leq t_{\crit}$
and $t_{crit} \leq t < 1$. For definiteness we assume from now on
\begin{equation}\label{hypothesist}
    0 < t \leq t_{\crit} = \frac{a_1-a_2}{(a_1-a_2) + (b_1-b_2)}.
\end{equation}
The case where $t_{\crit} \leq t< 1$ is similar and the corresponding
modifications will be briefly commented on later.

\subsection{First transformation: Gaussian elimination in the jump matrix}
\label{subsectiontransfo1large}

The first transformation of the RH problem is a Gaussian elimination step for the
jump matrix, serving to annihilate some of the undesired entries in
\eqref{defWblockbis}. This elimination step will be at the price of introducing
new jump matrices on certain contours $\Gamma_1$, $\Gamma_2$ in the complex plane.
This transformation is similar to the so-called global opening of the
lens discussed in \cite{ABK,AKV,DKV,LW}.

Similar to \eqref{defalpha}--\eqref{defbeta} we define
\begin{align} \label{defalphavarying}
    \alpha_j & = \alpha_j(t) =  (1-t)a_j+tb_j-\sqrt{4p_j t(1-t)},\\
    \label{defbetavarying}
    \beta_j & =  \beta_j(t) = (1-t)a_j+tb_j+\sqrt{4p_j t(1-t)},
\end{align}
for $j=1,2$. Since $p_j$ is varying with $n$, the quantities
\eqref{defalphavarying}--\eqref{defbetavarying} are also varying with $n$. For
$n \to \infty$ they tend to $\alpha_j^*$ and $\beta_j^*$ given by
\eqref{defalpha}--\eqref{defbeta} with $T=1$. Define the corresponding
semicircle laws
\begin{align}
    \label{wignerdensitiesvarying}
    \frac{1}{2\pi t(1-t)}\sqrt{(\beta_j-x)(x-\alpha_j)},\quad
    x\in[\alpha_j,\beta_j],
        \quad j=1,2,
\end{align}
which are also (slightly) varying with $n$.

\begin{figure}[t]
\begin{center}
   \setlength{\unitlength}{1truemm}
   \begin{picture}(100,70)(-5,2)
       \qbezier(32,34)(35,40)(32,46)
       \put(32,34){\line(-1,-2){12}}
       \put(32,46){\line(-1,2){12}}
       \put(18,62){$\Gamma_2$}
       \put(33.7,40){\thicklines\circle*{1}}
       \put(34,36.6){$x_2$}
       \qbezier(52,34)(49,40)(52,46)
       \put(52,34){\line(1,-2){12}}
       \put(52,46){\line(1,2){12}}
       \put(63,62){$\Gamma_1$}
       \put(50.2,40){\thicklines\circle*{1}}
       \put(47,36.6){$x_1$}

       \put(56,54){\thicklines\vector(1,2){1}}
       \put(28,54){\thicklines\vector(-1,2){1}}
       \put(28,26){\thicklines\vector(1,2){1}}
       \put(56,26){\thicklines\vector(-1,2){1}}

       \put(-15,40){\line(1,0){110}}
       \put(93,36.6){$\mathbb R $}
       \put(95,40){\thicklines\vector(1,0){.0001}}
       \put(9,36.6){$\alpha_2$}
       \put(27,36.6){$\beta_2$}
       \put(39,36.6){$x_0$}
       \put(57,36.6){$\alpha_1$}
       \put(84,36.6){$\beta_1$}
       \put(10,40){\thicklines\circle*{1}}
       \put(28,40){\thicklines\circle*{1}}
       \put(40,40){\thicklines\circle*{1}}
       \put(58,40){\thicklines\circle*{1}}
       \put(85,40){\thicklines\circle*{1}}

       \put(34.5,44.2){$\small{\begin{pmatrix}I_2 &W \\ 0 & I_2 \end{pmatrix}}$}
       \put(62.5,14){$\small{\begin{pmatrix}1&0&0&0\\ 0&1&0&0  \\
       0&0&1&\frac{w_{2,2}}{w_{2,1}}\\ 0&0&0&1
       \end{pmatrix}}$}
       \put(-8.5,14){$\small{\begin{pmatrix}1&0&0&0\\ 0&1&0&0  \\
       0&0&1&0\\ 0&0&-\frac{w_{2,1}}{w_{2,2}}&1
       \end{pmatrix}}$}
       \put(-15,47){$\small{\begin{pmatrix}1&0&0&w_{1,1}w_{2,2}\\ 0&1&0&w_{1,2}w_{2,2}  \\
       0&0&1&0\\ 0&0&0&1
       \end{pmatrix}}$}
       \put(65,47){$\small{\begin{pmatrix}1&0&w_{1,1}w_{2,1}&0\\ 0&1&w_{1,2}w_{2,1}&0  \\
       0&0&1&0\\ 0&0&0&1
       \end{pmatrix}}$}
   \end{picture}
   \caption{The figure shows the curves $\Gamma_1$, $\Gamma_2$,
   the intervals $(-\infty,x_2)$, $(x_2,x_1)$, $(x_1,\infty)$,
   and the jump matrices on each of these five curves in the RH problem for
   $X=X(z)$.}
   \label{figgamma12}
\end{center}
\end{figure}
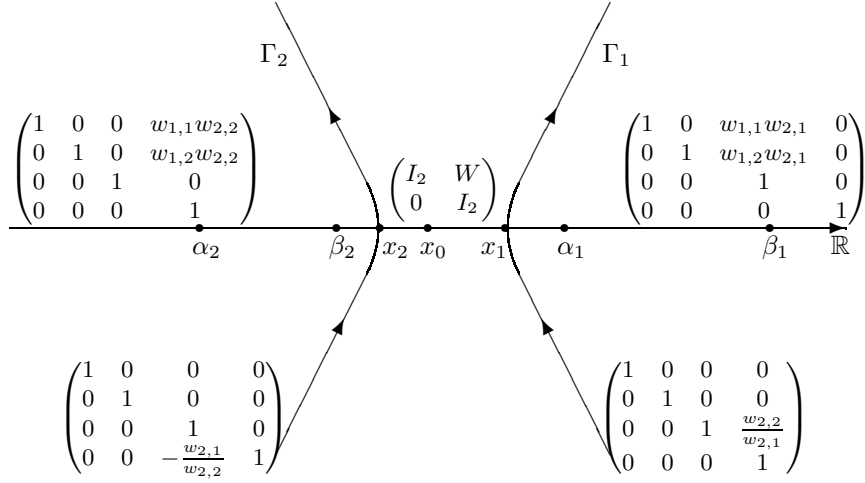

We take a reference point $x_0 \in (\beta_2, \alpha_1)$ and we choose
unbounded contours $\Gamma_1$, $\Gamma_2$ in the complex plane, crossing the
real axis in points $x_1$, $x_2$ so that
\[ \beta_2 < x_2 < x_0 < x_1 < \alpha_1. \]
The contour $\Gamma_1$ stretches out to infinity in the right half-plane (i.e.,
$\Re z \to \infty$ as $z \to \infty$ on $\Gamma_1$), and $\Gamma_2$ stretches
out to infinity in the left half-plane. The precise way to choose $x_0$ and the
curves $\Gamma_1$, $\Gamma_2$ will be described later. We orient these curves
in the upward direction as in Figure \ref{figgamma12}. We then define a new
$4\times 4$ matrix-valued function $X = X(z)$ by
\begin{align}
\label{defX1} X&= Y\begin{pmatrix}1&0&0&0\\ 0&1&0&0  \\
0&0&1&-\frac{w_{2,2}}{w_{2,1}}\\
0&0&0&1\end{pmatrix},\quad\textrm{to the right of }\Gamma_1,\\
\label{defX2} X&= Y\begin{pmatrix}1&0&0&0\\ 0&1&0&0  \\
0&0&1&0\\ 0&0&-\frac{w_{2,1}}{w_{2,2}}&1
\end{pmatrix},\quad\textrm{to the left of }\Gamma_2,\\
\label{defX3} X&= Y, \qquad \qquad\textrm{between }\Gamma_1 \textrm{ and }\Gamma_2.
\end{align}
The matrix function $X$ satisfies a new RH problem, with jumps on the contour
$\mathbb R \cup\Gamma_1\cup\Gamma_2$. The jump matrices are different on each of the
five pieces $(-\infty,x_2)$, $(x_2,x_1)$, $(x_1,\infty)$, $\Gamma_1$ and
$\Gamma_2$. They are shown in Figure \ref{figgamma12}.

Thus $X$ satisfies the following RH problem.
\begin{itemize}
\item[(1)] $X(z)$ is analytic in $\mathbb C \setminus (\mathbb R  \cup \Gamma_1 \cup \Gamma_2)$;
 \item[(2)] On $\mathbb R  \cup \Gamma_1 \cup \Gamma_2$ we have that $X_+ = X_- J_X$
 with jump matrices $J_X$ as shown  in Figure \ref{figgamma12}.
  \item[(3)] As $z\to\infty$, we have that
\begin{equation}\label{asymptoticconditionX} X(z) =
    (I_{4}+O(1/z))\diag(z^{n_1},z^{n_2},z^{-n_1},z^{-n_2}).
\end{equation}
\end{itemize}

For the asymptotic condition \eqref{asymptoticconditionX} we note that
\begin{align} \label{quotientw2}
    \frac{w_{2,2}(z)}{w_{2,1}(z)} &=
    \exp \left(-\frac{n}{1-t}(b_1-b_2)z\right).
\end{align}
Since $b_1 > b_2$ we have that \eqref{quotientw2} tends to $0$ exponentially fast
as $\Re z \to \infty$, and so in particular as $z \to \infty$ to the right of $\Gamma_1$.
Similarly, the inverse of \eqref{quotientw2} tends to $0$ exponentially fast
to the left of $\Gamma_2$.
Thus the transformation \eqref{defX1}, \eqref{defX2}, \eqref{defX3} leads to
the same asymptotic behavior for $X$ as we had for $Y$.

Note also that the jump matrices on $\Gamma_1$ and $\Gamma_2$ tend to
the identity matrix as $z \to \infty$.

What we have gained is that in the jump matrices on the
intervals $(-\infty,x_2)$ and $(x_1,\infty)$, the first and second columns of
the top right $2\times 2$ block \eqref{defWblockbis} of the jump matrix have
been eliminated, respectively.

\begin{remark}\label{remarkeliminationothercase} The description of the Gaussian
elimination step above has been done under the assumption \eqref{hypothesist},
i.e., $0<t\leq t_{\crit}$. The case where $t_{\crit} <t<1$ is different since it requires
other entries of the jump matrix to be eliminated in the different regions of
the complex plane. In order to do so one would define $X$ differently as
(compare with \eqref{defX1}, \eqref{defX2})
\begin{align*}
X&= Y\begin{pmatrix}1&0&0&0\\ \frac{w_{1,2}}{w_{1,1}}&1&0&0  \\
0&0&1&0\\
0&0&0&1\end{pmatrix},\quad\textrm{to the right of }\Gamma_1,\\
X&= Y\begin{pmatrix}1&\frac{w_{1,1}}{w_{1,2}}&0&0\\ 0&1&0&0  \\
0&0&1&0\\ 0&0&0&1
\end{pmatrix},\quad\textrm{to the left of }\Gamma_2,\\
\nonumber X&= Y,\qquad \qquad\textrm{in the region between }\Gamma_1 \textrm{ and } \Gamma_2.
\end{align*}
The further steps in the steepest descent analysis will then be similar to the
case where $0<t\leq t_{\crit}$ and we will not discuss this any further.

Alternatively, the results for $t_{\crit} <t<1$ could be
reduced to those for $0<t\leq t_{\crit}$ by virtue of the involution symmetry
to be discussed in Section \ref{subsectioninvolution}; a further note about
this will be given in Remark~\ref{remarkuseinvolutionsymmetry}.
\end{remark}

\subsection{Second transformation: $g$-functions}
\label{subsectiontransfo2large}

The second transformation is the normalization of the RH problem
at infinity. To this end we use the so-called $g$-functions.

As said before, in the limit $n\to\infty$ we expect the Brownian particles to
be distributed on two separate intervals $[\alpha_1^*,\beta_1^*]$ and
$[\alpha_2^*,\beta_2^*]$ (recall $t$ is fixed), with limiting densities given by the
two Wigner semicircle laws \eqref{wignerdensities}.

For finite $n$, we have defined $\alpha_j$ and $\beta_j$ by \eqref{defalphavarying}--\eqref{defbetavarying}
and the densities \eqref{wignerdensitiesvarying} that are varying with $n$.
We use the $g$-functions corresponding to these varying semicircle densities.
We define
\begin{align} \label{defg1g2}
g_1(z) &= \frac{1}{2\pi t(1-t)}\int_{\alpha_1}^{\beta_1}
\log(z-x)\sqrt{(\beta_1-x)(x-\alpha_1)}\ dx,\\
g_2(z) &= \frac{1}{2\pi t(1-t)}\int_{\alpha_2}^{\beta_2}
\log(z-x)\sqrt{(\beta_2-x)(x-\alpha_2)}\ dx,
\end{align}
where the logarithms are chosen with a branch cut on the positive and on the
negative real axis, respectively. Hence $g_1(z)$ is an analytic function on
$\mathbb C \setminus [\alpha_1,\infty)$ while $g_2(z)$ is analytic on $\mathbb C \setminus
(-\infty,\beta_2]$.
Note that
\[ g_j(z) = p_j \log z + O(1/z), \qquad z \to \infty. \]

It follows by contour integration that
\begin{align}
\label{diffg1} g_1'(z) &= \frac{1}{2 t(1-t)}\left[z-\frac{\alpha_1+\beta_1}{2}-((z-\alpha_1)(z-\beta_1))^{1/2}\right],\\
\label{diffg2} g_2'(z) &=
\frac{1}{2 t(1-t)}\left[z-\frac{\alpha_2+\beta_2}{2}-((z-\alpha_2)(z-\beta_2))^{1/2}\right],
\end{align}
for all $z$ in the domain of analyticity of $g_1$ and $g_2$. Here the branches of
the square roots are taken which are defined in
$\mathbb C \setminus[\alpha_1,\beta_1]$ and $\mathbb C \setminus[\alpha_2,\beta_2]$,
respectively, and which behave as $z$ when $z \to\infty$.

On the interval $[\alpha_1,\beta_1]$, the square root in
\eqref{diffg1} has two distinct boundary values, one from the upper and one from the lower half plane of
$\mathbb C $. These boundary values are complex conjugate and purely imaginary.
Similar statements hold for the behavior of
the square root in \eqref{diffg2} on the interval $[\alpha_2,\beta_2]$.

It follows from the above observations and from the definitions \eqref{defV1},
\eqref{defV2}, \eqref{defalpha}--\eqref{defbeta} that
\begin{align}
    g_{1,+}'(x)+g_{1,-}'(x) & = \frac{1}{t(1-t)} \nonumber
    \left(x-\frac{\alpha_1+\beta_1}{2} \right)  \\
    & = p_1 V_1'(x),& \quad x\in[\alpha_1,\beta_1], \label{beforeEulerlagrange1} \\
 g_{2,+}'(x)+g_{2,-}'(x)  &= \nonumber
    \frac{1}{t(1-t)}\left(x-\frac{\alpha_2+\beta_2}{2}\right) \\
    & = p_2 V_2'(x),& \quad
x\in[\alpha_2,\beta_2]. \label{beforeEulerlagrange2}
\end{align}
After integration, we see that there exist constants $l_1$, $l_2\in\mathbb R $ such that
\begin{align}
\label{Eulerlagrange1} g_{1,+}+g_{1,-} - p_1 V_1 = l_1 &\quad \textrm{on }[\alpha_1,\beta_1],\\
\label{Eulerlagrange2} g_{2,+}+g_{2,-} - p_2 V_2 = l_2 &\quad \textrm{on
}[\alpha_2,\beta_2].
\end{align}
Moreover, from the signs of the square roots in \eqref{diffg1} and
\eqref{diffg2} we have
\begin{align}
\label{Eulerlagrange3} g_{1,+}+g_{1,-} - p_1 V_1 < l_1 & \quad \textrm{on }\mathbb R \setminus[\alpha_1,\beta_1],\\
\label{Eulerlagrange4} g_{2,+}+g_{2,-} - p_2 V_2 < l_2 & \quad \textrm{on
}\mathbb R \setminus[\alpha_2,\beta_2].
\end{align}

The above equations \eqref{Eulerlagrange1}--\eqref{Eulerlagrange4} are the
Euler-Lagrange variational conditions for the equilibrium measure  under an
external field. Indeed, the Wigner semicircle laws
\eqref{wignerdensitiesvarying} (after normalization by a factor $1/p_j$) are
the equilibrium measures in the presence of the quadratic external fields
$V_1$, $V_2$, respectively, see \cite{Dei,SaffTotik}. The fact that $V_1$ in
\eqref{defV1} has a factor $p_1$ in its denominator can be interpreted by
noting that for $p_1\to 0$, the external field $V_1$ gets stronger and stronger
and hence the (potential theoretic) electrostatic particles are pushed together
onto a narrower and narrower interval, which is confirmed by
\eqref{defalphavarying}--\eqref{defbetavarying}.

Now we use the $g$-functions to normalize the RH problem at infinity.
We define a new $4\times 4$ matrix-valued function $T = T(z)$ by
\begin{equation}\label{defT}
T=  L^{-n}X G^nL^n ,
\end{equation}
where
\begin{align}
\label{defG} G &= \diag\left(e^{-g_1},e^{-g_2},e^{g_1},e^{g_2}\right),\\
\label{defLmx} L &= \diag\left(e^{l_1/2},e^{l_2/2+\kappa},e^{-l_1/2},e^{-l_2/2
+\kappa}\right).
\end{align}
Here $g_1$, $g_2$ are the $g$-functions, $l_1$, $l_2$ are the variational
constants defined in \eqref{Eulerlagrange1} and \eqref{Eulerlagrange2}, and
$\kappa\in\mathbb C $ is a constant to be determined.

The RH problem for $T=T(z)$ is normalized at infinity in the sense
that
\begin{equation}\label{asymptoticsT}T(z) = I_4 +O(1/z)\end{equation} as $z\to\infty$.
This follows from the fact that $g_j(z) = p_j \log z + O(1/z)$ as $z\to\infty$ for $j=1,2$,
so that
\[ e^{n g_j(z)} = z^{np_j}(1 + O(1/z)) = z^{n_j}(1+O(1/z)) \]
as $z \to \infty$. Hence the factor $G^n$ in \eqref{defT} cancels the powers of
$z$ appearing in \eqref{asymptoticconditionX}. Also note that the similarity
relation with the matrix $L$ in \eqref{defT} does not change the asymptotics
\eqref{asymptoticsT}.

\begin{figure}[t]
\begin{center}
   \setlength{\unitlength}{1truemm}
   \begin{picture}(100,70)(-5,2)
       \qbezier(32,34)(35,40)(32,46)
       \put(32,34){\line(-1,-2){12}}
       \put(32,46){\line(-1,2){12}}
       \put(19,62){$\Gamma_2$}
       \put(33.7,40){\thicklines\circle*{1}}
       \put(34,36.6){$x_2$}
       \qbezier(52,34)(49,40)(52,46)
       \put(52,34){\line(1,-2){12}}
       \put(52,46){\line(1,2){12}}
       \put(62,62){$\Gamma_1$}
       \put(50.2,40){\thicklines\circle*{1}}
       \put(47,36.6){$x_1$}

       \put(56,54){\thicklines\vector(1,2){1}}
       \put(28,54){\thicklines\vector(-1,2){1}}
       \put(28,26){\thicklines\vector(1,2){1}}
       \put(56,26){\thicklines\vector(-1,2){1}}

       \put(-15,40){\line(1,0){110}}
       \put(93,36.6){$\mathbb R $}
       \put(95,40){\thicklines\vector(1,0){.0001}}
       \put(9,36.6){$\alpha_2$}
       \put(27,36.6){$\beta_2$}
       \put(39,36.6){$x_0$}
       \put(57,36.6){$\alpha_1$}
       \put(84,36.6){$\beta_1$}
       \put(10,40){\thicklines\circle*{1}}
       \put(28,40){\thicklines\circle*{1}}
       \put(40,40){\thicklines\circle*{1}}
       \put(58,40){\thicklines\circle*{1}}
       \put(85,40){\thicklines\circle*{1}}

       \put(34.5,44.2){$\small{\begin{pmatrix}I_2 &\widetilde{W}\\ 0&I_2 \end{pmatrix}}$}
       \put(62.5,14){$\small{\begin{pmatrix}1&0&0&0\\ 0&1&0&0  \\
       0&0&1&e^{n(\lambda_4-\lambda_3)}\\ 0&0&0&1
       \end{pmatrix}}$}
       \put(-16,14){$\small{\begin{pmatrix}1&0&0&0\\ 0&1&0&0  \\
       0&0&1&0\\ 0&0&-e^{n(\lambda_3-\lambda_4)}&1
       \end{pmatrix}}$}
       \put(-32.5,47){$\small{\begin{pmatrix}1&0&0&e^{n(\lambda_{4,+}-\lambda_1\phantom{{}_{-}})}\\ 0&e^{n(\lambda_{2,+}-\lambda_{2,-})}&0&e^{n(\lambda_{4,+}-\lambda_{2,-})}  \\
       0&0&1&0\\ 0&0&0&e^{n(\lambda_{4,+}-\lambda_{4,-})}
       \end{pmatrix}}$}
       \put(59,47){$\small{\begin{pmatrix}e^{n(\lambda_{1,+}-\lambda_{1,-})}&0&e^{n(\lambda_{3,+}-\lambda_{1,-})}&0\\ 0&1&e^{n(\lambda_{3,+}-\lambda_2\phantom{{}_{-}})}&0  \\
       0&0&e^{n(\lambda_{3,+}-\lambda_{3,-})}&0\\ 0&0&0&1
       \end{pmatrix}}$}
     \end{picture}
   \caption{The figure shows the jump matrices in the RH problem
   for $T=T(z)$.
   Here $\widetilde{W}$ is defined in \eqref{deftildeW}, and $\lambda_1(z),\ldots,\lambda_4(z)$
   are defined in \eqref{deflambda14}.}
   \label{figjumpslambdas}
\end{center}
\end{figure}
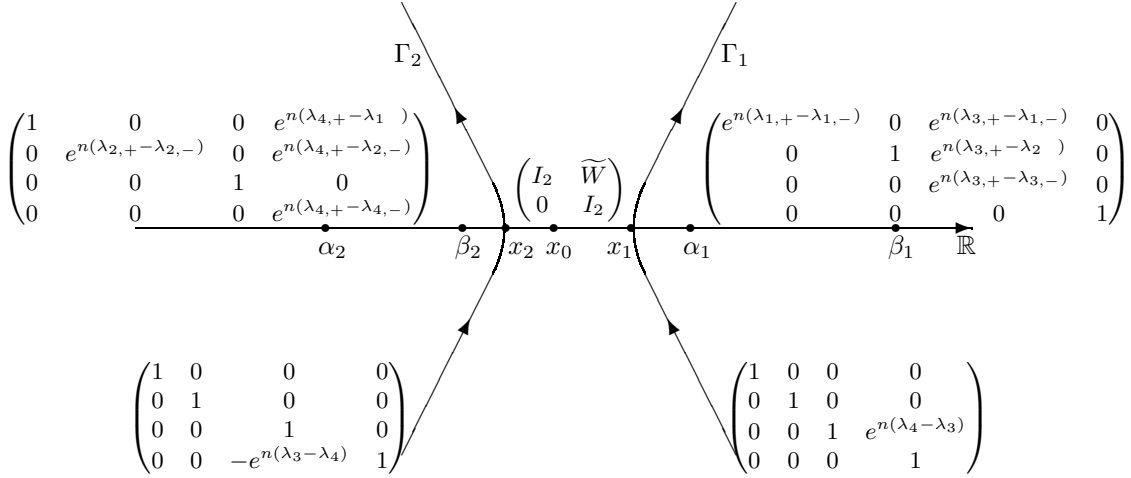

The jumps for $T(z)$ are on the same contours $\mathbb R \cup \Gamma_1 \cup \Gamma_2$
as those for $X(z)$, but with
different jump matrices.
To state the resulting jumps, we find it convenient to work with the following
functions (\lq $\lam$-functions\rq)
\begin{align}\label{deflambda14}
    \begin{array}{ll}   \lambda_1(z) &:= - g_1(z)+\frac{l_1}{2}+\frac{1}{2t}(z^2-2a_1z), \\
    \lambda_2(z) &:= -g_2(z)+\frac{l_2}{2} +\kappa+\frac{1}{2t}(z^2-2a_2z), \\
    \lambda_3(z) &:=  g_1(z)-\frac{l_1}{2} -\frac{1}{2(1-t)}(z^2-2b_1z), \\
    \lambda_4(z) &:=  g_2(z)-\frac{l_2}{2} +\kappa-\frac{1}{2(1-t)}(z^2-2b_2z).
    \end{array}
\end{align}
Using \eqref{deflambda14} and after a little calculation,
we find that we can write all jump matrices in
the RH problem for $T(z)$ as shown in Figure \ref{figjumpslambdas}.
On $(x_2,x_1)$ the jump matrix is
\begin{equation}
\label{deftildeW}
\begin{pmatrix} I_2 & \widetilde{W} \\ 0 & I_2 \end{pmatrix} \qquad \textrm{with} \quad
 \widetilde{W} = \begin{pmatrix}e^{-n\lambda_1(z)} \\ e^{-n\lambda_2(z)}\end{pmatrix}
 \begin{pmatrix}e^{n\lambda_3(z)} & e^{n\lambda_4(z)}\end{pmatrix}. \end{equation}

Thus $T$ satisfies the following RH problem.
\begin{itemize}
\item[(1)] $T(z)$ is analytic in $\mathbb C \setminus ( \mathbb R  \cup \Gamma_1 \cup \Gamma_2)$;
 \item[(2)] On $\mathbb R  \cup \Gamma_1 \cup \Gamma_2$ we have that $T_+ = T_- J_T$
 with jump matrices $J_T$ as shown  in Figure~\ref{figjumpslambdas}.
  \item[(3)] As $z\to\infty$, we have that
\begin{equation}\label{asymptoticconditionT} T(z) =
    I_{4}+O(1/z).
\end{equation}
\end{itemize}

\subsection{Asymptotic behavior of the jump matrices}
\label{subsectionasymptoticsjumpslarge}

In this subsection, we investigate the asymptotic behavior in the large $n$
limit of each of the jump matrices in the RH problem for $T=T(z)$ in Figure
\ref{figjumpslambdas}. Our goal is to show that all the jump matrices are
exponentially close (as $n\to\infty$) to the identity matrix, except for the
jump matrices on the intervals $(\alpha_1, \beta_1)$ and $(\alpha_2, \beta_2)$
which have an  oscillatory behavior. To be able do so we still have the freedom
to take the reference point $x_0$ and the contours $\Gamma_1$ and $\Gamma_2$ in
an appropriate way. We also have the constant $\kappa$ in \eqref{defLmx} at our
disposal.

\subsubsection{Preliminaries on $\xi$-functions}
In what follows we will make extensive use of the $\xi$-functions
which are essentially the derivatives of the $\lambda$-functions.
We define
\begin{equation}
    \label{defxij}
    \xi_j(z) = \lambda_j'(z) - \frac{(1-2t)}{2t(1-t)} z, \qquad j=1,\ldots, 4.
    \end{equation}
From \eqref{deflambda14} we then see that
\begin{align} \label{defxifunctions}
    \begin{array}{ll}
\xi_1(z) & =  - g'_1(z)+\frac{z}{t}-\frac{a_1}{t} - \frac{(1-2t)}{2t(1-t)} z, \\
\xi_2(z) & =  - g'_2(z)+\frac{z}{t}-\frac{a_2}{t} - \frac{(1-2t)}{2t(1-t)} z, \\
\xi_3(z) & =  g'_1(z)-\frac{z}{1-t}+\frac{b_1}{1-t} - \frac{(1-2t)}{2t(1-t)} z, \\
\xi_4(z) & =  g'_2(z)-\frac{z}{1-t}+\frac{b_2}{1-t} - \frac{(1-2t)}{2t(1-t)} z.
\end{array}
\end{align}
Inserting \eqref{diffg1}--\eqref{diffg2} and
\eqref{defalphavarying}--\eqref{defbetavarying} into these expressions, we
obtain after some simplification,
\begin{equation}\label{defxi14}
    \begin{array}{lll}
    \xi_1(z) = & \frac{1}{2t(1-t)}\left(-(1-t)a_1+tb_1+((z-\alpha_1)(z-\beta_1))^{1/2}\right), \\
    \xi_2(z) = & \frac{1}{2t(1-t)}\left(-(1-t)a_2+tb_2+((z-\alpha_2)(z-\beta_2))^{1/2}\right), \\
    \xi_3(z) = & \frac{1}{2t(1-t)}\left(-(1-t)a_1+tb_1-((z-\alpha_1)(z-\beta_1))^{1/2}\right), \\
    \xi_4(z) = & \frac{1}{2t(1-t)}\left(-(1-t)a_2+tb_2-((z-\alpha_2)(z-\beta_2))^{1/2}\right).
\end{array}
\end{equation}

\begin{figure}[tb]
\centering \begin{overpic}[scale=0.5,angle=270]%
{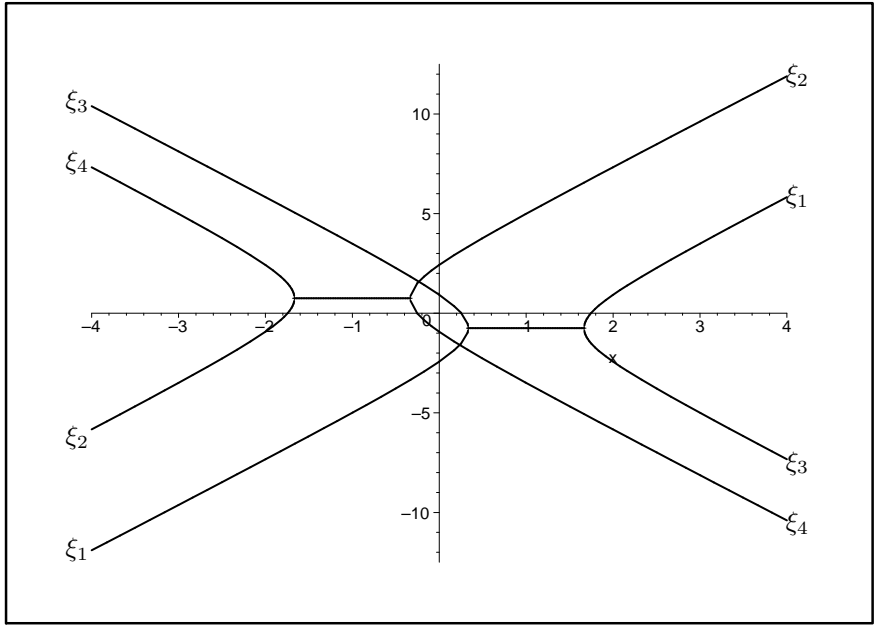}%
      \put(90,63){$\xi_2$}
      \put(90,49){$\xi_1$}
      \put(90,18){$\xi_3$}
      \put(90,11){$\xi_4$}
      \put(7,60){$\xi_3$}
      \put(7,53){$\xi_4$}
      \put(7,21){$\xi_2$}
      \put(7,8){$\xi_1$}
\end{overpic}
\caption{The figure shows the functions $\xi_1(x)$, $\xi_3(x)$ for
$x\in\mathbb R \setminus [\alpha_1,\beta_1]$ and $\xi_2(x)$, $\xi_4(x)$ for
$x\in\mathbb R \setminus [\alpha_2,\beta_2]$, together with the two horizontal line
segments $\Re(\xi_1(x))=\Re(\xi_3(x)) = \Xi_1$ for $x\in [\alpha_1,\beta_1]$ and
$\Re(\xi_2(x))=\Re(\xi_4(x)) = \Xi_2$ for $x\in [\alpha_2,\beta_2]$. We have here
$p_1=p_2=1/2$, $t=1/3$, $-a_2=a_1=1$, $-b_2=b_1=1$ and (hence) $\alpha_1\approx
0.33$, $\beta_1\approx 1.66$, $\alpha_2\approx -1.66$, $\beta_2\approx -0.33$
and $x_0=0$. The graphs of $\xi_2$ and $\xi_3$ intersect at the $x$-value
$-\sqrt{2}/6\approx -0.23$ and those of $\xi_1$ and $\xi_4$ intersect at the
$x$-value $\sqrt{2}/6\approx 0.23$.} \label{figxi1}
\end{figure}

From either \eqref{defxifunctions} or \eqref{defxi14} we see that
$\xi_j$ is defined and analytic in the cut plane $\mathbb C \setminus [\alpha_k, \beta_k]$
(with $k=j$ or $k=j-2$),
but $\Re \xi_j(z)$ is well-defined on the cut as well. From
\eqref{defxi14} we easily deduce the following behavior on the real line,
see also Figure \ref{figxi1}.

\begin{lemma} \label{lemmaonxi1}
\begin{enumerate}
\item[\rm (a)] We have that $\Re \xi_1$ and $\Re \xi_3$ are constant
on the interval $[\alpha_1, \beta_1]$:
\[ \Re \xi_1(x) = \Re \xi_3(x) = \Xi_1 := \frac{1}{2t(1-t)}
    \left(-(1-t)a_1+tb_1 \right), \]
for $x \in  [\alpha_1,\beta_1]$.
On $\mathbb R \setminus [\alpha_1, \beta_1]$
we have that $\xi_1$ and $\xi_3$
are real, $\xi_1$ is strictly increasing,
$\xi_3$ is strictly decreasing, and
\[  \frac{1}{2} \left( \xi_1(x) + \xi_3(x) \right) = \Xi_1, \qquad
    x \in \mathbb R \setminus [\alpha_1, \beta_1]. \]
\item[\rm (b)] We have that $\Re \xi_2$ and $\Re \xi_4$ are constant
on the interval $[\alpha_2, \beta_2]$:
\[ \Re \xi_2(x) = \Re \xi_4(x) = \Xi_2 := \frac{1}{2t(1-t)}
    \left(-(1-t)a_2+tb_2 \right), \]
for $x \in [\alpha_2,\beta_2]$.
On $\mathbb R \setminus [\alpha_2,\beta_2]$ we have that $\xi_2$ and $\xi_4$
are real, $\xi_2$ is strictly increasing,
$\xi_4$ is strictly decreasing, and
\[ \frac{1}{2} \left( \xi_2(x) + \xi_4(x) \right) = \Xi_2, \qquad
    x \in \mathbb R \setminus [\alpha_2, \beta_2]. \]
\item[\rm (c)] We have
\[ \Xi_2 \geq \Xi_1. \]
\end{enumerate}
\end{lemma}
\bewijs.
All statements immediately follow from \eqref{defxi14}, except maybe the
statement of part (c). Part (c) follows from the fact that
\[ \Xi_2 - \Xi_1 =  \frac{1}{2t(1-t)}
   \left( (1-t)(a_1-a_2) - t(b_1-b_2) \right) \]
which is indeed non-negative because of our assumption \eqref{hypothesist}.
$\bol$

The main property of the $\xi$-functions is contained in the following lemma.
It will only be valid under the large separation assumption, as we will see in
the proof.

\begin{lemma} \label{lemmaonxi2}
There is a value $x_0 \in (\beta_2, \alpha_1)$ such that
\begin{equation} \label{xiinequalities}
    \xi_2(x_0) \geq \xi_3(x_0) > \xi_4(x_0) \geq \xi_1(x_0).
    \end{equation}
\end{lemma}

\bewijs.
We prove the lemma in two steps.
In the first step we show
that there exists $x_0 \in (\beta_2, \alpha_1)$ such that
\[ \xi_3(x_0) > \xi_4(x_0). \]
This we can do by giving $x_0$ the explicit value
\begin{equation}\label{defx0}
    x_0 = \frac{\sqrt{p_1}}{\sqrt{p_1}+\sqrt{p_2}}\frac{\alpha_2+\beta_2}{2}+
    \frac{\sqrt{p_2}}{\sqrt{p_1}+\sqrt{p_2}}\frac{\alpha_1+\beta_1}{2}
\end{equation}
so that clearly
\[ \frac{\alpha_2 + \beta_2}{2} < x_0 < \frac{\alpha_1 + \beta_1}{2} \]
which means that $x_0$ lies between the midpoints of the two intervals, and so
in particular $\alpha_2 < x_0 < \beta_1$.

We have by \eqref{defalphavarying}--\eqref{defbetavarying} and
\eqref{defx0} that
\begin{align} \nonumber
     & (\alpha_1-x_0)  (\beta_1-x_0) \\
     & = \nonumber
    \left(\frac{\alpha_1+\beta_1}{2}-x_0\right)^2-\left(\frac{\beta_1-\alpha_1}{2}\right)^2  \\
    & = \nonumber \frac{p_1}{(\sqrt{p_1}+\sqrt{p_2})^2} \left(\frac{\alpha_1+\beta_1-\alpha_2-\beta_2}{2}\right)^2
        -4p_1 t(1-t) \\
    & = \frac{p_1}{(\sqrt{p_1} + \sqrt{p_2})^2}
        \left( ((1-t)(a_1-a_2) + t(b_1-b_2))^2 - 4 t(1-t) (\sqrt{p_1} + \sqrt{p_2})^2 \right).
        \label{productx00}
\end{align}
Because we are assuming that we are in a situation of large separation
we have $(\sqrt{p_1} + \sqrt{p_2})^2 < (a_1-a_2)(b_1-b_2)$
so that we obtain from \eqref{productx00}
\begin{align} \nonumber
     & (\alpha_1-x_0)  (\beta_1-x_0) \\ \nonumber
    & > \frac{p_1}{(\sqrt{p_1} + \sqrt{p_2})^2}
        \left( ((1-t)(a_1-a_2) + t(b_1-b_2))^2 - 4 t(1-t) (a_1-a_2)(b_1-b_2) \right) \\
    & = \frac{p_1}{(\sqrt{p_1} + \sqrt{p_2})^2}
         ((1-t)(a_1-a_2) - t(b_1-b_2))^2.
        \label{productx01}
        \end{align}
Similarly,
\begin{align} \nonumber
    &(x_0-\alpha_2)(x_0-\beta_2) \\
    & = \nonumber
    \frac{p_2}{(\sqrt{p_1}+\sqrt{p_2})^2}\left(((1-t)(a_1-a_2)+t(b_1-b_2))^2
    -4t(1-t)(\sqrt{p_1}+\sqrt{p_2})^2\right) \\
    & >  \frac{p_2}{(\sqrt{p_1} + \sqrt{p_2})^2}
         ((1-t)(a_1-a_2) - t(b_1-b_2))^2.
        \label{productx02}
\end{align}
It then follows that both factors on the left-hand sides of \eqref{productx01} and \eqref{productx02}
are positive, so that
\[ \beta_2 < x_0 < \alpha_1, \]
as claimed.

Furthermore, we obtain from  \eqref{productx01} and \eqref{productx02}
\begin{multline} \label{derivativenegative2}
    \sqrt{(\alpha_1- x_0)(\beta_1 - x_0)}+\sqrt{(x_0-\alpha_2)(x_0-\beta_2)} \\
    > \left| (1-t)(a_1-a_2) - t(b_1-b_2) \right|.
\end{multline}
Since by \eqref{defxi14} we have
\begin{multline}
    2t(1-t) \left( \xi_3(x_0) - \xi_4(x_0) \right) \\
    =
    \sqrt{(\alpha_1- x_0)(\beta_1 - x_0)}+\sqrt{(x_0-\alpha_2)(x_0-\beta_2)} \\
    - \left( (1-t)(a_1-a_2) - t(b_1-b_2) \right),
    \end{multline}
    it follows from \eqref{derivativenegative2} that
$\xi_3(x_0) > \xi_4(x_0)$, as claimed as well.

\medskip

In the second step we prove that the string of inequalities
\eqref{xiinequalities} holds
for some value $x_0 \in (\beta_2, \alpha_1)$ which, however, could be different
from $x_0$ used in the first step.
We distinguish three cases, depending on whether
$\xi_3(\beta_2) < \Xi_2$, $\xi_3(\beta_2) = \Xi_2$, or $\xi_3(\beta_2) > \Xi_2$ where
$\Xi_2$ is the value of $\Re \xi_2(x) = \Re \xi_4(x)$ for $x \in [\alpha_2, \beta_2]$,
see part (b) of Lemma \ref{lemmaonxi1}.

If $\xi_3(\beta_2) < \Xi_2$, then $\xi_3(\beta_2) < \xi_4(\beta_2)$,
and since by what we already proved, the inequality $\xi_3(x) > \xi_4(x)$
holds for some $x  \in (\beta_2, \alpha_1)$, there also exists $x^* \in (\beta_2, \alpha_1)$
where equality $\xi_3(x^*) = \xi_4(x^*)$ holds.
Then by the monotonicity properties of the $\xi$-functions (stated
in parts (a) and (b) of Lemma \ref{lemmaonxi1})
\[ \xi_2(x^*) > \xi_3(x^*) = \xi_4(x^*) > \xi_1(x^*). \]
Then by slightly increasing $x^*$ we find $x_0$ such that the inequalities
\eqref{xiinequalities} hold (with strict inequalities).

If $\xi_3(\beta_2) = \Xi_2$, then as above we have
\[ \xi_2(\beta_2) = \xi_3(\beta_2) = \xi_4(\beta_2) > \xi_1(\beta_2). \]
Then for $x_0$ slightly larger than $\beta_2$, we have
 \eqref{xiinequalities} with strict inequalities,  since
$\xi_2'(x_0) = +\infty$ and $\xi_4'(x_0) = -\infty$, see
formulas \eqref{defxi14}.

If $\xi_3(\beta_2) > \Xi_2$, then $\xi_3(\beta_2) > \xi_2(\beta_2)$.
Since $\xi_2$ is increasing we find  $\xi_2(\alpha_1) > \xi_2(\beta_2) = \Xi_2 \geq \Xi_1
= \xi_3(\alpha_1)$, see Lemma \ref{lemmaonxi1}. Thus there exists $x_0 \in (\beta_2, \alpha_1)$
with $\xi_2(x_0) = \xi_3(x_0)$.
We have $\xi_4(x_0) < \xi_2(x_0)$, $\xi_1(x_0) < \xi_3(x_0)$,
and
\[ \xi_1(x_0) + \xi_3(x_0) = 2\Xi_1 \leq 2\Xi_2 = \xi_2(x_0) + \xi_4(x_0). \]
again by  Lemma \ref{lemmaonxi1}, so that $\xi_1(x_0) \leq \xi_4(x_0)$
and the inequalities \eqref{xiinequalities} hold.
$\bol$

It follows from the proof that in most cases we may assume
that the inequalities \eqref{xiinequalities} are strict. Only if $t= t_{\crit}$
then $\Xi_1=\Xi_2$ and then we can only obtain
\[ \xi_2(x_0) = \xi_3(x_0) > \xi_4(x_0) = \xi_1(x_0). \]

\subsubsection{Jump matrices on $\Gamma_1$ and $\Gamma_2$}

After the preliminaries on the $\xi$-functions we consider the jump matrices on
$\Gamma_1$ and $\Gamma_2$. We want that the jump matrices  on these curves,
shown in Figure \ref{figjumpslambdas}, are exponentially close to the identity
matrix as $n\to\infty$. Thus we want  to choose $\Gamma_1$ and $\Gamma_2$ so
that
\begin{align} \label{signlambda34}
    \Re( \lambda_4 - \lambda_3)
    \left\{ \begin{array}{ll}
        < -c_1 < 0,  &  \textrm{ on } \Gamma_1, \\
        > c_2 > 0, & \textrm{ on } \Gamma_2
        \end{array} \right.
\end{align}
for certain constants $c_1$, $c_2$ that do not depend on $n$.

We take the point $x_0$ satisfying the inequalities \eqref{xiinequalities}
of Lemma \ref{lemmaonxi2}.
We have by \eqref{deflambda14}
\begin{equation}\label{realpartinteresting}
    \lambda_4(z)-\lambda_3(z) =  g_{2}(z)- \frac{l_2}{2} + \kappa
    - g_{1}(z)+\frac{l_1}{2} - \frac{b_1-b_2}{1-t}z.
\end{equation}
The constant $\kappa$ is still at our disposal. We choose
it here so
that \eqref{realpartinteresting} vanishes for $z = x_0$.
Thus
\[ \lambda_4(x_0) - \lambda_3(x_0) = 0. \]

Since the derivative is
\[ \left. \frac{d}{dx} \left(\lambda_4(x) - \lambda_3(x) \right) \right|_{x=x_0}
    = \xi_4(x_0) - \xi_3(x_0) < 0, \]
    see \eqref{defxij}, we can then find $x_1$ and $x_2$ sufficiently close
    to $x_0$ so that $\beta_2 < x_2 < x_0 < x_1 < \alpha_1$ and
\[ \lambda_4(x_1) - \lambda_3(x_1) < 0 < \lambda_4(x_2) - \lambda_3(x_2). \]

We will choose $\Gamma_1$ and $\Gamma_2$ so that they cross the
real axis in $x_1$ and $x_2$, respectively. It remains to show that
they can be extended to infinity, while remaining in
the open sets
\begin{align}\label{defOmega1}
\Omega_1 & := \{z\in\mathbb C  \mid \Re(\lambda_4(z)-\lambda_3(z))<0\}, \\
    \label{defOmega2}
\Omega_2 & := \{z\in\mathbb C  \mid \Re(\lambda_4(z)-\lambda_3(z))>0\},
\end{align}
respectively. Note that $\Omega_1$ and $\Omega_2$ are open
subsets of $\mathbb C$, symmetric with respect to the real axis and
such that
\[ x_1 \in \Omega_1, \qquad x_2 \in \Omega_2. \]
The fact that we can indeed choose $\Gamma_1$ and $\Gamma_2$ this way,
follows from the following lemma.

\begin{lemma}\label{lemmacanchoosecurves}
For $j=1,2$, let $\Omega_j^o$ denote the connected component of $\Omega_j$ that
contains $x_j$. Then $\Omega_j^o$ is unbounded. In addition, for each $\varepsilon > 0$
there exists $R > 0$ so that
\begin{align} \label{Omega1unbounded}
    \{ z \in \mathbb C \mid  |z| > R, \,  -\pi/2 + \varepsilon < \arg z < \pi/2 - \varepsilon \}
    & \subset \Omega_1^o,
    \end{align}
    and
    \begin{align} \label{Omega2unbounded}
     \{ z \in \mathbb C \mid |z| > R, \, \pi/2 + \varepsilon < \arg z < 3 \pi/2 - \varepsilon \}
    & \subset \Omega_2^o. \end{align}
\end{lemma}

\bewijs. We will establish these properties by using the maximum principle for
subharmonic functions, where we use that $\Re g_j$ as given by \eqref{defg1g2}
is subharmonic on $\mathbb C$ and harmonic on $\mathbb C \setminus [\alpha_j, \beta_j]$
for $j=1,2$.

Suppose that $\overline{\Omega_1^o}$ (the closure of $\Omega_1^o$) has
nonempty intersection with the interval $[\alpha_2, \beta_2]$.
Since $\Omega_1^o$ is connected, and symmetric with respect to
the real axis, it will then surround the point $x_2$, so that
$\Omega_2^o$ must be bounded and
\[ \overline{\Omega_2^o} \cap [\alpha_1, \beta_1] = \emptyset. \]
Then $\Re g_1$ is harmonic on $\overline{\Omega_2^o}$,
$\Re g_2$ is subharmonic and therefore by \eqref{realpartinteresting}
\begin{equation} \label{realpartinteresting2}
    \Re(\lambda_4 - \lambda_3)(z) =
    \Re g_2(z) - \Re g_1(z)
    - \frac{b_1-b_2}{1-t} \Re z + \textrm{const}, \end{equation}
is subharmonic on $\overline{\Omega_2^o}$.
By definition we have that $\Re (\lambda_4 - \lambda_3) > 0$ on
$\Omega_2^o$ with equality on $\partial \Omega_2^o$.
This is a contradiction with the maximum principle for
subharmonic functions and it follows that
\[ \overline{\Omega_1^o} \cap [\alpha_2, \beta_2] = \emptyset. \]
Then if $\Omega_1^o$ is bounded we obtain a contradiction with
the maximum principle in the same way.
Thus $\Omega_1^o$ is unbounded, and likewise $\Omega_2^o$ is
unbounded as well.

From \eqref{realpartinteresting2} we further see that
\begin{align}
     \Re(\lambda_4(z)-\lambda_3(z))
    \label{realpartinfty} &= -\frac{b_1-b_2}{1-t}\Re z +O(\log(|z|))
\end{align}
as $z\to\infty$, since $g_j(z) = p_j \log z +O(1/z)$ as $z \to \infty$.
In addition, we have
\begin{align} \label{realpartinfty2}
    \frac{d}{dz} (\lambda_4(z)-\lambda_3(z)) &=
    \xi_4(z) - \xi_3(z) =
    -\frac{b_1-b_2}{1-t}+O(1/z)
\end{align}
as $z \to \infty$.
From \eqref{realpartinfty}--\eqref{realpartinfty2} it easily follows that
for $|\Im z|$ large enough, we have that $\Re(\lambda_4(z)-\lambda_3(z))$
monotonically decreases as $\Re z$ increases.
This implies that the domains $\Omega_1$ and $\Omega_2$ given in
\eqref{defOmega1}--\eqref{defOmega2} both have only one
unbounded component, which then coincide with $\Omega_1^o$ and $\Omega_2^o$
since these are unbounded. It now also follows from \eqref{realpartinfty}
that \eqref{Omega1unbounded}
and \eqref{Omega2unbounded} hold. $\bol$

We conclude from Lemma \ref{lemmacanchoosecurves} that the curves
$\Gamma_1$, $\Gamma_2$ can be extended to infinity so that
\[ \Re( \lambda_4 - \lambda_3)
    \left\{ \begin{array}{ll}
        < -c_1 < 0,  &  \textrm{ on } \Gamma_1, \\
        > c_2 > 0, & \textrm{ on } \Gamma_2
        \end{array} \right.
\]
for some constants $c_1, c_2 > 0$. Since as $n \to \infty$
we are in a non-critical situation, the contours $\Gamma_j$
and the constants $c_j$ can be taken independently of $n$
for $n$ large enough.
Then the off-diagonal
entries in the jump matrices on $\Gamma_1$ and $\Gamma_2$ are
uniformly exponentially small
as $n \to \infty$, as required.

\subsubsection{Jump matrices on $\mathbb R$}

We next investigate the jump matrices in the RH problem for $T$ on the various
parts of the real line. Recall that the jump matrices are given in Figure
\ref{figjumpslambdas} and so we are interested in the behavior of $
\lambda_{k,+}- \lambda_{j,-}$ on the real line for various combinations of $j$
and $k$.

We prove:
\begin{lemma} \label{lemmajumpsonR}
The following hold.
\begin{enumerate}
\item[\rm (a)] For $j=1,2$, we have
\begin{align}
\label{Eulerlagrange5}
    \lambda_{j+2,+}-\lambda_{j,-}=0 &  \quad \textrm{ on } [\alpha_j,\beta_j],\\
\label{Eulerlagrange6}
    \lambda_{j+2,+}-\lambda_{j,-}<0 &  \quad \textrm{ on } \mathbb R \setminus [\alpha_j,\beta_j].
\end{align}
\item[\rm (b)]
We have
\begin{align}
\label{lambdaineq41}
    \Re(\lambda_{4,+} - \lambda_1) < 0 & \quad \textrm{ on } (-\infty, x_1), \\
\label{lambdaineq32}
    \Re(\lambda_{3,+} - \lambda_2) < 0 &  \quad \textrm{ on } (x_2, \infty).
\end{align}
\end{enumerate}
\end{lemma}

\bewijs.
(a) From the definitions
\eqref{deflambda14} and \eqref{defV1}--\eqref{defV2} it follows that for $j=1,2$,
\[ \lambda_{j+2,+}-\lambda_{j,-} = g_{j,+} - g_{j,-} - l_j - p_j V_j, \]
so that part (a) is a restatement of the Euler-Lagrange conditions
\eqref{Eulerlagrange1}--\eqref{Eulerlagrange4}.

(b) Recall that the constant $\kappa$ was chosen so that
$\lambda_3(x_0) = \lambda_4(x_0)$.
Then in view of part (a) we know that
\begin{align*}
    \lambda_4(x_0)-\lambda_1(x_0) = \lambda_3(x_0)-\lambda_1(x_0)<0,\\
    \lambda_3(x_0)-\lambda_2(x_0) = \lambda_4(x_0)-\lambda_2(x_0)<0.
\end{align*}
Taking $x_1$ and $x_2$ sufficiently close to $x_0$ (which we can
do without loss of generality), we
then have the inequalities \eqref{lambdaineq41}--\eqref{lambdaineq32}
on the interval $(x_2,x_1)$.

From Lemma \ref{lemmaonxi1} we know that $\Re \xi_4$ is decreasing while $\Re \xi_1$ is strictly
increasing on $(-\infty, \alpha_1]$. Then we have for $x < x_0$,
\[ \Re (\xi_4 - \xi_1)(x) > \Re (\xi_4 - \xi_1)(x_0)
    = \xi_4(x_0) - \xi_1(x_0) \geq 0 \]
where the last inequality holds because of Lemma \ref{lemmaonxi2}.
Since
\[ \frac{d}{dx} \Re (\lambda_4 - \lambda_1)
    = \Re ( \xi_4 - \xi_1) \]
it then follows that $\Re (\lambda_4 - \lambda_1)$ is strictly increasing
on $(-\infty, x_0)$. Since we already proved that the
inequality \eqref{lambdaineq41} holds for $x \in (x_2,x_1)$,
it then also follows for $x < x_0$.

The inequality \eqref{lambdaineq32} for $x > x_0$ follows in the
same way.
$\bol$

It follows from Lemma \ref{lemmajumpsonR} that, up to exponentially
small corrections, the jump matrices on the real line in the RH problem for $T$ take
the following form
\begin{align} \label{jumpTapprox1}
    \begin{pmatrix} 1 & 0 & 0 & 0 \\
   0 & e^{n(\lambda_{2,+}-\lambda_{2,-})} & 0 & e^{n(\lambda_{4,+}-\lambda_{2,-})} \\
   0 & 0 & 1 & 0 \\
   0 & 0 & 0 & e^{n(\lambda_{4,+}-\lambda_{4,-})}
   \end{pmatrix} \qquad \textrm{on } (-\infty, x_2) \end{align}
\begin{align} \label{jumpTapprox2}
    \begin{pmatrix}
   e^{n(\lambda_{1,+}-\lambda_{1,-})} & 0 & e^{n(\lambda_{3,+}-\lambda_{1,-})} & 0 \\
    0 & 1 & 0 & 0 \\
    0 & 0 & e^{n(\lambda_{3,+}-\lambda_{3,-})} & 0 \\
    0 & 0 & 0 & 1
   \end{pmatrix} \qquad \textrm{on } (x_1, \infty), \end{align}
and $\begin{pmatrix} I_2 & 0 \\ 0 & I_2 \end{pmatrix}$ on $(x_2,x_1)$.

The matrices \eqref{jumpTapprox1} and \eqref{jumpTapprox2} are
in a standard form.
By standard arguments we have that the non-trivial diagonal entries
of \eqref{jumpTapprox1} and \eqref{jumpTapprox2} are rapidly
oscillating on the intervals $[\alpha_2,\beta_2]$ and $[\alpha_1, \beta_1]$, respectively,
and they are equal to $1$ outside these intervals.

From part (a) of Lemma \ref{lemmajumpsonR} we can further deduce that
\eqref{jumpTapprox1} reduces to
\begin{align} \label{jumpTapprox3}
    \begin{pmatrix} 1 & 0 & 0 & 0 \\
   0 & e^{n(\lambda_{2,+}-\lambda_{2,-})} & 0 & 1 \\
   0 & 0 & 1 & 0 \\
   0 & 0 & 0 & e^{n(\lambda_{4,+}-\lambda_{4,-})}
   \end{pmatrix} \qquad \textrm{on } [\alpha_2,\beta_2] \end{align}
and \eqref{jumpTapprox2} reduces to
\begin{align} \label{jumpTapprox4} \begin{pmatrix}
   e^{n(\lambda_{1,+}-\lambda_{1,-})} & 0 & 1 & 0 \\
    0 & 1 & 0 & 0 \\
    0 & 0 & e^{n(\lambda_{3,+}-\lambda_{3,-})} & 0 \\
    0 & 0 & 0 & 1
   \end{pmatrix} \qquad \textrm{on } [\alpha_1, \beta_1]. \end{align}

Outside these intervals, the jump matrices \eqref{jumpTapprox1} and \eqref{jumpTapprox2}
are exponentially close to the identity.

\subsubsection{Summary}
Summarizing, we have shown that we can take $\Gamma_1$ and $\Gamma_2$
so that all jumps in the RH problem for $T$ tend to the identity
matrix as $n \to \infty$, except for the jump matrices on the two
intervals $[\alpha_j, \beta_j]$,
$j=1,2$, which are given by \eqref{jumpTapprox3}--\eqref{jumpTapprox4}
plus an exponentially small term.

Observe that each of the asymptotic jump matrices above has the sparsity
pattern
\begin{equation}\label{sparsity}
\begin{pmatrix} \times & 0 & \times & 0 \\
0 & \times & 0 & \times \\
\times & 0 & \times & 0 \\
0 & \times & 0 & \times
 \end{pmatrix}.
\end{equation}
It follows
that the $4\times 4$ RH problem for $T=T(z)$ asymptotically decouples into two
RH problems of size $2\times 2$, one involving rows and columns $1$ and $3$ and
the other involving rows and columns $2$ and $4$. The coupling between the
two RH problems is exponentially small as $n \to \infty$. Thus we are dealing now
essentially with two RH problems of size $2\times 2$. The remaining steps in
the Deift-Zhou steepest descent analysis can then be done in a standard way on
these decoupled $2\times 2$ problems. This will be described in the next
subsection.

\subsection{Remaining steps of the steepest descent analysis}

\subsubsection{Third transformation: Opening of the lenses}
\label{subsectiontransfo3large}

 The next step in the steepest descent analysis
is to open lenses around the intervals $(\alpha_j,\beta_j)$, $j=1,2$. This
operation serves to transform the oscillating jumps on these intervals into
exponentially decaying or constant ones. This can be done in the usual way
\cite{Dei} on each of the decoupled $2\times 2$ problems. We obtain in this way
a new matrix function $S=S(z)$ obtained from $T(z)$ by multiplication on the
right with a suitable transformation matrix, the precise form of which depends
on the different regions of the complex plane.

Let us describe this in detail. First consider the interval
$[\alpha_1,\beta_1]$. We have here the matrix factorization
\begin{eqnarray*} & & \begin{pmatrix} e^{n(\lambda_{1,+}-\lambda_{1,-})} & 1 \\ 0 & e^{n(\lambda_{3,+}-\lambda_{3,-})}
\end{pmatrix}\\
&=& \begin{pmatrix} e^{n(\lambda_{1,+}-\lambda_{3,+})} & 1 \\ 0 &
e^{n(\lambda_{1,-}-\lambda_{3,-})} \end{pmatrix}\\
&=& \begin{pmatrix} 1 & 0 \\ e^{n(\lambda_{1,-}-\lambda_{3,-})} & 1
\end{pmatrix}
\begin{pmatrix} 0 & 1 \\ -1 & 0 \end{pmatrix}
\begin{pmatrix} 1 & 0 \\ e^{n(\lambda_{1,+}-\lambda_{3,+})} & 1 \end{pmatrix}.
\end{eqnarray*}
We can then open a lens around $[\alpha_1,\beta_1]$ and define
\begin{equation}\label{defS1} S(z) = \left\{ \begin{array}{ll}
T(z) \begin{pmatrix} 1 & 0 & 0 & 0 \\ 0 & 1 & 0 & 0 \\
-e^{n(\lambda_{1}-\lambda_{3})} & 0 & 1 & 0 \\ 0 & 0 & 0 & 1 \end{pmatrix}
\textrm{ in upper lens region around $[\alpha_1,\beta_1]$}\\
T(z) \begin{pmatrix} 1 & 0 & 0 & 0 \\ 0 & 1 & 0 & 0 \\
e^{n(\lambda_{1}-\lambda_{3})} & 0 & 1 & 0 \\ 0 & 0 & 0 & 1
\end{pmatrix}\textrm{ in lower lens region around $[\alpha_1, \beta_1]$.}
\end{array}\right.
\end{equation}
Similarly we open a lens around $[\alpha_2,\beta_2]$ and define
\begin{equation}\label{defS2} S(z) = \left\{ \begin{array}{ll}
T(z) \begin{pmatrix} 1 & 0 & 0 & 0 \\ 0 & 1 & 0 & 0 \\
0 & 0 & 1 & 0 \\ 0 & -e^{n(\lambda_{2}-\lambda_{4})} & 0 & 1 \end{pmatrix}\textrm{ in upper lens region around $[\alpha_2,\beta_2]$}\\
T(z) \begin{pmatrix} 1 & 0 & 0 & 0 \\ 0 & 1 & 0 & 0 \\
0 & 0 & 1 & 0 \\ 0 & e^{n(\lambda_{2}-\lambda_{4})} & 0 & 1
\end{pmatrix}\textrm{ in lower lens region around $[\alpha_2, \beta_2]$.}
\end{array}\right.
\end{equation}
We also set
\begin{equation}\label{defS3} S(z) = T(z),\quad \textrm{ outside the lenses.}
\end{equation}
The matrix function $S=S(z)$ satisfies the following RH problem
\begin{itemize}
\item[(1)] $S(z)$ is analytic in $\mathbb C \setminus (\mathbb R \cup \Gamma_1\cup\Gamma_2\cup \textrm{lips of the lenses})$;
\item[(2)] $S$ satisfies the following jumps:
\begin{itemize}
\item[] $S_+(x) = S_-(x)\begin{pmatrix} 0 & 0 & 1 & 0 \\ -e^{n(\lambda_{1,+}-\lambda_2)} & 1 & e^{n(\lambda_{3,+}-\lambda_2)} & 0 \\
-1 & 0 & 0 & 0 \\ 0 & 0 & 0 & 1 \end{pmatrix}$ on $(\alpha_1,\beta_1)$
\item[] $S_+(x) = S_-(x)\begin{pmatrix} 1 & -e^{n(\lambda_{2,+}-\lambda_1)} & 0 & e^{n(\lambda_{4,+}-\lambda_1)} \\ 0 & 0 & 0 & 1 \\
0 & 0 & 1 & 0 \\ 0 & -1 & 0 & 0 \end{pmatrix}$ on $(\alpha_2,\beta_2)$
\item[] $S_+(x) = S_-(x)\begin{pmatrix} 1 & 0 & 0 & 0 \\ 0 & 1 & 0 & 0 \\
e^{n(\lambda_1-\lambda_3)} & 0 & 1 & 0 \\ 0 & 0 & 0 & 1 \end{pmatrix}$ on the lips of
the lens around $[\alpha_1,\beta_1]$,
\item[] $S_+(x) = S_-(x)\begin{pmatrix} 1 & 0 & 0 & 0 \\ 0 & 1 & 0 & 0 \\
0 & 0 & 1 & 0 \\ 0 & e^{n(\lambda_2-\lambda_4)} & 0 & 1 \end{pmatrix}$ on the lips of
the lens around $[\alpha_2, \beta_2]$.
\item[] On the remaining contours
$\Gamma_1,\Gamma_2,(-\infty,\alpha_2),(\beta_2,\alpha_1)$ and
$(\beta_1,\infty)$, the jumps for $S$ are exactly the same as those for $T$.
\end{itemize}
  \item[(3)] As $z\to\infty$, we have that
\begin{equation*}S(z) =
    I+O(1/z).
\end{equation*}
\end{itemize}
As already mentioned before, the entries in positions $(1,2)$, $(1,4)$, $(2,1)$
and $(2,3)$ in the jump matrices for $S$ are all exponentially small when
$n\to\infty$ (use parts (a) and (b) of
Lemma \ref{lemmajumpsonR}).

\subsubsection{Model RH problem: Parametrix away from the branch points}
\label{subsectionPinftylarge}

We consider now the model RH problem obtained from the RH problem for $S=S(z)$
by ignoring all exponentially small entries in the jump matrices. The model RH
will be defined in the region \begin{equation*}\mathbb C \setminus
([\alpha_1,\beta_1]\cup[\alpha_2,\beta_2]).\end{equation*} The solution
$P^{(\infty)}(z)$ to this model RH problem can be constructed in the usual way
\cite{Dei} for each of the two $2\times 2$ problems into which the $4\times 4$
problem decouples. This leads to the parametrix
\begin{equation}\label{Pinfty}
P^{(\infty)}(z) = \begin{pmatrix} \frac{\gamma_1(z)+\gamma_1^{-1}(z)}{2} & 0 &
\frac{\gamma_1(z)-\gamma_1^{-1}(z)}{2i}& 0\\
0&\frac{\gamma_2(z)+\gamma_2^{-1}(z)}{2} & 0 & \frac{\gamma_2(z)-\gamma_2^{-1}(z)}{2i} \\
-\frac{\gamma_1(z)-\gamma_1^{-1}(z)}{2i} & 0 & \frac{\gamma_1(z)+\gamma_1^{-1}(z)}{2} & 0 \\
0&-\frac{\gamma_2(z)-\gamma_2^{-1}(z)}{2i} & 0 &
\frac{\gamma_2(z)+\gamma_2^{-1}(z)}{2}
\end{pmatrix}
\end{equation}
where
\begin{equation}\label{defBD}
\gamma_1(z):=\left(\frac{z-\beta_1}{z-\alpha_1}\right)^{1/4},\quad
\gamma_2(z):=\left(\frac{z-\beta_2}{z-\alpha_2}\right)^{1/4},
\end{equation}
and where we choose the principal branches of the $1/4$ powers.

\subsubsection{Local parametrices around the branch points}
\label{subsectionPAirylarge}

Consider disks around the branch points $\alpha_j,\beta_j$, $j=1,2$ with
sufficiently small radius. Inside these disks one can construct local
parametrices $P^{(\textrm{Airy})}(z)$ to the RH problem for $S(z)$ in terms of
Airy functions. Once again, these parametrices can be constructed in the usual
way \cite{Dei} on each of the two $2\times 2$ RH problems.
We omit further details.

\subsubsection{Fourth transformation and completion of the steepest descent analysis}

Define a final matrix-valued function $R(z)$ by
\begin{equation} \label{defR0}
    R(z) = \left\{
    \begin{array}{ll}
    S(z)(P^{(\textrm{Airy})})^{-1}(z),& \quad \textrm{in the disks around } \alpha_1,\beta_1,\alpha_2,\beta_2\\
    S(z)(P^{(\infty)})^{-1}(z),& \quad \textrm{elsewhere}.
\end{array}
\right.
\end{equation}

From the construction of the parametrices it then follows that $R$ satisfies a RH problem
\begin{itemize}
\item[(1)] $R(z)$ is analytic in $\mathbb C \setminus \Sigma_R$
where $\Sigma_R$ consists of the real line $\mathbb R$, the contours $\Gamma_1$ and $\Gamma_2$,
the lips of the lenses outside the disks, and the boundaries of the disks
around $\alpha_1$, $\beta_1$,    $\alpha_2$, and $\beta_2$,
\item[(2)] $R$ has jumps $R_+ = R_- J_R$ on $\Sigma_R$,
that satisfy
\begin{align*}
    J_R(z) & = I + O(1/n), \qquad \text{on the boundaries of the disks}, \\
    J_R(z) & = I + O(e^{-cn|z|}), \qquad \text{on the other parts of $\Sigma_R$},
    \end{align*}
    for some constant $c >0$.
\item[(3)] $R(z) = I + O(1/z)$ as $z \to \infty$.
\end{itemize}

Then, as in \cite{Dei,DKMVZ1,DKMVZ2} we may conclude that
\begin{equation} \label{eq:Rasymptotics}
    R(z) = I_4 + O\left(\frac{1}{n(|z|+1)}\right)
    \end{equation}
    as $n \to \infty$, uniformly for $z$ in the complex plane.
    This completes the RH steepest descent analysis.

\subsection{Proof of Theorem \ref{st2ellipses} in the case of large separation}
\label{subsectionproof2ellipses}

We will now establish Theorem \ref{st2ellipses} in the case of large
separation. The idea is that from the asymptotic decoupling of the $4\times 4$
RH problem for large $n$, there follows a similar decoupling for the kernel in
\eqref{correlationkernel} and hence for the associated non-intersecting
Brownian particles.

First assume that $x,y\in (\alpha_1^*,\beta_1^*)$ and consider the correlation
kernel \eqref{correlationkernel} $K_{n_1,n_2;n_1,n_2}$, which for short we
denote by $K_n$. By virtue of \eqref{defX1} we obtain
\begin{equation*} \label{correlationkernelX}
    K_n(x,y) = \frac{1}{2\pi i(x-y)}\begin{pmatrix} 0 & 0 & w_{2,1}(y) &
    0\end{pmatrix} X_{+}^{-1}(y)X_{+}(x)\begin{pmatrix} w_{1,1}(x)\\
    w_{1,2}(x)\\ 0 \\ 0 \end{pmatrix}.
\end{equation*}
Using \eqref{defT} and \eqref{deflambda14} this becomes
\begin{equation*} \label{correlationkernelT}
    K_n(x,y) = \frac{1}{2\pi i(x-y)}\begin{pmatrix} 0 & 0 & e^{n\lambda_{3,+}(y)} &
    0\end{pmatrix} T_{+}^{-1}(y)T_{+}(x)\begin{pmatrix} e^{-n\lambda_{1,+}(x)}\\
    e^{-n\lambda_{2}(x)}\\ 0 \\ 0 \end{pmatrix}.
\end{equation*}
From \eqref{defS1} we get
\begin{equation} \label{correlationkernelS}
    K_n(x,y) = \frac{1}{2\pi i(x-y)}\begin{pmatrix} -e^{n\lambda_{1,+}(y)} & 0 & e^{n\lambda_{3,+}(y)} &
    0\end{pmatrix} S_{+}^{-1}(y)S_{+}(x)\begin{pmatrix} e^{-n\lambda_{1,+}(x)}\\
    e^{-n\lambda_{2}(x)}\\ e^{-n\lambda_{3,+}(x)} \\ 0 \end{pmatrix}.
\end{equation}
Now it follows from \eqref{eq:Rasymptotics} by standard arguments (e.g.\
\cite[Section 9]{BK2}) that
\begin{equation*}
S^{-1}(y)S(x) = I+O(x-y),\quad \textrm{as }y\to x
\end{equation*}
uniformly in $n$. Inserting this into \eqref{correlationkernelS} yields
\begin{multline} \label{correlationkernelS2}
    K_n(x,y)=
    \frac{-e^{n(\lambda_{1,+}(y)-\lambda_{1,+}(x))}+e^{n(\lambda_{3,+}(y)-\lambda_{3,+}(x))}}{2\pi
    i(x-y)}+O(1) \\+O(e^{n(\lambda_{1,+}(y)-\lambda_{2}(x))})+O(e^{n(\lambda_{3,+}(y)-\lambda_{2}(x))}),\quad
    y\to x,
\end{multline}
uniformly in $n$. Now in the limit when $y\to x$, the last two terms in
\eqref{correlationkernelS2} become exponentially small by virtue of Lemma
\ref{lemmajumpsonR}(b). Then by letting $y \to x$ and using l'H\^opital's rule and find
\begin{align*}K_n(x,x) &= \frac{n(\xi_{1,+}-\xi_{3,+})}{2\pi i}+O(1)\\
&= \frac{n}{2\pi t(1-t)}\sqrt{(\beta_1-x)(x-\alpha_1)}+O(1)
\end{align*}
as $n \to \infty$,
where the last step follows from \eqref{defxi14}. We conclude that
    \begin{equation*}\lim_{n\to\infty} \frac{1}{n} K_n(x,x) = \frac{1}{2\pi
    t(1-t)}\sqrt{(\beta_1^*-x)(x-\alpha_1^*)}, \qquad
    x \in (\alpha_1^*, \beta_1^*).
\end{equation*}
For  $x\in(\alpha_2^*,\beta_2^*)$ we obtain in a similar way that $\frac{1}{n}
K_n(x,x)$ tends to the semi-circle density on the interval $[\alpha_2^*,
\beta_2^*]$. Thus by \eqref{correlationkernel} we have completed the proof of
Theorem \ref{st2ellipses} in the case of large separation (recall that $T=1$).
$\bol$

\section{Steepest descent analysis in the case of critical separation}
\label{sectionsteepestdescentcritical}

In this section we study the non-intersecting Brownian motions in case of
critical separation between the endpoints. We will work under the assumption of
the double scaling regime \eqref{deftau}--\eqref{doublescalinglimit}. This
means that we take the endpoints $a_j,b_j$, $j=1,2$ fixed such that
\begin{equation}\label{criticalseparationbis}
(a_1-a_2)(b_1-b_2) = \left(\sqrt{p_1^*}+\sqrt{p_2^*}\right)^2,
\end{equation} and that we consider the
temperature $T$ to be varying with $n$ as
\begin{equation}\label{doublescalinglimitbis}
    T = T_n = 1 + L n^{-2/3}.
\end{equation} Recall also the assumptions
\eqref{balanced}--\eqref{defp2}.

Before proceeding further, let us recall the main objects needed for the
steepest descent analysis. The points $\alpha_j,\beta_j$, $j=1,2$ are as
defined in \eqref{defalphavarying}--\eqref{defbetavarying}, but now with $T$ in
\eqref{doublescalinglimitbis} not identically equal to 1:
\begin{align} \label{defalphavaryingbis}
    \alpha_j & = \alpha_j(t) =  (1-t)a_j+tb_j-\sqrt{4p_jT t(1-t)},\\
    \label{defbetavaryingbis}
    \beta_j & =  \beta_j(t) = (1-t)a_j+tb_j+\sqrt{4p_j Tt(1-t)},
\end{align}
for $j=1,2$. The limiting values for $n\to\infty$ of these points are denoted
by (cf.\ \eqref{defalpha}--\eqref{defbeta})
\begin{align} \label{defalphastar}
    \alpha_j^* & = \alpha_j^*(t) =  (1-t)a_j+tb_j-\sqrt{4p_j^* t(1-t)},\\
    \label{defbetastar}
    \beta_j^* & =  \beta_j^*(t) = (1-t)a_j+tb_j+\sqrt{4p_j^* t(1-t)},
\end{align}
$j=1,2$. Recall also the $\lam$-functions and $\xi$-functions which are given
by \eqref{deflambda14} and \eqref{defxi14}. We will denote the limiting values
for $n\to\infty$ of these functions by $\lambda_k^*(z)$, $\xi_k^*(z)$,
$k=1,\ldots,4$. For example, the functions $\xi_k^*(z)$ are given by
\begin{equation}\label{defxi14star}
    \begin{array}{lll}
    \xi_1^*(z) & =  \frac{1}{2t(1-t)}\left(-(1-t)a_1+tb_1+((z-\alpha_1^*)(z-\beta_1^*))^{1/2}\right), \\
    \xi_2^*(z) & =  \frac{1}{2t(1-t)}\left(-(1-t)a_2+tb_2+((z-\alpha_2^*)(z-\beta_2^*))^{1/2}\right), \\
    \xi_3^*(z) & =  \frac{1}{2t(1-t)}\left(-(1-t)a_1+tb_1-((z-\alpha_1^*)(z-\beta_1^*))^{1/2}\right), \\
    \xi_4^*(z) & = \frac{1}{2t(1-t)}\left(-(1-t)a_2+tb_2-((z-\alpha_2^*)(z-\beta_2^*))^{1/2}\right).
\end{array}
\end{equation}

Throughout this section we will assume again that the hypothesis
\eqref{hypothesist} holds, but now with strict inequality, i.e.,
\begin{equation}\label{hypothesistbis}
    0 < t < t_{\crit} = \frac{a_1-a_2}{(a_1-a_2) + (b_1-b_2)}.
\end{equation}
The time $t_{\crit}$ is now the time where the two ellipses in Figure
\ref{figcriticalsep} are tangent to each other. The case where $t_{\crit}<t<1$
can be handled in a similar way; cf.\ Remark \ref{remarkeliminationothercase}.
The behavior at the tangent time $t_{\crit}$ itself leads to a multi-critical
situation which is outside the scope of this paper.

\subsection{Modifications in the choice of $x_0$, $\Gamma_1$ and $\Gamma_2$}

The main technical tool that made the analysis in Section
\ref{subsectionasymptoticsjumpslarge} work was the existence of a point
$x_0\in(\beta_2,\alpha_1)$ such that the string of inequalities
\eqref{xiinequalities} holds, cf.\ Lemma \ref{lemmaonxi2}. It will turn out
that in the present context, this can be achieved by the point $x_0^*$ given by
the explicit formula
\begin{equation}\label{defx0bis} x_0^* = \frac{\sqrt{p_1^*}}{\sqrt{p_1^*}+\sqrt{p_2^*}}\frac{\alpha^*_2+\beta^*_2}{2}+
    \frac{\sqrt{p_2^*}}{\sqrt{p_1^*}+\sqrt{p_2^*}}\frac{\alpha^*_1+\beta^*_1}{2}.
\end{equation}
Note that we already encountered the analogue of the point \eqref{defx0bis} in
the proof of Lemma~\ref{lemmaonxi2}, cf.\ \eqref{defx0}. In particular, in the
proof of Lemma~\ref{lemmaonxi2} we derived inequalities \eqref{productx01} and
\eqref{productx02} in case of large separation of the endpoints. Running
through the proof of these inequalities, we find that in case of
critical separation these inequalities become equalities:
\begin{equation*}
(\alpha^*_1-x_0^*)(\beta^*_1-x_0^*) =
\frac{p_1^*}{(\sqrt{p_1^*}+\sqrt{p_2^*})^2}\left((1-t)(a_1-a_2)-t(b_1-b_2)\right)^2
\end{equation*}
and
\begin{equation*}
(x_0^*-\alpha^*_2)(x_0^*-\beta^*_2) =
\frac{p_2^*}{(\sqrt{p_1^*}+\sqrt{p_2^*})^2}\left((1-t)(a_1-a_2)-t(b_1-b_2)\right)^2.
\end{equation*}
These equalities can be restated in the form
\begin{align}\label{productx0crit1}
\sqrt{(\alpha^*_1-x_0^*)(\beta^*_1-x_0^*)} =
\frac{\sqrt{p_1^*}}{\sqrt{p_1^*}+\sqrt{p_2^*}}\left((1-t)(a_1-a_2)-t(b_1-b_2)\right)\\
\label{productx0crit2} \sqrt{(x_0^*-\alpha^*_2)(x_0^*-\beta^*_2)} =
\frac{\sqrt{p_2^*}}{\sqrt{p_1^*}+\sqrt{p_2^*}}\left((1-t)(a_1-a_2)-t(b_1-b_2)\right).
\end{align}
The positivity of the right-hand sides of
\eqref{productx0crit1}--\eqref{productx0crit2} follows from our assumption
\eqref{hypothesistbis}.

In what follows we will also need the identities
\begin{align}\label{sumx01}
\frac{\alpha_1^*+\beta_1^*}{2}-x_0^* &=
\frac{\sqrt{p_1^*}}{\sqrt{p_1^*}+\sqrt{p_2^*}}\left(
(1-t)(a_1-a_2)+t(b_1-b_2)\right),\\
\label{sumx02}
x_0^*-\frac{\alpha_2^*+\beta_2^*}{2} &=
\frac{\sqrt{p_2^*}}{\sqrt{p_1^*}+\sqrt{p_2^*}}\left(
(1-t)(a_1-a_2)+t(b_1-b_2)\right).
\end{align}
These identities follow immediately from the definitions \eqref{defx0bis} and
\eqref{defalphastar}--\eqref{defbetastar}.

\begin{lemma}\label{lemmaonxicrit1}
Under the double scaling regime
\eqref{criticalseparationbis}--\eqref{doublescalinglimitbis}, the point $x_0^*$
defined in \eqref{defx0bis} satisfies the inequalities
\begin{equation} \label{xiinequalitiescrit1}
\xi_2^*(x_0^*) > \xi_3^*(x_0^*),\quad \xi_4^*(x_0^*) > \xi_1^*(x_0^*).
\end{equation}
\end{lemma}

\bewijs. From the definitions \eqref{defxi14star} and keeping track of the
correct branches of the $1/2$ powers, we find
\begin{multline}\label{derivativenegative7}
\xi_4^*(x_0^*)-\xi_1^*(x_0^*)  = \frac{1}{2t(1-t)}(
(1-t)(a_1-a_2)-t(b_1-b_2)\\
- \sqrt{(x_0^*-\alpha_2)(x_0^*-\beta_2)}
 +\sqrt{(x_0^*-\alpha_1)(x_0^*-\beta_1)}).
\end{multline}
By \eqref{productx0crit1} and \eqref{productx0crit2}, this is
\begin{align*}
\frac{1}{2t(1-t)}\left((1-t)(a_1-a_2)-t(b_1-b_2)\right)(1-\frac{\sqrt{p_2^*}}{\sqrt{p_1^*}+\sqrt{p_2^*}}+\frac{\sqrt{p_1^*}}{\sqrt{p_1^*}+\sqrt{p_2^*}})\\
= \frac{1}{2t(1-t)}(
(1-t)(a_1-a_2)-t(b_1-b_2))\frac{2\sqrt{p_1^*}}{\sqrt{p_1^*}+\sqrt{p_2^*}}\
>\ 0.
\end{align*}
Similarly, $\xi_2^*(x_0^*)-\xi_3^*(x_0^*)>0$. $\bol$

Note that the relations \eqref{xiinequalitiescrit1} correspond to the two
outermost inequalities in \eqref{xiinequalities}, evaluated asymptotically for
$n\to\infty$. By continuity, these inequalities then also hold for $n$ finite
but sufficiently large.

In a similar way we would like to have the middle inequality in
\eqref{xiinequalities}, i.e., $\xi_3^*(x_0^*)>\xi_4^*(x_0^*)$.
However, this is not the case, since instead we have equality
$\xi_4^*(x_0^*)-\xi_3^*(x_0^*)=0$. We need a more detailed statement of the
zero behavior.

\begin{lemma}\label{lemmathirdder} We have
\begin{equation}\label{thirdrootbehavior}
\lambda_4^*(z)-\lambda_3^*(z) = c \frac{(z-x_0^*)^3}{3!}+O((z-x_0^*)^4),\quad z\to
x_0^*,
\end{equation}
where $c$ is the positive constant given by
\begin{equation}\label{thirdderivativeexplicit}
c = \frac{2(\sqrt{p_1^*}+\sqrt{p_2^*})^4}{\sqrt{p_1^*p_2^*}}
\left((1-t)(a_1-a_2)-t(b_1-b_2)\right)^{-3}.
\end{equation}
\end{lemma}

\bewijs. Recall from Section \ref{sectionsteepestdescentlargesep} that we
choose the constant $\kappa$ such that
\begin{equation}\label{lambdasvanish}\lambda_4(x_0^*)=\lambda_3(x_0^*)\end{equation}
and hence in particular
\begin{equation*}\lambda_4^*(x_0^*)=\lambda_3^*(x_0^*).\end{equation*}

Now we consider the subsequent derivatives of $\lambda_4^*-\lambda_3^*$ at $x_0^*$.
We start with the first derivative $\xi_4^*(x_0^*)-\xi_3^*(x_0^*)$. Recall that
in the proof of Lemma~\ref{lemmaonxi2} we showed that this derivative is
negative provided there is large separation between the endpoints. The same
proof shows that in case of critical separation, this derivative is
zero. This means that
\begin{equation*}\xi_4^*(x_0^*)=\xi_3^*(x_0^*),\end{equation*} as we already
alluded to in the paragraph before the statement of Lemma \ref{lemmathirdder}.

Consider then the second derivative
\begin{equation}\label{secondderivative}(\xi_4^*)'(x_0^*)-(\xi_3^*)'(x_0^*)\end{equation}
of $\lambda_4^*-\lambda_3^*$ at $x_0^*$. It follows immediately from
\eqref{defxi14star} that
\begin{align*} (\xi_3^*)'(z) &=
    -\frac{1}{2t(1-t)} \left(z-\tfrac{\alpha_1^*+\beta_1^*}{2}\right)((z-\alpha_1^*)(z-\beta_1^*))^{-1/2}\\
(\xi_4^*)'(z) &=
    -\frac{1}{2t(1-t)}\left(z-\tfrac{\alpha_2^*+\beta_2^*}{2}\right)((z-\alpha_2^*)(z-\beta_2^*))^{-1/2}.
\end{align*}
Inserting \eqref{productx0crit1}--\eqref{sumx02}, and keeping track of the
branches of the 1/2 powers, we see that both terms in \eqref{secondderivative}
cancel each other when $z=x_0^*$ and hence \eqref{secondderivative} vanishes.

Consider then the third derivative
\begin{equation}\label{thirdderivative}(\xi_4^*)''(x_0^*)-(\xi_3^*)''(x_0^*)
\end{equation}
of $\lambda_4^*-\lambda_3^*$ at $x_0^*$. We have
\begin{multline*} (\xi_3^*)''(z) = \frac{1}{2t(1-t)} \left(
    \left(z-\tfrac{\alpha_1^*+\beta_1^*}{2}\right)^2
    ((z-\alpha_1^*)(z-\beta_1^*))^{-3/2}\right.\\ \left.-((z-\alpha_1^*)(z-\beta_1^*))^{-1/2}
    \right)\end{multline*}\vspace{-7mm}
    \begin{multline*}
(\xi_4^*)''(z) = \frac{1}{2t(1-t)} \left(
    \left(z-\tfrac{\alpha_2^*+\beta_2^*}{2}\right)^2
    ((z-\alpha_2^*)(z-\beta_2^*))^{-3/2}\right.\\ \left.-((z-\alpha_2^*)(z-\beta_2^*))^{-1/2}
    \right),
\end{multline*}
or equivalently
\begin{multline*} (\xi_3^*)''(z) =
    \frac{1}{2t(1-t)}((z-\alpha_1^*)(z-\beta_1^*))^{-3/2}\times\\ \left(
    \left(z-\tfrac{\alpha_1^*+\beta_1^*}{2}\right)^2-(z-\alpha_1^*)(z-\beta_1^*) \right)
    \end{multline*}\vspace{-7mm}
    \begin{multline*}
    (\xi_4^*)''(z) = \frac{1}{2t(1-t)}((z-\alpha_2^*)(z-\beta_2^*))^{-3/2}\times\\ \left(
    \left(z-\tfrac{\alpha_2^*+\beta_2^*}{2}\right)^2-(z-\alpha_2^*)(z-\beta_2^*) \right),
\end{multline*}
which can be further simplified to
\begin{align*} (\xi_3^*)''(z) &=
    \frac{1}{2t(1-t)}((z-\alpha_1^*)(z-\beta_1^*))^{-3/2}\left(
    \tfrac{\alpha_1^*-\beta_1^*}{2} \right)^2
    \\
    (\xi_4^*)''(z) &= \frac{1}{2t(1-t)}((z-\alpha_2^*)(z-\beta_2^*))^{-3/2}\left(
    \tfrac{\alpha_2^*-\beta_2^*}{2} \right)^2.
\end{align*}
By inserting \eqref{defalphastar}--\eqref{defbetastar} this yields
\begin{align*} (\xi_3^*)''(z) &=
2p_1^*((z-\alpha_1^*)(z-\beta_1^*))^{-3/2}
\\
(\xi_4^*)''(z) &= 2p_2^*((z-\alpha_2^*)(z-\beta_2^*))^{-3/2}.
\end{align*}

Using \eqref{productx0crit1}--\eqref{productx0crit2}, we then find
\begin{eqnarray*}
    & & (\xi_4^*)''(x_0^*) - (\xi_3^*)''(x_0^*) \\ &=&
    2\left(\frac{1}{\sqrt{p_1^*}}+\frac{1}{\sqrt{p_2^*}}\right)
    \left(\sqrt{p_1^*}+\sqrt{p_2^*}\right)^3
    \left((1-t)(a_1-a_2)-t(b_1-b_2)\right)^{-3}\\
    &=& \frac{2(\sqrt{p_1^*}+\sqrt{p_2^*})^4}{\sqrt{p_1^*p_2^*}}
\left((1-t)(a_1-a_2)-t(b_1-b_2)\right)^{-3}.
\end{eqnarray*}
This proves the lemma. $\bol$

From Lemma \ref{lemmathirdder} we see that the behavior of
$\Re(\lambda_4^*(z)-\lambda_3^*(z))$ in the neighborhood of $x_0^*$ depends on the
different sectors in the complex plane: in particular this real part is
negative in the sectors
\begin{align*}
 -\frac{\pi}{2}< \arg(z-x_0^*)<-\frac{\pi}{6}, \qquad \frac{\pi}{6} <\arg(z-x_0^*)<\frac{\pi}{2}
\end{align*}
and it is positive in the sectors
\begin{align*}
-\frac{2\pi}{3} < \arg (z-x_0^*)<-\frac{\pi}{2}, \qquad \frac{\pi}{2} < \arg(z-x_0^*)<\frac{2\pi}{3}.
\end{align*}
We will then choose $\Gamma_1$ so that it lies in the sectors with negative
real part and $\Gamma_2$ so that it lies in the sectors with positive real
part. In particular, note that $\Gamma_1$ and $\Gamma_2$ now both pass
through $x_0^*$. Just as in Section
\ref{sectionsteepestdescentlargesep}, we will choose these curves independent
of $n$.

\subsection{Transformations of the RH problem}

We are now ready to describe the steepest descent analysis in case of
critical separation, assuming the double scaling limit
\eqref{criticalseparationbis}--\eqref{doublescalinglimitbis} with $n$ fixed and
sufficiently large.

\subsubsection{First transformation: Gaussian elimination in the jump matrix}
\label{subsectiontransfo1critical}
 The first transformation $Y \mapsto X$  of the steepest descent analysis is a
Gaussian elimination in the jump matrix. This is done in the same way as in
Section \ref{subsectiontransfo1large}, except that we choose the contours
$\Gamma_1$ and $\Gamma_2$ in a different way. As explained at the end of the
previous subsection, these curves should be chosen fixed (independent of $n$)
and such that they pass through the point $x_{0}^*$ in certain sectors of the
complex plane.

\subsubsection{Second transformation: Normalization at infinity}
\label{subsectiontransfo2critical}

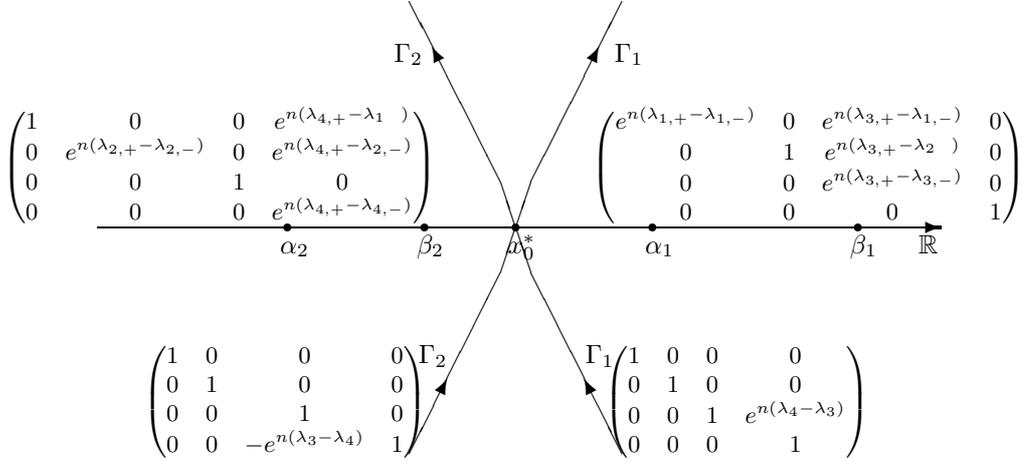
\begin{figure}[t]
\begin{center}
   \setlength{\unitlength}{1truemm}
   \begin{picture}(100,70)(-5,2)
       \put(40,40){\line(-1,-3){2}}
       \put(40,40){\line(-1,3){2}}
       \put(38,34){\line(-1,-2){12}}
       \put(38,46){\line(-1,2){12}}
       \put(24,62){$\Gamma_2$}
       \put(27.2,22){$\Gamma_2$}
       \put(40,40){\line(1,-3){2}}
       \put(40,40){\line(1,3){2}}
       \put(42,34){\line(1,-2){12}}
       \put(42,46){\line(1,2){12}}
       \put(53,62){$\Gamma_1$}
       \put(49.2,22){$\Gamma_1$}
       \put(50,62){\thicklines\vector(1,2){1}}
       \put(30,62){\thicklines\vector(-1,2){1}}
       \put(50,18){\thicklines\vector(-1,2){1}}
       \put(30,18){\thicklines\vector(1,2){1}}

       \put(-15,40){\line(1,0){110}}
       \put(93,36.6){$\mathbb R $}
       \put(95,40){\thicklines\vector(1,0){1}}
       \put(9,36.6){$\alpha_2$}
       \put(27,36.6){$\beta_2$}
       \put(39,36.6){$x_{0}^*$}
       \put(57,36.6){$\alpha_1$}
       \put(84,36.6){$\beta_1$}
       \put(10,40){\thicklines\circle*{1}}
       \put(28,40){\thicklines\circle*{1}}
       \put(40,40){\thicklines\circle*{1}}
       \put(58,40){\thicklines\circle*{1}}
       \put(85,40){\thicklines\circle*{1}}

       \put(51.7,16){$\small{\begin{pmatrix}1&0&0&0\\ 0&1&0&0  \\
       0&0&1&e^{n(\lambda_4-\lambda_3)}\\ 0&0&0&1
       \end{pmatrix}}$}
       \put(-9,16){$\small{\begin{pmatrix}1&0&0&0\\ 0&1&0&0  \\
       0&0&1&0\\ 0&0&-e^{n(\lambda_3-\lambda_4)}&1
       \end{pmatrix}}$}
       \put(-27.5,47){$\small{\begin{pmatrix}1&0&0&e^{n(\lambda_{4,+}-\lambda_1\phantom{{}_{-}})}\\ 0&e^{n(\lambda_{2,+}-\lambda_{2,-})}&0&e^{n(\lambda_{4,+}-\lambda_{2,-})}  \\
       0&0&1&0\\ 0&0&0&e^{n(\lambda_{4,+}-\lambda_{4,-})}
       \end{pmatrix}}$}
       \put(50,47){$\small{\begin{pmatrix}e^{n(\lambda_{1,+}-\lambda_{1,-})}&0&e^{n(\lambda_{3,+}-\lambda_{1,-})}&0\\ 0&1&e^{n(\lambda_{3,+}-\lambda_2\phantom{{}_{-}})}&0  \\
       0&0&e^{n(\lambda_{3,+}-\lambda_{3,-})}&0\\ 0&0&0&1
       \end{pmatrix}}$}

   \end{picture}
   \caption{The jump matrices in the RH problem
   for $T=T(z)$ in the case of critical separation. Compare with Figure \ref{figjumpslambdas}.}
   \label{figjumpslambdascrit}
\end{center}
\end{figure}

The transformation $X \mapsto T$ is the same as
in Section \ref{subsectiontransfo2large}. The new matrix
valued function $T=T(z)$ is normalized at infinity in the sense that
$T(z)=I+O(1/z)$ as $z\to\infty$. The jump matrices in the RH problem
for $T(z)$ are shown in Figure \ref{figjumpslambdascrit}.

\subsubsection{Asymptotic behavior of the jump matrices}
\label{subsectiontransfo3critical}

By Lemmas \ref{lemmaonxicrit1} and \ref{lemmathirdder} and the choice of the
contours $\Gamma_1$ and $\Gamma_2$, we see that the conclusions of Section
\ref{subsectionasymptoticsjumpslarge} all remain valid for $n$ sufficiently
large, except that in a neighborhood of the point $x_{0}^*$ (through which now
both $\Gamma_1$ and $\Gamma_2$ pass), the jump matrices are not exponentially
close to the identity matrix. This means that, except for a small neighborhood
of $x_{0}^*$, the RH problem again asymptotically decouples into two $2\times
2$ problems involving rows and columns $1$, $3$ and $2$, $4$, respectively, in
exactly the same way as in \eqref{jumpTapprox3}--\eqref{jumpTapprox4}.

\begin{figure}[t]
\begin{center}
   \setlength{\unitlength}{1truemm}
   \begin{picture}(100,70)(-5,2)
       \put(40,40){\line(2,3){20}}
       \put(40,40){\line(2,-3){20}}
       \put(40,40){\line(-2,3){20}}
       \put(40,40){\line(-2,-3){20}}
       \put(40,40){\thicklines\circle*{1}}
       \put(38,35){$x_{0}^*$}
       \put(50,55){\thicklines\vector(2,3){.0001}}
       \put(50,25){\thicklines\vector(-2,3){.0001}}
       \put(30,55){\thicklines\vector(-2,3){.0001}}
       \put(30,25){\thicklines\vector(2,3){.0001}}
       \put(48,60){$\Gamma_1$}
       \put(48,20){$\Gamma_1$}
       \put(28,60){$\Gamma_2$}
       \put(30,20){$\Gamma_2$}

       \put(55,60){$\small{\begin{pmatrix}1&e^{n(\lambda_4-\lambda_3)}\\ 0&1 \end{pmatrix}}$}
       \put(54,20){$\small{\begin{pmatrix}1&e^{n(\lambda_4-\lambda_3)}\\ 0&1\end{pmatrix}}$}
       \put(-1,60){$\small{\begin{pmatrix}1&0\\ -e^{n(\lambda_3-\lambda_4)}&1 \end{pmatrix}}$}
       \put(0.5,20){$\small{\begin{pmatrix}1&0\\ -e^{n(\lambda_3-\lambda_4)}&1\end{pmatrix}}$}
      \end{picture}
   \caption{After ignoring the exponentially small entries of the jump matrices  in Figure
   \ref{figjumpslambdascrit}, the jump matrices near $x_0^*$ in the RH problem
   for $T(z)$ have a $2\times 2$ structure as shown in the figure.}
   \label{figlocalRH}
\end{center}
\end{figure}
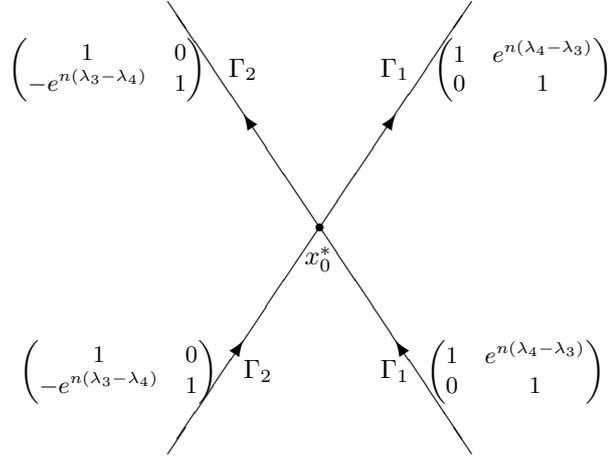

The only place where the above decoupling is not valid is in a neighborhood of
$x_{0}^*$. In fact, since $x_{0}^*$ is away from the intervals
$[\alpha_1,\beta_1]$ and $[\alpha_2,\beta_2]$ for $n$ sufficiently large, the
jump matrices on the real line close to $x_{0}^*$ are all exponentially close
to the identity matrix. Ignoring these jump matrices, the only jump conditions
that remain are those on the curves $\Gamma_1$ and $\Gamma_2$. From
Figure~\ref{figjumpslambdascrit} we see that the latter constitute essentially
a $2\times 2$ RH problem involving rows and columns $3$ and $4$ only. This
leads to the jump matrices shown in Figure~\ref{figlocalRH}.

\subsubsection{Third transformation: Opening of the lenses}
\label{subsectionlenses}

Around the intervals $[\alpha_1,\beta_1]$ and $[\alpha_2,\beta_2]$, we are in
the region where the RH problem decouples into two $2\times 2$ problems. We can
then define a transformation $T \mapsto S$ by opening lenses around
these intervals. This is done in exactly the same way as in Section
\ref{subsectiontransfo3large}. Since these operations have no influence on the
behavior around $x_{0}^*$ (which is our main point of interest here), we omit
a detailed description.

\subsubsection{Model RH problem: Parametrix away from the special points}
\label{subsectionPinftycrit}

We describe now the solution to the model RH problem where we ignore all
exponentially small entries in the jump matrices, and where we stay away from
the special points $\alpha_1,\beta_1,\alpha_2,\beta_2$ and $x_{0}^*$. Since we
are considering the region away from the special point $x_{0}^*$, we are
essentially dealing with two $2\times 2$ matrix valued RH problems. The construction of the
parametrix $P^{(\infty)}(z)$ is then exactly the same as in Section
\ref{subsectionPinftylarge}.

\subsubsection{Local parametrix around $\alpha_1,\beta_1,\alpha_2,\beta_2$}
\label{subsectionlocalparametrixalphabeta}

In small disks around the endpoints of the intervals $[\alpha_1,\beta_1]$ and
$[\alpha_2,\beta_2]$, one can construct local parametrices
$P^{(\textrm{Airy})}(z)$ to the RH problem in terms of Airy functions. The
construction is exactly the same as in Section \ref{subsectionPAirylarge}.

\subsection{Local parametrix around $x_0^*$}
\label{subsectionlocalparametrixx0star}

\subsubsection{Model RH problem associated to the Painlev\'e II equation}
\label{subsectionRHPII}

Our next task is to build a local parametrix in a neighborhood of the special
point $x_{0}^*$. This will be the main technical part of the analysis.
In this subsection we will first
recall the model RH problem associated with the Hastings-McLeod solution of the
Painlev\'e II equation. We describe this RH problem and a few basic
deformations of it, serving to bring it closer to the form we will need.

\begin{figure}[t]
\begin{center}
   \setlength{\unitlength}{1truemm}
   \begin{picture}(100,70)(-5,2)
       \put(40,40){\line(3,2){30}}
       \put(40,40){\line(3,-2){30}}
       \put(40,40){\line(-3,2){30}}
       \put(40,40){\line(-3,-2){30}}
       \put(40,40){\line(0,1){30}}
        \put(40,40){\line(0,-1){30}}
       \put(40,40){\thicklines\circle*{1}}
       \put(38,36){$0$}
       \put(55,50){\thicklines\vector(3,2){.0001}}
       \put(55,30){\thicklines\vector(3,-2){.0001}}
       \put(25,50){\thicklines\vector(-3,2){.0001}}
       \put(25,30){\thicklines\vector(-3,-2){.0001}}
       \put(40,58){\thicklines\vector(0,1){.0001}}
       \put(40,22){\thicklines\vector(0,-1){.0001}}
       \qbezier(43,42)(44,40)(43,38)
       \put(45,40){$\pi/3$}

       \put(61.5,50){$\small{\begin{pmatrix}1&0\\ s_1&1 \end{pmatrix}}$}
       \put(39.5,67){$\small{\begin{pmatrix}1&s_2\\ 0&1 \end{pmatrix}}$}
       \put(4,51){$\small{\begin{pmatrix}1&0\\ s_3&1 \end{pmatrix}}$}
       \put(6.5,29){$\small{\begin{pmatrix}1&s_1\\ 0&1\end{pmatrix}}$}
       \put(39.5,13){$\small{\begin{pmatrix}1&0\\ s_2&1 \end{pmatrix}}$}
       \put(58.5,30){$\small{\begin{pmatrix}1&s_3\\ 0&1\end{pmatrix}}$}
  \end{picture}
   \caption{The jump matrices in the $\zeta$-plane
   in the RH problem for  $\Psi = \Psi(\zeta,s)$ associated
   with a general solution of the Painlev\'e II equation. The Stokes multipliers
   $s_1$, $s_2$, $s_3$  satisfy $s_1 + s_2 + s_3 + s_1s_2s_3 = 0$.}
   \label{figRHPsigen}
\end{center}
\end{figure}
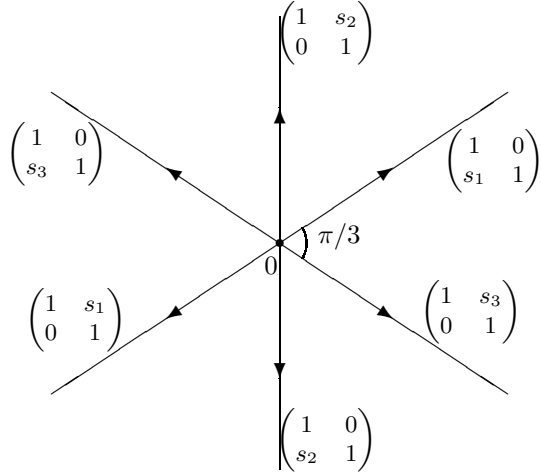

We start with the RH problem for the \lq $\Psi$-functions\rq\ associated to a
general solution of the Painlev\'e II equation \cite{FN,FIKN}. To this end,
consider a matrix function depending on two variables $\Psi = \Psi(\zeta,s)$,
of which the variable $s$ is considered fixed at this moment. The jump
conditions of $\Psi = \Psi(\zeta,s)$ in the $\zeta$-plane are shown in Figure
\ref{figRHPsigen}. The jumps are determined by three numbers
$s_1,s_2,s_3\in\mathbb C $ which are called Stokes multipliers; these may be
any complex numbers satisfying the relation
\begin{equation}\label{stokesmultipliers}
    s_1 + s_2 + s_3 + s_1 s_2 s_3 = 0.
\end{equation}
The $2 \times 2$ matrix valued function $\Psi$ then satisfies the following RH
problem:
\begin{itemize}
\item[(1)] $\Psi(\zeta,s)$ is analytic for $\zeta\in\mathbb C \setminus\bigcup_{k=0}^5\{\arg\zeta = \pi/6+k\pi/3\}$;
\item[(2)] For $\zeta\in\bigcup_{k=0}^5\{\arg\zeta = \pi/6+k\pi/3\}$, $\Psi(\zeta,s)$ has jumps as shown in Figure \ref{figRHPsigen};
\item[(3)] As $\zeta\to\infty$ we have that
\begin{equation}\label{asymptoticHML}\Psi(\zeta,s)e^{i(\frac{4}{3}\zeta^3+s\zeta)\sigma_3} =
I+O(\zeta^{-1}),\end{equation} where $\sigma_3 := \diag(1,-1)$ denotes the
third Pauli matrix.
\end{itemize}
The corresponding Painlev\'e II function is recovered from the RH problem by
the formula \cite{FN,FIKN}
\[ q(s) = \lim_{\zeta \to \infty} 2i \zeta \Psi_{12}(\zeta,s) e^{i(\frac{4}{3}\zeta^3+s\zeta)\sigma_3}. \]

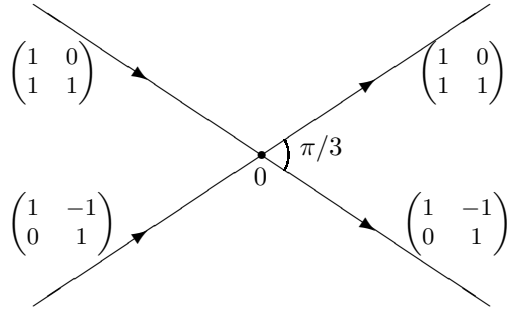
\begin{figure}[t]
\begin{center}
   \setlength{\unitlength}{1truemm}
   \begin{picture}(100,70)(-5,20)
       \put(40,40){\line(3,2){30}}
       \put(40,40){\line(3,-2){30}}
       \put(40,40){\line(-3,2){30}}
       \put(40,40){\line(-3,-2){30}}
       \put(40,40){\thicklines\circle*{1}}
       \put(39,36){$0$}
       \put(55,50){\thicklines\vector(3,2){.0001}}
       \put(55,30){\thicklines\vector(3,-2){.0001}}
       \put(25,50){\thicklines\vector(3,-2){.0001}}
       \put(25,30){\thicklines\vector(3,2){.0001}}
       \qbezier(43,42)(44,40)(43,38)
       \put(45,40){$\pi/3$}

       \put(60.5,50){$\small{\begin{pmatrix}1&0\\ 1&1 \end{pmatrix}}$}
       \put(58.5,30){$\small{\begin{pmatrix}1&-1\\ 0&1\end{pmatrix}}$}
       \put(6.5,50){$\small{\begin{pmatrix}1&0\\ 1&1 \end{pmatrix}}$}
       \put(6.5,30){$\small{\begin{pmatrix}1&-1\\ 0&1\end{pmatrix}}$}
      \end{picture}
   \caption{The jump matrices in the $\zeta$-plane
   in the RH problem for  $\Psi = \Psi(\zeta,s)$ associated
   with the Hastings-McLeod solution of the Painlev\'e II equation. The Stokes
   multipliers are $s_1 = 1$, $s_2 = 0$, $s_3 = -1$.}
   \label{figRHPsi}
\end{center}
\end{figure}

The \emph{Hastings-McLeod solution} to the Painlev\'e II equation corresponds
to the special choice of Stokes multipliers
\[ s_1=1, \qquad s_2=0, \qquad s_3=-1. \]
Since $s_2=0$, we see that the jump on the imaginary axis in Figure
\ref{figRHPsigen} disappears. The jump
matrices of the RH problem then reduce to the situation in Figure
\ref{figRHPsi}. Note that we have reversed the orientation of the two leftmost
rays.

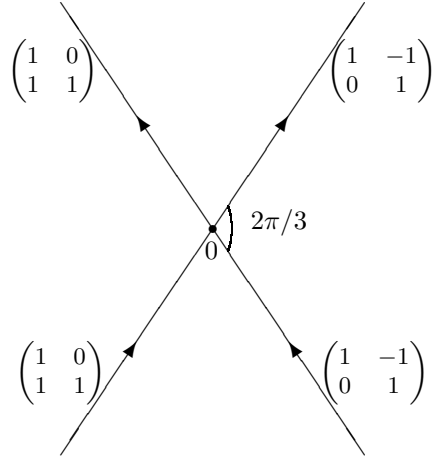
\begin{figure}[t]
\begin{center}
   \setlength{\unitlength}{1truemm}
   \begin{picture}(100,70)(-5,10)
       \put(40,40){\line(2,3){20}}
       \put(40,40){\line(2,-3){20}}
       \put(40,40){\line(-2,3){20}}
       \put(40,40){\line(-2,-3){20}}
       \put(40,40){\thicklines\circle*{1}}
       \put(39,36){$0$}
       \put(50,55){\thicklines\vector(2,3){.0001}}
       \put(50,25){\thicklines\vector(-2,3){.0001}}
       \put(30,55){\thicklines\vector(-2,3){.0001}}
       \put(30,25){\thicklines\vector(2,3){.0001}}
       \qbezier(42,43)(43,40)(42,37)
       \put(45,40){$2\pi/3$}

       \put(55,60){$\small{\begin{pmatrix}1&-1\\ 0&1 \end{pmatrix}}$}
       \put(54,20){$\small{\begin{pmatrix}1&-1\\ 0&1\end{pmatrix}}$}
       \put(13,60){$\small{\begin{pmatrix}1&0\\ 1&1 \end{pmatrix}}$}
       \put(14,20){$\small{\begin{pmatrix}1&0\\ 1&1\end{pmatrix}}$}
  \end{picture}
   \caption{The jump matrices in the $\zeta$-plane
   in the RH problem for $\widetilde{\Psi} = \widetilde{\Psi}(\zeta,s)$.}
   \label{figRHtildePsi}
\end{center}
\end{figure}

Now we apply a rotation over $90$ degrees, i.e., we set
\begin{equation}\label{deftildePsi}
\widetilde{\Psi}(\zeta,s) := \Psi(-i\zeta,s).
\end{equation}
Then $\widetilde{\Psi}$ satisfies the following RH problem:
\begin{itemize}
\item[(1)] $\widetilde{\Psi}(\zeta,s)$ is analytic for
    $\zeta\in \mathbb C \setminus\bigcup_{k=1}^2\{\arg\zeta = \pm k\pi/3\}$;
\item[(2)] On the rays $\arg\zeta=\pm \pi/3$, $\widetilde{\Psi}(\zeta,s)$ has jump matrix $\begin{pmatrix} 1 &
-1 \\ 0 & 1
\end{pmatrix}$,\\
On the rays $\arg\zeta=\pm 2\pi/3$, $\widetilde{\Psi}(\zeta,s)$ has jump matrix
$\begin{pmatrix} 1 & 0 \\ 1 & 1
\end{pmatrix}$, see Figure~\ref{figRHtildePsi};
\item[(3)] As $\zeta\to\infty$ we have that
\begin{equation*}\widetilde{\Psi}(\zeta,s)e^{(-\frac{4}{3}\zeta^3+s\zeta)\sigma_3} =
I+O(\zeta^{-1}).\end{equation*}
\end{itemize}

\begin{figure}[t]
\begin{center}
   \setlength{\unitlength}{1truemm}
   \begin{picture}(100,70)(-5,10)
       \put(40,40){\line(2,3){20}}
       \put(40,40){\line(2,-3){20}}
       \put(40,40){\line(-2,3){20}}
       \put(40,40){\line(-2,-3){20}}
       \put(40,40){\thicklines\circle*{1}}
       \put(39,36){$0$}
       \put(50,55){\thicklines\vector(2,3){.0001}}
       \put(50,25){\thicklines\vector(-2,3){.0001}}
       \put(30,55){\thicklines\vector(-2,3){.0001}}
       \put(30,25){\thicklines\vector(2,3){.0001}}
       \qbezier(42,43)(43,40)(42,37)
       \put(45,40){$2\pi/3$}

       \put(55,60){$\small{\begin{pmatrix}1&e^{\frac{8}{3}\zeta^3-2s\zeta}\\ 0&1 \end{pmatrix}}$}
       \put(54,20){$\small{\begin{pmatrix}1&e^{\frac{8}{3}\zeta^3-2s\zeta}\\ 0&1\end{pmatrix}}$}
       \put(-4,60){$\small{\begin{pmatrix}1&0\\ -e^{-(\frac{8}{3}\zeta^3-2s\zeta)}&1 \end{pmatrix}}$}
       \put(-3,20){$\small{\begin{pmatrix}1&0\\ -e^{-(\frac{8}{3}\zeta^3-2s\zeta)}&1\end{pmatrix}}$}
      \end{picture}
   \caption{The figure shows the jump matrices in the $\zeta$-plane
   in the RH problem for $M = M(\zeta,s)$.}
   \label{figRHM}
\end{center}
\end{figure}
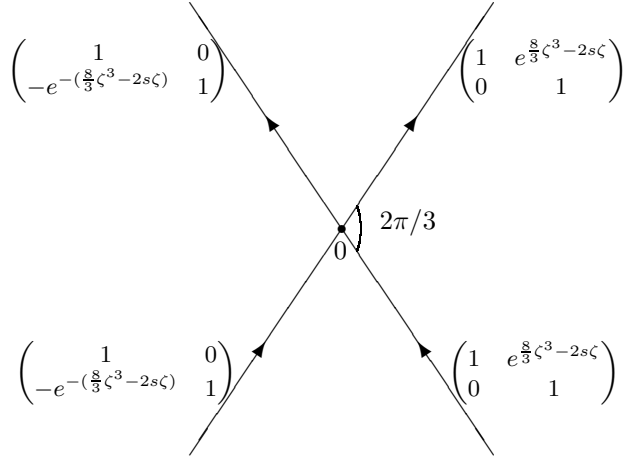

Now we apply one final modification. Define a new matrix function
\begin{equation}\label{defM}
M(\zeta,s) :=
\sigma_3\widetilde{\Psi}(\zeta,s)e^{(-\frac{4}{3}\zeta^3+s\zeta)\sigma_3}\sigma_3.
\end{equation}
Then $M$ satisfies the following RH problem:
\begin{itemize}
\item[(1)] $M(\zeta,s)$ is analytic for $\zeta\in\mathbb C \setminus\bigcup_{k=1}^2\{\arg\zeta = \pm k\pi/3\}$;
\item[(2)] On the rays $\arg\zeta=\pm \pi/3$, $M(\zeta,s)$ has jump matrix $\begin{pmatrix} 1 &
e^{\frac{8}{3}\zeta^3-2s\zeta} \\ 0 & 1
\end{pmatrix}$,\\
On the rays $\arg\zeta=\pm 2\pi/3$, $M(\zeta,s)$
has jump matrix $\begin{pmatrix} 1 & 0 \\
-e^{-(\frac{8}{3}\zeta^3-2s\zeta)} & 1
\end{pmatrix}$, see Figure~\ref{figRHM};
\item[(3)] As $\zeta\to\infty$ we have that
\begin{equation}\label{asymptoticM}
M(\zeta,s) = I+O(\zeta^{-1}).\end{equation}
\end{itemize}

Note in particular that the RH problem for $M=M(\zeta,s)$ is normalized at
infinity in the sense of \eqref{asymptoticM}. Moreover, the $O$-term in
\eqref{asymptoticM} holds uniformly for $s$ in compact subsets of
$\mathbb C \setminus\Pee$ with $\Pee$ the set
of poles of the Hastings-McLeod solution to the Painlev\'e II equation. Of
particular interest for us is the fact \cite{HML} that $\Pee\cap\mathbb R  =
\emptyset$, i.e., there are no poles of the Hastings-Mcleod solution
on the real line. The RH problem for $M=M(\zeta,s)$ will be
used in the construction of a local parametrix near $x_{0}^*$.

\subsubsection{Construction of local parametrix}
\label{subsectionlocalparametrixx0}

Now we construct a local parametrix near the special point $x_{0}^*$. It turns
out that the construction can be done in a similar way as done by Claeys and
Kuijlaars in \cite{Claeys1}; see also \cite{Claeys2}. We also note that the
construction of a local parametrix with $\Psi$-functions associated with
Painlev\'e II played a role \cite{BDJ} and \cite{BI2}.

We are going to construct the local parametrix in the neighborhood
\begin{equation}\label{defUdelta}
    U_{\delta}:=\{z\in\mathbb C |\ |z-x_{0}^*|<\delta\}
    \end{equation}
of $x_0^*$, where the radius $\delta>0$ is fixed but sufficiently small.

Recall that for $n$ sufficiently large the jumps near $x_{0}^*$ are essentially
those of a $2\times 2$ matrix-valued RH problem involving rows and columns $3$, $4$,
with jump matrices shown in Figure
\ref{figlocalRH}. To construct the local parametrix near $x_0^*$, we
construct similar jumps via the model RH problem for $M(\zeta,s)$ in Figure
\ref{figRHM}. To this end, we propose a local parametrix $P^{(x_{0}^*)}(z)$ of
the form
\begin{equation}\label{defPx0}
P^{(x_{0}^*)}(z) := P^{(\infty)}(z)\begin{pmatrix} I & 0\\
0 & M(n^{1/3}f(z),n^{2/3}s_n(z))
\end{pmatrix},\quad  z\in U_{\delta}
\end{equation}
where $P^{(\infty)}(z)$ is the solution to the model RH problem in Section
$\ref{subsectionPinftycrit}$ and where we set
\begin{align}
\label{deffconformal} f(z) &:=\left(\frac{3}{8}(\lambda_{4}^*-\lambda_{3}^*)(z)\right)^{1/3},\\
\label{defsn} s_n(z) &:=
\frac{(\lambda_{4}^*-\lambda_{3}^*)(z)-(\lambda_{4}-\lambda_{3})(z)}{2 f(z)}.
\end{align}

We will show in the next lemma that $f(z)$ is a conformal map mapping the
neighborhood $U_{\delta}$ of $x_{0}^*$ (provided $\delta$ is small enough)
onto a neighborhood of the origin in the $\zeta$-plane. It follows that $s_n(z)$
in \eqref{defsn} is a well-defined analytic function in $U_{\delta}$,
due to pole-zero cancelation at $x_{0}^*$. Indeed,
this follows from the fact that both terms in the numerator of
\eqref{defsn} vanish at $z=x_0^*$, cf.\ \eqref{lambdasvanish}.

\begin{lemma}\label{lemmafz} (Conformal map)
The function $f(z)$ is analytic in a neighborhood  $U_{\delta}$
of $x_0^*$ and satisfies
\begin{equation}\label{expansionf}
    f(z) = \frac{z-x_0^*}{2K((1-t)(a_1-a_2)-t(b_1-b_2))}+O(z-x_{0}^*)^2,
\end{equation}
as $z \to x_0^*$, where the constant $K$ is defined in \eqref{defK}.
\end{lemma}
\bewijs.
This follows from \eqref{deffconformal} and Lemma \ref{lemmathirdder}.
$\bol$

From the above definitions \eqref{defPx0}--\eqref{defsn} we see that if
$\zeta=n^{1/3}f(z)$ and $s=n^{2/3}s_n(z)$, then the exponent occurring in the
jump matrices in Figure \ref{figRHM} reduces to
\begin{align}
\frac{8}{3}\zeta^3-2s\zeta &= \nonumber n\left(\frac{8}{3}f^3(z)-2s_n(z)
f(z)\right) = n(\lambda_{4}-\lambda_{3})(z).
\end{align}
This shows that $M(n^{1/3}f(z),n^{2/3}s_n(z))$ has precisely the required jumps
in the RH problem near $x_0^*$, cf.\ Figure \ref{figlocalRH}, provided the
contours $\Gamma_1$ and $\Gamma_2$ near $x_{0}^*$ are chosen in such a way that
they are mapped by the conformal map $f$ to the straight lines in Figure \ref{figRHM}.
We have indeed the freedom to choose $\Gamma_1$ and $\Gamma_2$ in that way near
$x_0^*$.

A technical issue that remains is showing that the RH problem for\\
$M(n^{1/3}f(z),n^{2/3}s_n(z))$ is solvable. This is equivalent with the fact
that $n^{2/3}s_n(z)$ stays away from $\Pee$, the set of poles of the
Hastings-McLeod solution to the Painlev\'e II equation. This will indeed follow
from the next lemma. Recall that we are working under the double scaling
assumption
\[ T=  1 + Ln^{-2/3}. \]

\begin{lemma} (Asymptotics of $s_n(z)$) We have as $n \to \infty$,
\begin{equation}\label{snasymptotics}
    n^{2/3} s_n(z) =
    \frac{L}{8\sqrt{t(1-t)} f(z)}
     \left[F_1(z) - F_2(z) \right]  + O(n^{-1/3})
     \end{equation}
     where
     \begin{equation} \label{defFj}
        F_j(z) = \sqrt{p_j^*} (\beta_j^*-\alpha_j^*)
        \int_{x_0^*}^z \left[(y-\alpha_j^*)(y-\beta_j^*)\right]^{-1/2} dy,
        \qquad j=1,2, \end{equation}
    and the $O$-term in \eqref{snasymptotics} is uniform for $z \in U_{\delta}$.
\end{lemma}

\bewijs. Since $T = 1 + Ln^{-2/3}$ and $p_j = p_j^* + O(n^{-1})$ it
follows from
\eqref{defalphavaryingbis}--\eqref{defbetavaryingbis} that
\begin{align} \label{alphajbehavior}
    \alpha_j & = \alpha_j^* - L \sqrt{p_j^*t(1-t)} n^{-2/3} + O(n^{-1}), \\
    \beta_j  &= \beta_j^* + L \sqrt{p_j^*t(1-t)} n^{-2/3} + O(n^{-1}),
    \label{betajbehavior}
    \end{align}
    as $n \to \infty$.
Then by \eqref{alphajbehavior} and \eqref{betajbehavior} we have
uniformly for $y \in \overline{U_{\delta}}$,
\begin{multline} \label{sqrtbehavior}
    \left[(y-\alpha_j)(y-\beta_j)\right]^{1/2} - \left[(y-\alpha_j^*)(y-\beta_j^*)\right]^{1/2} \\
    =
    - \frac{L \sqrt{p_j^*t(1-t)} (\beta_j^*-\alpha_j^*)}{2\left[(y-\alpha_j^*)(y-\beta_j^*)\right]^{1/2}} n^{-2/3}
        + O(n^{-1}). \end{multline}

Since for $z \in U_{\delta}$,
\begin{align*}
    (\lambda_{4}-\lambda_{3})(z) &= \int_{x_{0}^*}^{z} (\xi_{4}-\xi_{3})(y)\, dy, \\
    (\lambda_4^*-\lambda_3^*)(z) &= \int_{x_0^*}^z (\xi_4^*-\xi_3^*(y)) \, dy,
\end{align*}
we obtain from \eqref{defxi14}  that
\begin{multline}\label{intrepr1}
    (\lambda_{4}-\lambda_{3})(z) -
    (\lambda_4^*-\lambda_3^*)(z) \\
    = \frac{1}{2t(1-t)} \int_{x_{0}^*}^{z}
      \left(  \left[(y-\alpha_1)(y-\beta_1)\right]^{1/2} -
        \left[(y-\alpha_1^*)(y-\beta_1^*)\right]^{1/2} \right) dy \\
    -  \frac{1}{2t(1-t)} \int_{x_{0}^*}^{z}
   \left( \left[(y-\alpha_2)(y-\beta_2)\right]^{1/2}
    -     \left[(y-\alpha_2^*)(y-\beta_2^*)\right]^{1/2} \right)dy
\end{multline}
which by \eqref{sqrtbehavior} and \eqref{defFj} leads to
\[
    (\lambda_{4}-\lambda_{3})(z) -
    (\lambda_4^*-\lambda_3^*)(z)
    = -\frac{L}{4\sqrt{t(1-t)}}
    \left[F_1(z) - F_2(z) \right] n^{-2/3} + O(n^{-1}(z-x_0^*)) \]
and the $O$ term is uniform for $z \in U_{\delta}$.

The lemma now follows because of the definition of $s_n(z)$ in \eqref{defsn}
and the fact that $f(z)$ has a simple zero at $z=x_0^*$. $\bol$

\begin{corollary}
\begin{enumerate}
\item[\rm (a)] We have as $z \to x_0^*$ and $n \to \infty$,
\begin{equation} \label{sntosstrong}
    n^{2/3} s_n(z) = s + O(z-x_0^*) + O(n^{-1/3}) \end{equation}
where the real constant $s$ is defined in \eqref{defs0}.
\item[\rm (b)] If $\delta$ is sufficiently small, then
there is a compact subset $\mathcal K$ of $\mathbb C \setminus \mathcal P$,
where $\mathcal P$ denotes the set of poles of the Hastings-McLeod solution
of Painlev\'e II, such that
\[ n^{2/3} s_n(z) \in \mathcal K \]
for every $z \in U_{\delta}$ and for all large enough $n$.
\end{enumerate}
\end{corollary}
\bewijs.
(a) From \eqref{snasymptotics}
it follows that
\begin{equation} \label{snproof1}
    n^{2/3} s_n(z) =
    \frac{L}{8\sqrt{t(1-t)}} \lim_{z \to x_0^*}
     \frac{F_1(z) - F_2(z)}{f(z)}   + O(z-x_0^*) + O(n^{-1/3}).
     \end{equation}
By \eqref{defFj}, \eqref{defalphastar}--\eqref{defbetastar},
and \eqref{productx0crit1}--\eqref{productx0crit2}
we have for $j=1,2$,
\begin{align} \nonumber
\lim_{z \to x_0^*} \frac{F_j(z)}{z-x_0^*} & = \sqrt{p_j^*} (\beta_j^*-\alpha_j^*)
    \left[(x_0^* - \alpha_j^*)(x_0^*-\beta_j^*)\right]^{-1/2}  \\
       & = \label{snproof2}
        (-1)^j  \frac{4  \sqrt{p_j^*t(1-t)} (\sqrt{p_1^*} + \sqrt{p_2^*})}
            {(1-t)(a_1-a_2) - t(b_1-b_2)},
                \end{align}
and by \eqref{expansionf} we have
\begin{align} \label{snproof3}
    \lim_{z \to x_0^*} \frac{z-x_0^*}{f(z)}
        & = 2K ((1-t)(a_1-a_2)-t(b_1-b_2)).
\end{align}
Multiplying \eqref{snproof2} and \eqref{snproof3} we find
\begin{equation} \label{snproof4}
    \frac{L}{8\sqrt{t(1-t)}} \lim_{z \to x_0^*}
     \frac{F_1(z) - F_2(z)}{f(z)}
     = - LK \left(\sqrt{p_1^*} + \sqrt{p_2^*}\right)^2 = s
     \end{equation}
where for the last equality we used the definition of $K$, see \eqref{defK},
and the relation \eqref{defs0} between $L$ and $s$. Now \eqref{sntosstrong}
follows from \eqref{snproof1} and \eqref{snproof4}.

(b) The equation \eqref{sntosstrong} implies that there exist
constants $K_1$ and $K_2$, independent of $n$ and $z$, such that
for $z \in U_{\delta}$ and $n$ large enough,
\[ |n^{2/3} s_n(z) - s| \leq K_1 |z-x_0^*| + K_2 n^{-1/3}. \]
Since $s$ is real, and since the Hastings-McLeod solution has no poles on the
real line, part (b) follows. $\bol$

Now we can check that the local parametrix $P^{(x_0^*)}$ defined by \eqref{defPx0}
satisfies the following `local RH problem'

\begin{itemize}
\item[(1)] $P^{(x_0^*)}$ is analytic in $U_{\delta} \setminus (\Gamma_1 \cup
\Gamma_2)$;
\item[(2)] On $\Gamma_1$ and $\Gamma_2$, $P^{(x_0^*)}$ has jumps with jump matrices
of the form $\begin{pmatrix} I_2 & 0 \\ 0 & * \end{pmatrix}$ where the $2\times
2$ matrices indicated by $*$ are given in Figure \ref{figlocalRH};
\item[(3)] As $n \to \infty$,
\begin{equation}\label{compatiblearoundx0}
    P^{(x_{0}^*)}(z) = P^{(\infty)}(z)(I+O(n^{-1/3})),
\end{equation}
uniformly for $z$ on the circle $\partial U_{\delta}$.
\end{itemize}
The matching condition (3) follows from the definition \eqref{defPx0}
and the fact that
$M(\zeta,s) = I + O(1/\zeta)$ as $\zeta \to \infty$
uniformly for all $s$ of the form $s=n^{2/3} s_n(z)$.


\subsection{Fourth transformation and completion of the proof of Theorem \ref{st2ellipses}}

\subsubsection{Fourth transformation}
Using the global parametrix $P^{(\infty)}$ and the
local parametrices $P^{(\textrm{Airy})}$ and $P^{(x_0^*)}$,
we define the final transformation $S \mapsto R$ by
\begin{equation}\label{defR}
R(z) = \left\{
    \begin{array}{ll}
    S(z)(P^{(x_{0}^*)})^{-1}(z),& \quad \textrm{inside the disk } U_{\delta} \textrm{ around }x_0^*, \\
    S(z)(P^{(\textrm{Airy})})^{-1}(z),& \quad \textrm{in the disks around } \alpha_1,\beta_1,\alpha_2,\beta_2, \\
    S(z)(P^{(\infty)})^{-1}(z),& \quad \textrm{elsewhere}.
\end{array}
\right.
\end{equation}

From the constructions  in this section it follows that $R$ satisfies
the following RH problem
\begin{itemize}
\item[(1)] $R(z)$ is analytic in $\mathbb C \setminus \Sigma_R$
where $\Sigma_R$ consists of the real line $\mathbb R$, the contours $\Gamma_1$
and $\Gamma_2$ outside $U_{\delta}$, the lips of the lenses outside the disks,
and the boundaries of the disks around $\alpha_1$, $\beta_1$,  $\alpha_2$,
$\beta_2$ and $x_0^*$;
\item[(2)] $R$ has jumps $R_+ = R_- J_R$ on $\Sigma_R$,
that satisfy
\begin{align*}
    J_R(z) & = I + O(n^{-1/3}), \qquad \text{on the boundary of $U_{\delta}$}, \\
    J_R(z) & = I + O(1/n), \qquad \text{on the boundaries of the other disks}, \\
    J_R(z) & = I + O(e^{-cn|z|}), \qquad \text{on the other parts of $\Sigma_R$},
    \end{align*}
    for some constant $c >0$.
\item[(3)] $R(z) = I + O(1/z)$ as $z \to \infty$.
\end{itemize}

Then, as before we may conclude that
\begin{equation} \label{eq:Rasymptoticscritical}
    R(z) = I_4 + O\left(\frac{1}{n^{1/3}(|z|+1)}\right)
    \end{equation}
    as $n \to \infty$, uniformly for $z$ in the complex plane.
    This completes the RH steepest descent analysis.

We will explicitly compute this $O(n^{-1/3})$ term in Section
\ref{sectionPIIasymptotics}.

\subsubsection{Proof of Theorem \ref{st2ellipses} in the critical case}

We have now completed the steepest descent analysis of the RH
problem under the double scaling regime
\eqref{criticalseparationbis}--\eqref{doublescalinglimitbis}. In a completely
similar way as in Section~\ref{subsectionproof2ellipses}
we can now use the result \eqref{eq:Rasymptoticscritical} to
show that the non-intersecting Brownian particles are asymptotically
distributed on the two intervals according to two semicircle laws.
This completes the proof of Theorem~\ref{st2ellipses}.

\section{Proofs of Theorems \ref{stp2asymptotics} and \ref{stp2asymptoticsB}}
\label{sectionPIIasymptotics}

\subsection{Formulas for recurrence coefficients}
In this section we investigate the large $n$ asymptotics of the recurrence
coefficients of the multiple Hermite polynomials (Sections
\ref{subsectionintorecHerm} and \ref{subsectionintroPIIasymptotics}) under the
assumptions in Section~\ref{sectionsteepestdescentcritical}. This will lead to
the proofs of Theorems \ref{stp2asymptotics} and \ref{stp2asymptoticsB}.

We have not yet proved Proposition \ref{proprecurrence4x4}.
This will be done in Section \ref{sectionlaxpair}, but anticipating
this result, we call the combinations
\begin{equation} \label{cijcji}
    c_{i,j} c_{j,i}, \qquad 1 \leq i < j \leq 4,
    \end{equation}
the off-diagonal recurrence coefficients, and
\begin{equation} \label{cijcjkcik}
    (1-t)a_i+tb_k-\frac{c_{i,j} c_{j,k+2}}{c_{i,k+2}}, \qquad 1 \leq i,j,k \leq 2,
    \quad i\neq j,
    \end{equation}
the diagonal recurrence coefficients.

Recall that the $c_{i,j}$ are the entries of the matrix $Y_1$ in
the expansion
\[ Y(z) = \left(I + \frac{Y_1}{z} + O\left(\frac{1}{z^2}\right)\right)
    \diag(z^{n_1}, z^{n_2}, z^{-n_1}, z^{-n_2})
    \]
as $z \to \infty$.

Following the transformations \eqref{defX3}, \eqref{defT}, \eqref{defS3} and
\eqref{defR} in the steepest descent analysis, we have the following
representation for  $z$ in the region between the contours $\Gamma_1$
and $\Gamma_2$
\begin{align}
    Y(z) &= L^n T(z) G^{-n}(z) L^{-n}
    \label{tracebacksteepestdescent}
    = L^n R(z)P^{(\infty)}(z)G^{-n}(z) L^{-n}.
\end{align}

The following lemma is easy to check.
\begin{lemma} \label{lemmaY1}
We have
\begin{equation} \label{Y1formula}
    Y_1 = L^n \left(G_1 + P_1^{(\infty)} + R_1 \right) L^{-n},
    \end{equation}
where $G_1$, $P_1^{(\infty)}$ and $R_1$ are matrices  from the expansions
as $z \to \infty$,
\begin{align} \label{asymptoticG}
    G^{-n}(z) & =
    \left(I+\frac{G_1}{z}+O\left(\frac{1}{z^2}\right)\right)\diag(z^{n_1},z^{n_2},z^{-n_1},z^{-n_2}), \\
    \label{asymptoticPinfty}
    P^{(\infty)}(z) & = I+\frac{P_1^{(\infty)}}{z}+O\left(\frac{1}{z^2}\right),\\
    \label{asymptoticR0}
    R(z) & = I+\frac{R_1}{z}+ O\left(\frac{1}{z^2}\right).
\end{align}
\end{lemma}

We are only interested in the combinations \eqref{cijcji} and \eqref{cijcjkcik}
of entries of $Y_1$. Since $L$ is a diagonal matrix, the factors $L^n$ and $L^{-n}$
in \eqref{Y1formula} will not play a role for these combinations. Also, since
$G_1$ is a diagonal matrix (which is clear from \eqref{asymptoticG}, since
$G(z)$ is diagonal), this does not play a role either.
Therefore we have for $i < j$,
\begin{equation} \label{cijcjiformula}
    c_{i,j} c_{j,i} =
    \left(P_1^{(\infty)} + R_1 \right)_{i,j}
    \left(P_1^{(\infty)} + R_1 \right)_{j,i} \end{equation}
and for distinct $i, j, k$,
\begin{equation} \label{cijcjkcikformula}
    \frac{c_{i,j} c_{j,k}}{c_{i,k}}
    = \frac{\left(P_1^{(\infty)} + R_1 \right)_{i,j}
     \left(P_1^{(\infty)} + R_1 \right)_{j,k}}{\left(P_1^{(\infty)} + R_1 \right)_{i,k}}
     \end{equation}

In what follows we evaluate $P_1^{(\infty)}$ and $R_1$.

The evaluation of $P_1^{(\infty)}$ is straightforward.
\begin{lemma} \label{lem:Pinfty1}
We have
\begin{equation}
\label{Pinfty1formula} P_1^{(\infty)} =
    i \sqrt{Tt(1-t)} \begin{pmatrix}
    0 & 0 & \sqrt{p_1} & 0 \\
    0 & 0 & 0 & \sqrt{p_2} \\
    -\sqrt{p_1} & 0 & 0 & 0 \\
    0 & -\sqrt{p_2} & 0 & 0
\end{pmatrix}.
\end{equation}
\end{lemma}

\bewijs.
Recall that $P^{(\infty)}(z)$ is defined by means of \eqref{Pinfty} and
\eqref{defBD}. For $z \to \infty$ we have the expansions
\begin{align*}
    \gamma_j(z) & =  \left(1 - \frac{\beta_j-\alpha_j}{z-\alpha_j}\right)^{1/4} =
    1 - \frac{1}{4}\frac{\beta_j-\alpha_j}{z}+O\left(\frac{1}{z^2}\right), \\
    \gamma_j^{-1}(z) & = 1+\frac{1}{4} \frac{\beta_j-\alpha_j}{z}+O\left(\frac{1}{z^2}\right).
\end{align*}
Since $\beta_j - \alpha_j = 4 \sqrt{p_j Tt(1-t)}$, see
\eqref{defalphavaryingbis}--\eqref{defbetavaryingbis},
we obtain  \eqref{Pinfty1formula}  from \eqref{Pinfty} and \eqref{asymptoticPinfty}. $\bol$

\subsection{Second term in the expansion of the jump matrix for $R(z)$}

We use the clockwise orientation on the circle $\partial U_{\delta}$ around $x_{0}^*$.
Thus the outside of the circle is the
$+$-side and the inside of the circle is the $-$-side. The matrix valued function
$R(z)$ in \eqref{defR} then satisfies the jump condition
\begin{equation}\label{jumpR}
    R_+(z) = R_-(z)\Delta(z),
    \qquad z \in \partial U_{\delta}
    \end{equation}
on $\partial U_{\delta}$, with jump matrix $\Delta(z)$ given by
\begin{align}
    \nonumber\Delta(z) &= P^{(x_{0}^*)}(z)(P^{(\infty)})^{-1}(z)\\ &=
    P^{(\infty)}(z)\begin{pmatrix} I & 0\\
     0 & M(n^{1/3}f(z),n^{2/3}s_n(z))
    \end{pmatrix}(P^{(\infty)})^{-1}(z),
    \label{defDelta}
\end{align}
see \eqref{defPx0}. To prepare for the evaluation of $R_1$,
we first compute the second term in the expansion of the
jump matrix $\Delta(z)$ as $n \to \infty$.

\begin{lemma}\label{lemmaY2} (Asymptotics of jump matrix for $R(z)$)
The jump matrix $\Delta(z)$ in \eqref{jumpR}, \eqref{defDelta} has the
asymptotics
\begin{equation}\label{expansionDelta}
\Delta(z) = I+\Delta^{(1)}(z)n^{-1/3}+O(n^{-2/3}),\qquad  z \in \partial U_{\delta},
\end{equation}
as $n \to \infty$, where
\begin{equation}
\label{defDelta1} \Delta^{(1)}(z) = \frac{1}{2f(z)} P^{(\infty)}(z)
\begin{pmatrix}0&0&0&0\\ 0&0&0&0\\
0&0& u(n^{2/3} s_n(z))& -q(n^{2/3} s_n(z)) \\
0&0& q(n^{2/3} s_n(z))& -u(n^{2/3} s_n(z))\end{pmatrix}(P^{(\infty)})^{-1}(z),
\end{equation}
uniformly for $z$ on the circle $\partial U_{\delta}$. Here $q(s)$ denotes the
Hastings-McLeod solution to Painlev\'e II, $u(s)$ denotes the Hamiltonian
\begin{equation} \label{defus}
    u(s) := (q'(s))^2-sq^2(s)-q^4(s),
    \end{equation}
 and $f(z)$ is the conformal map in
\eqref{deffconformal}.
\end{lemma}

\bewijs. First, we calculate the large $n$ asymptotics of the matrix
$M(n^{1/3}f(z),$\\ $n^{2/3}s_n(z))$ which occurs in the expression
\eqref{defDelta}. Recall from Section \ref{subsectionRHPII} that $M(\zeta,s)$
is defined as an easy modification of the model RH matrix $\Psi(\zeta,s)$
associated to the Hastings-McLeod solution to the Painlev\'e II equation. It is
known (see e.g.\ \cite{FIKN}) that the asymptotic expansion
\eqref{asymptoticHML} of $\Psi(\zeta,s)$ in powers of $\zeta^{-1}$ can be
refined to
\begin{equation*}
    \Psi(\zeta,s) = \left(I+\frac{1}{2i\zeta}\begin{pmatrix}
    u(s)&q(s)\\-q(s)&-u(s)\end{pmatrix}
    +O\left(\frac{1}{\zeta^2}\right)\right)e^{-i(\frac{4}{3}\zeta^3+s\zeta)\sigma_3}, \qquad
    \text{as } \zeta\to\infty,
\end{equation*}
where $q(s)$ denotes the Hastings-McLeod solution of Painlev\'e II and $u(s)$
denotes the Hamiltonian \eqref{defus}. Following the sequence of transformations
$\Psi\mapsto\widetilde{\Psi}\mapsto M$ in Section \ref{subsectionRHPII} we
obtain
\begin{equation*}
    M(\zeta,s) = I+\frac{1}{2\zeta}\begin{pmatrix}
    u(s)&-q(s)\\q(s)&-u(s)\end{pmatrix}
    +O\left(\frac{1}{\zeta^2}\right), \qquad \text{as } \zeta \to \infty.
\end{equation*}
It follows that
\begin{multline}\label{expansionM}
M(n^{1/3}f(z),n^{2/3}s_n(z)) = I+\frac{1}{2f(z)}\begin{pmatrix}
u(n^{2/3}s_n(z))&-q(n^{2/3}s_n(z))\\q(n^{2/3}s_n(z))&-u(n^{2/3}s_n(z))\end{pmatrix}
    n^{-1/3}
\\+O(n^{-2/3}),
\end{multline}
as $n\to\infty$, uniformly for all $z$ on the circle $\partial U_{\delta}$.

Then \eqref{defDelta1} follows from \eqref{expansionM}, \eqref{defDelta},
and \eqref{expansionDelta}.$\bol$

Note that $\Delta^{(1)}(z)$ depends on $n$, as is obvious
from the appearance of $n^{2/3} s_n(z)$ in
\eqref{defDelta1}. Also $P^{(\infty)}$ depends on $n$.
As $n \to \infty$, the $n$-dependent entries have
limits, and therefore we can still obtain
from the expansion \eqref{expansionDelta} of $\Delta(z)$
a similar expansion of the RH matrix $R(z)$
\begin{equation} \label{Rforlargen}
    R(z) = I +  R^{(1)}(z) n^{-1/3} + O(n^{-2/3}) \qquad \text{as } n \to \infty.
\end{equation}

The coefficient $R^{(1)}(z)$, in this expansion can be found by inserting
\eqref{expansionDelta} and \eqref{Rforlargen} into the jump condition $R_+(z) =
R_-(z)\Delta(z)$ and collecting terms of order $n^{-1/3}$. This leads to the
following additive RH problem for $R^{(1)}$:
\begin{itemize}
\item[(1)] $R^{(1)}(z)$ is analytic in $\mathbb C \setminus\partial U_{\delta}$,
\item[(2)] $R^{(1)}_{+}(z) = R^{(1)}_{-}(z) + \Delta^{(1)}(z)$ for $z\in\partial
U_{\delta}$,
\item[(3)] $R^{(1)}(z) = O(1/z)$ as $z\to\infty$.
\end{itemize}
Note that $R^{(1)}$ also depends on $n$.

The jump matrix \eqref{defDelta1}  has a simple pole at $z=x_{0}^*$.
As in \cite{Claeys2,KMVV}, the RH problem for $R^{(1)}$ then has the explicit solution
\begin{equation}\label{solutionR1}
    R^{(1)}(z) = \left\{\begin{array}{ll}
    \ds \frac{1}{z-x_{0}^*}\Res(\Delta^{(1)}(z),z=x_{0}^*)-\Delta^{(1)}(z),& z\in U_{\delta}\\[5pt]
    \ds \frac{1}{z-x_{0}^*}\Res(\Delta^{(1)}(z),z=x_{0}^*),& z\in \mathbb C \setminus
    U_{\delta}.
 \end{array} \right.
\end{equation}

As in other works (see e.g. \cite{DKMVZ1,KMVV}) the expansions of $R(z)$ as
$z \to \infty$ and $n \to \infty$ commute with each other. It thus follows from
\eqref{Rforlargen} and \eqref{solutionR1} that
\begin{align} \nonumber
     R_1 & = \lim_{z \to \infty} z R^{(1)}(z) n^{-1/3} + O(n^{-2/3}) \\
         & = \Res(\Delta^{(1)}(z),z=x_{0}^*) n^{-1/3}  + O(n^{-2/3}).
         \label{R1andDelta1}
         \end{align}

\subsection{Evaluation of $R_1$}

In view of \eqref{R1andDelta1} the evaluation of $R_1$ comes down
to the determination of the residue of $\Delta^{(1)}$ at $z=x_0^*$.
The result is the following.

\begin{lemma} \label{lem:R1evaluation}
    We have
    \begin{equation} \label{R1blocks}
        R_1 = K \begin{pmatrix} A & B \\ C & D \end{pmatrix} n^{-1/3} + O(n^{-2/3}),
        \qquad \text{ as } n \to \infty,
    \end{equation}
    where $K$ is the constant defined in \eqref{defK}
    and where the $2 \times 2$ blocks $A$, $B$, $C$, $D$ are
    given by
    \begin{align} \label{defA}
    A & = t (b_1 - b_2) \begin{pmatrix} -u(s) & - q(s) \\ q(s) & u(s) \end{pmatrix} \\
    \label{defB}
    B & = i\sqrt{t(1-t)(a_1-a_2)(b_1-b_2)} \begin{pmatrix} -u(s) & q(s) \\ q(s) & -u(s) \end{pmatrix} \\
    \label{defC}
    C & = -i\sqrt{t(1-t)(a_1-a_2)(b_1-b_2)} \begin{pmatrix} u(s) & q(s) \\ q(s) & u(s) \end{pmatrix} \\
    \label{defD}
    D & = (1-t)(a_1-a_2) \begin{pmatrix} u(s) & -q(s) \\ q(s) & -u(s) \end{pmatrix}
    \end{align}
    with the value  $s$ given by \eqref{defs0}.
\end{lemma}

\bewijs.
From  \eqref{expansionf} we obtain
\[
    \lim_{z \to x_0^*} \frac{z-x_0^*}{2f(z)} = K ((1-t)(a_1-a_2)-t(b_1-b_2)),
\]
and from \eqref{sntosstrong} we have
\[ q(n^{2/3}s_n(x_0^*)) = q(s) + O(n^{-1/3}), \qquad u(n^{2/3} s_n(x_0^*)) = u(s) + O(n^{-1/3}) \]
as $n \to \infty$, so that
\eqref{defDelta1}  yields
\begin{multline} \label{Delta1formula}
   \Res(\Delta^{(1)}(z),z=x_{0}^*)
   = K ((1-t)(a_1-a_2)-t(b_1-b_2)) \\
   \times
    P^{(\infty)}(x_0^*)
        \begin{pmatrix} 0 & 0 & 0 & 0 \\ 0 & 0 & 0 & 0 \\
        0 & 0 & u(s) & -q(s) \\
        0 & 0 & q(s) & -u(s) \end{pmatrix}
        (P^{(\infty)})^{-1}(x_0^*)
        + O(n^{-1/3}).
        \end{multline}

To evaluate $P^{(\infty)}(x_0^*)$ we observe that by \eqref{Pinfty}
we need the entries $\gamma_j(x_0^*) \pm \gamma_j(x_0^*)$, for $j=1,2$.
We compute the squares of these
\begin{align} \nonumber
    (\gamma_j(x_{0}^*) \pm \gamma_j^{-1}(x_{0}^*))^2 &=
    \gamma_j^{2}(x_{0}^*) \pm 2 +\gamma_j^{-2}(x_{0}^*)\\
     &= \pm 2 + \sqrt{\frac{x_0^* -\beta_j}{x_0^*-\alpha_j}} +
    \sqrt{\frac{x_0^*-\alpha_j}{x_0^*-\beta_j}} \nonumber \\
     & = \label{Pinftyentries}
    \pm 2+ (-1)^j \frac{2x_0^* - \alpha_j-\beta_j}{\sqrt{(x_0^* - \alpha_j)(x_0^*-\beta_j)}},
    \qquad j=1,2,
\end{align}
where we used \eqref{defBD}. We emphasize that we use $\sqrt{\cdot}$ to denote
the positive square root of a positive real number. Now recall that $\alpha_j$
and $\beta_j$ are varying with $n$, and satisfy $\alpha_j = \alpha_j^* +
O(n^{-2/3})$, $\beta_j = \beta_j^* + O(n^{-2/3})$. For the starred quantities
we can use the identities \eqref{productx0crit1}--\eqref{sumx02}, so that
\eqref{Pinftyentries} yields
\begin{align*}
    (\gamma_j(x_{0}^*) \pm \gamma_j^{-1}(x_{0}^*))^2 &=
    \pm 2 + 2\frac{(1-t)(a_1-a_2) + t(b_1-b_2)}{(1-t)(a_1-a_2)-t(b_1-b_2)} +O(n^{-2/3}).
    \end{align*}
Since $0 < \gamma_2(x_0^*) < 1 < \gamma_1(x_0^*)$ (as can easily be
seen from \eqref{defBD} and the fact that $\alpha_2 < \beta_2 < x_0^* < \alpha_1 < \beta_1$),
we obtain after simple calculations
\begin{align}
\frac{\gamma_j(x_{0}^*) +\gamma_j^{-1}(x_{0}^*)}{2}
\label{asymptoticsBD1} &=
\frac{\sqrt{(1-t)(a_1-a_2)}}{\sqrt{(1-t)(a_1-a_2)-t(b_1-b_2)}}+O(n^{-2/3}), \\
\label{asymptoticsBD2}
\frac{\gamma_j(x_{0}^*) -\gamma_j^{-1}(x_{0}^*)}{2} &= (-1)^{j-1}
\frac{\sqrt{t(b_1-b_2)}}{\sqrt{(1-t)(a_1-a_2)-t(b_1-b_2)}}+O(n^{-2/3}),
\end{align}
for $j=1,2$.

Putting \eqref{asymptoticsBD1} and \eqref{asymptoticsBD2} into \eqref{Pinfty}
we find
\begin{multline}\label{Pinftyx0}
P^{(\infty)}(x_{0}^*) =
\frac{1}{\sqrt{(1-t)(a_1-a_2)-t(b_1-b_2)}}\times \\
\begin{pmatrix}\sqrt{(1-t)(a_1-a_2)} &0& -i\sqrt{t(b_1-b_2)} & 0\\
0& \sqrt{(1-t)(a_1-a_2)} &0& i\sqrt{t(b_1-b_2)}\\
i\sqrt{t(b_1-b_2)} &0& \sqrt{(1-t)(a_1-a_2)} &0\\
0&-i\sqrt{t(b_1-b_2)} &0& \sqrt{(1-t)(a_1-a_2)}
 \end{pmatrix}\\+O(n^{-2/3}).
\end{multline}
Using the fact that $P^{(\infty)}(x_0^*)$ has a $2\times 2$ block structure
with blocks having determinant one, we then also obtain
\begin{multline}\label{Pinftyx0inverse}
(P^{(\infty)})^{-1}(x_{0}^*) =
\frac{1}{\sqrt{(1-t)(a_1-a_2)-t(b_1-b_2)}}\times \\
\begin{pmatrix}\sqrt{(1-t)(a_1-a_2)} &0& i\sqrt{t(b_1-b_2)} & 0\\
0& \sqrt{(1-t)(a_1-a_2)} &0& -i\sqrt{t(b_1-b_2)}\\
-i\sqrt{t(b_1-b_2)} &0& \sqrt{(1-t)(a_1-a_2)} &0\\
0&i\sqrt{t(b_1-b_2)} &0& \sqrt{(1-t)(a_1-a_2)}
 \end{pmatrix}\\+O(n^{-2/3}).
\end{multline}
Then we insert \eqref{Pinftyx0} and \eqref{Pinftyx0inverse} into
\eqref{Delta1formula} and after simple calculations we obtain
\eqref{R1blocks}--\eqref{defD}, see also \eqref{R1andDelta1}. $\bol$

\subsection{Proofs of Theorems \ref{stp2asymptotics} and \ref{stp2asymptoticsB}}

Having Lemmas \ref{lem:Pinfty1} and \ref{lem:R1evaluation} it is now
straightforward to compute the asymptotic behavior of the recurrence
coefficients from \eqref{cijcjiformula} and \eqref{cijcjkcikformula}.

As an example, let us compute the expression $\frac{c_{1,2}c_{2,4}}{c_{1,4}}$
in \eqref{examplerecurrenceintro2}. By \eqref{cijcjkcikformula} and from the
fact that $(P_1^{(\infty)})_{i,j} = 0$ whenever $i+j$ is odd, see
\eqref{Pinfty1formula}, we have
\begin{align} \label{examplecalculation1}
    \frac{c_{1,2}c_{2,4}}{c_{1,4}} &=
    \frac{(R_1)_{1,2}\left((P_1^{(\infty)})_{2,4}+(R_1)_{2,4}\right)}{(R_1)_{1,4}}
    \end{align}
By \eqref{R1blocks} and \eqref{defA}--\eqref{defB} we have
\begin{align} \nonumber
    \frac{(R_1)_{1,2}}{(R_1)_{1,4}} & =  \frac{A_{1,2}}{B_{1,2}} + O(n^{-1/3})
    = -\frac{t(b_1-b_2)}{i\sqrt{t(1-t)(a_1-a_2)(b_1-b_2)}} + O(n^{-1/3}) \\
    \label{examplecalculation2}
    & = - t \sqrt{\frac{b_1-b_2}{a_1-a_2}} \frac{1}{i \sqrt{t(1-t)}} + O(n^{-1/3}).
    \end{align}
Furthermore, by \eqref{Pinfty1formula} and \eqref{R1blocks},
\[ (P_1^{(\infty)})_{2,4} + (R_1)_{2,4}
    = i \sqrt{p_2} \sqrt{Tt(1-t)} + O(n^{-1/3}), \]
which since $p_2 = p_2^* + O(n^{-1})$ and $T = 1 + O(n^{-1/3})$,
reduces to
\[ (P_1^{(\infty)})_{2,4} + (R_1)_{2,4}
    = i \sqrt{p_2^*} \sqrt{t(1-t)} + O(n^{-1/3}). \]
Using this and \eqref{examplecalculation2} in \eqref{examplecalculation1},
we obtain
\begin{align} \label{examplecalculation3}
    \frac{c_{1,2}c_{2,4}}{c_{1,4}} &= - t \sqrt{p_2^*\frac{b_1-b_2}{a_1-a_2}} + O(n^{-1/3}),
    \end{align}
which proves the second formula in Theorem \ref{stp2asymptoticsB}. The
other formulas in Theorems \ref{stp2asymptotics} and \ref{stp2asymptoticsB}
are established in a similar way. $\bol$

\begin{remark}\label{remarkuseinvolutionsymmetry} (The case where $t_{\crit}<t<1$)
Recall that the above derivations have all been made under the assumption that
$0<t<t_{\crit}$. The case where $t_{\crit}<t<1$ can be handled by means of
similar calculations; see also Remark \ref{remarkeliminationothercase}.
Alternatively, one can immediately reduce this to the case where
$0<t<t_{\crit}$ by virtue of the involution symmetry in Corollary
\ref{lemmainvolution2}. After a straightforward calculation, one sees then that
the expansion in Lemma \ref{lem:R1evaluation}, and hence the conclusions of
Theorems~\ref{stp2asymptotics} and \ref{stp2asymptoticsB}, remain valid for
$t_{\crit}<t<1$ as well.
\end{remark}

\section{Proofs of Propositions \ref{proprecurrence4x4} and \ref{propconnectionsrecurrencecoeff4x4}}
\label{sectionlaxpair}

In this final section, we prove Propositions \ref{proprecurrence4x4} and
\ref{propconnectionsrecurrencecoeff4x4}. We start by considering recurrence
relations for general multiple orthogonal polynomials. Next we specialize these
results to the multiple Hermite case (i.e., the case of Gaussian weight
functions \eqref{gaussianw1}--\eqref{gaussianw2}). An important tool will be
the Lax pair satisfied by the multiple Hermite polynomials. The compatibility
condition of this Lax pair allows us to derive more detailed information for
the multiple Hermite case. For example, we show how this condition implies
certain scalar product relations between row and column vectors in $Y_1$; these
turn out to be essentially given by the numbers of Brownian particles
associated with the different starting or ending points.

This section is organized as follows. Section \ref{subsectionrecmop} discusses
recurrence relations for multiple orthogonal polynomials in a general setting.
The rest of the section is devoted to the multiple Hermite case. In Section
\ref{subsectiondiffeq} we discuss a differential equation satisfied by multiple
Hermite polynomials. Section \ref{subsectionlaxpair} deals with the
compatibility condition of the Lax pair for multiple Hermite polynomials.
Sections \ref{subsectionscalarproducts} and \ref{subsectionspectralcurve} deal
with the induced scalar product relations. Finally, Section
\ref{subsectionlaxpairfurtherconsequences} completes the proof of Propositions
\ref{proprecurrence4x4} and \ref{propconnectionsrecurrencecoeff4x4}.

\subsection{Recurrence relations for general multiple orthogonal polynomials}
\label{subsectionrecmop}

We start by discussing the recurrence relations for general  multiple
orthogonal polynomials, leading in particular to the proof of the expression
for the off-diagonal recurrence coefficients as given in Proposition
\ref{proprecurrence4x4}.

It is a classical result that orthogonal polynomials on the real line satisfy a
three-term recurrence relation. In this subsection, we discuss the $p+q+1$ term
recurrence relations for the multiple orthogonal polynomials of
Definition~\ref{defmop}.

As in the classical case $p=q=1$ \cite{Dei}, it turns out that the recurrence
relations can be conveniently derived from the Riemann-Hilbert problem in
Section \ref{subsectionintroRH}. To this end, we introduce the next two terms
in the $z\to\infty$ asymptotics of the matrix $Y(z) = Y_{\vecn,\vecm}(z)$ in
\eqref{asymptoticconditionY0}:
\begin{equation}\label{asymptoticconditionY2}
Y(z) = \left(I+\frac{Y_1}{z}+\frac{Y_2}{z^2}+
O\left(\frac{1}{z^3}\right)\right)\diag(z^{n_1},\ldots,z^{n_p},z^{-m_1},\ldots,z^{-m_q}).\end{equation}
Note that the matrices $Y_1 = (Y_{1})_{\vecn,\vecm}$ and $Y_2 =
(Y_{2})_{\vecn,\vecm}$ are constant (independent of $z$).

Our goal is to find recurrence relations between the RH matrices with
multi-indices $\vecn+\vece_k,\vecm+\vece_l$ and $\vecn,\vecm$, respectively.

\begin{convention} \label{conventionek}
Above we use $\vece_k$ and $\vece_l$ as row vectors. In some of the matrix
calculations that follow, we also use these standard basis vectors, but then we
see them as column vectors, as is usual in matrix algebra. For example $\vece_j
\vece_k^T$ where $\cdot^T$ is the transpose, denotes a matrix whose only
nonzero entry is a one at the $(j,k)$th position, while $\vece_j^T \vece_k =
\delta_{j,k}$ is the scalar product of the two vectors. So we introduce the
convention that $\vece_k$ is a row vector when used as a multi-index. In matrix
computations it is considered as a column vector. We trust that this will not
lead to any confusion.
\end{convention}

To derive the recurrence relations, assume some fixed $k\in\{1,\ldots,p\}$ and
$l\in\{1,\ldots,q\}$. From the fact that the jump matrix of
$Y_{\vecn,\vecm}(z)$ on the real line is independent of $\vecn,\vecm$ (cf.\
\eqref{defjumpmatrix0}, \eqref{defWblock0}), it follows that the matrix valued
function
\begin{equation}\label{leadstorecurrence1}
U_{\vecn,\vecm}^{k,l}(z) := Y_{\vecn+\vece_k,\vecm+\vece_l}(z)
Y_{\vecn,\vecm}^{-1}(z)
\end{equation}
is entire. From \eqref{leadstorecurrence1} and the normalizations at infinity
in \eqref{asymptoticconditionY2}, it follows that for $z\to\infty$ we have
\begin{multline}\label{leadstorecurrence1bis}
U_{\vecn,\vecm}^{k,l}(z) = \left(I+\frac{1}{z}
(Y_1)_{\vecn+\vece_k,\vecm+\vece_l}+\cdots\right)\left(z \vece_k \vece_k^T +
\sum_{j\not\in\{k,p+l\}} \vece_j \vece_j^T \right) \times
\\
\left(I+\frac{1}{z} (Y_1)_{\vecn,\vecm}+\cdots\right)^{-1}+O(1/z).
\end{multline}
Here the middle factor in \eqref{leadstorecurrence1bis} is just the identity
matrix with its $k$th and $(p+l)$th diagonal entries replaced by $z$ and $0$,
respectively. By Liouville's theorem it then follows that
$U_{\vecn,\vecm}^{k,l}$ is a polynomial and so by \eqref{leadstorecurrence1bis}
\begin{equation}\label{forwardcoefficientmatrix}
U_{\vecn,\vecm}^{k,l}(z) = z \vece_k \vece_k^T + \sum_{j\not\in\{k,p+l\}}
\vece_j \vece_j^T + (Y_1)_{\vecn+\vece_k,\vecm+\vece_l} \vece_k \vece_k^T -
\vece_k \vece_k^T (Y_1)_{\vecn,\vecm}. \end{equation}
Summarizing, we obtain the following proposition.

\begin{proposition}\label{propforwardmatrixrecurrence} (Forward matrix recurrence relations)
Let $k\in\{1,\ldots,p\}$ and $l\in\{1,\ldots,q\}$. Then we have the matrix
recurrence relation
\begin{equation}\label{forwardmatrixrecurrence}
Y_{\vecn+\vece_k,\vecm+\vece_l}(z) =
U_{\vecn,\vecm}^{k,l}(z)Y_{\vecn,\vecm}(z),
\end{equation} with $U_{\vecn,\vecm}^{k,l}(z)$ given by \eqref{forwardcoefficientmatrix}.
\end{proposition}

Evaluating the different entries of the matrix recurrence relation
\eqref{forwardmatrixrecurrence} leads to a set of scalar recurrence relations
for the MOP. In fact, since the multiplication with the matrix
$U_{\vecn,\vecm}^{k,l}$ in \eqref{forwardmatrixrecurrence} is performed on the
left, it easily follows that the $k$th component functions $A_k(x)$ of the MOP,
$k=1,\ldots,p$ satisfy all the same recurrence relations. So the recurrence
relations for the individual components $A_1(x),\ldots,A_p(x)$ can be
conveniently stacked into vector recurrence relations for the vectors
$\vecA(x)$.

Let us illustrate this for $p=q=2$, $k=l=1$. Then the recurrence relation
\eqref{forwardmatrixrecurrence}, \eqref{forwardcoefficientmatrix} becomes
\begin{equation}\label{examplerecurrence}
Y_{\vecn+\vece_1,\vecm+\vece_1} = \begin{pmatrix} z+\tilde{c}_{1,1}-c_{1,1} & -c_{1,2} & -c_{1,3} & -c_{1,4}\\
\tilde{c}_{2,1} & 1 & 0 & 0 \\ \tilde{c}_{3,1} & 0 & 0 & 0 \\
\tilde{c}_{4,1} & 0 & 0 & 1
\end{pmatrix}Y_{\vecn,\vecm},
\end{equation}
where we use $c_{i,j}$ and $\til{c}_{i,j}$ to denote the entries of
$(Y_1)_{\vecn,\vecm}$ and $(Y_1)_{\vecn+\vece_1,\vecm+\vece_1}$, respectively.
Evaluating the first row of this expression
and using \eqref{RHmatrix4x40} leads to the vector recurrence relation
\begin{multline}\label{examplerecurrence2}
\vecA^{(II,1)}_{\vecn+\vece_1+\vece_1,\vecm+\vece_1} =
(z+\tilde{c}_{1,1}-c_{1,1}) \vecA^{(II,1)}_{\vecn+\vece_1,\vecm}
-c_{1,2}\vecA^{(II,2)}_{\vecn+\vece_2,\vecm}\\
-c_{1,3}(-2\pi i)\vecA^{(I,1)}_{\vecn,\vecm-\vece_1} -c_{1,4}(-2\pi
i)\vecA^{(I,2)}_{\vecn,\vecm-\vece_2}.
\end{multline}
This is the desired five-term recurrence relation. Note that the factors $-2\pi
i$ are due to the diagonal matrix $D$ in Theorem \ref{stsolutionRH}.

Instead of the first row, one could also evaluate one of the other rows of
\eqref{examplerecurrence}. Doing this leads to an additional set of (mostly
trivial) relations.

The relation \eqref{examplerecurrence2} has the undesirable feature that it
involves vectors of MOP with different types of normalizations. To express
everything in terms of one type of normalization, for example the $(II,1)$
normalization, we need the \emph{transition numbers} between different types of
MOP normalizations.

\begin{definition}(Transition numbers) \label{deftransitionnumber}
Let $\vecn,\vecm$ with $|\vecn|=|\vecm| +1 $ be such that the multiple
orthogonal polynomials exist. We define the \emph{transition number} $t_{\vecn,
\vecm}^{(II,k\to I,l)}\in\mathbb C $ to be the nonzero constant so that the vector
relation
\begin{equation}\label{defeqtransitionnumber}
\vecA_{\vecn,\vecm}^{(II,k)} = t_{\vecn, \vecm}^{(II,k\to I,l)}
\vecA_{\vecn,\vecm}^{(I,l)}
\end{equation}
holds. In a similar way we define the transition numbers
$t_{\vecn,\vecm}^{(II,k\to II,\til{k})}$, $t_{\vecn,\vecm}^{(I,l\to
I,\til{l})}$, and $t_{\vecn,\vecm}^{(I,l\to II,k)}$.
\end{definition}

The transition numbers are contained in the matrix $Y_1$ as follows.

\begin{proposition} \label{proptransitionnumbers} (Interpretation of $Y_1$ via transition numbers)
Let $\vecn$, $\vecm$ with $|\vecn|=|\vecm|$ be such that the solvability
condition holds. Let $D$ be the diagonal matrix in Theorem \ref{stsolutionRH}.
Then the off-diagonal entries of the matrix $D^{-1}(Y_1)_{\vecn,\vecm}D$ can be
expressed as transition numbers between different types of normalizations of
MOP. More precisely, the entries of $D^{-1}(Y_1)_{\vecn,\vecm}D$ are given as
follows:
\begin{itemize}
\item If $1\leq k,\til{k}\leq p$ with $k\neq \til{k}$ then the $(k,\til{k})$ entry is $t_{\vecn+\vece_k,\vecm}^{(II,k\to
II,\til{k})}$,
\item If $1\leq k\leq p$ and $1\leq l\leq q$ then the $(k,p+l)$ entry is $t_{\vecn+\vece_k,\vecm}^{(II,k\to
I,l)}$,
\item If $1\leq l\leq q$ and $1\leq k\leq p$ then the $(p+l,k)$ entry is $t_{\vecn,\vecm-\vece_l}^{(I,l\to
II,k)}$,
\item If $1\leq l,\til{l}\leq q$ with $l\neq\til{l}$ then the $(p+l,p+\til{l})$ entry is $t_{\vecn,\vecm-\vece_l}^{(I,l\to
I,\til{l})}$.
\end{itemize}
\end{proposition}

\bewijs. This follows easily from \eqref{asymptoticconditionY2},
\eqref{RHmatrix4x40}, and from the definition of the normalizations. The result
is most easily seen by evaluating \eqref{asymptoticconditionY2} column by
column; we do not provide the details here. $\bol$

Let us illustrate Proposition \ref{proptransitionnumbers} for the case where
$p=q=2$. In this case the proposition asserts that
\begin{equation}\label{transitionnumbers4x4} (Y_1)_{\vecn,\vecm} =
D\left(\begin{array}{cccc}
* & t_{\vecn+\vece_1,\vecm}^{(II,1\to II,2)} & t_{\vecn+\vece_1,\vecm}^{(II,1\to I,1)} & t_{\vecn+\vece_1,\vecm}^{(II,1\to I,2)} \\
t_{\vecn+\vece_2,\vecm}^{(II,2\to II,1)} & * & t_{\vecn+\vece_2,\vecm}^{(II,2\to I,1)} & t_{\vecn+\vece_2,\vecm}^{(II,2\to I,2)} \\
t_{\vecn,\vecm-\vece_1}^{(I,1\to II,1)} & t_{\vecn,\vecm-\vece_1}^{(I,1\to II,2)} & * & t_{\vecn,\vecm-\vece_1}^{(I,1\to I,2)} \\
t_{\vecn,\vecm-\vece_2}^{(I,2\to II,1)} & t_{\vecn,\vecm-\vece_2}^{(I,2\to
II,2)} & t_{\vecn,\vecm-\vece_2}^{(I,2\to I,1)} & *
\end{array}\right)D^{-1},
\end{equation}
where the diagonal entries denoted with $*$ are unspecified, and where the
diagonal matrix $D = \diag(1,1,-2\pi i, -2\pi i)$.

We use the information in Proposition \ref{proptransitionnumbers} to obtain
more meaningful expressions for the recurrence coefficients. To this end,
consider again the five-term recurrence relation \eqref{examplerecurrence2}.
This relation yields a connection between several MOP with different
multi-indices and different types of normalizations. To express everything in
terms of type $(II,1)$ normalizations it suffices to multiply with the
appropriate transition numbers from \eqref{transitionnumbers4x4}. Then
\eqref{examplerecurrence2} transforms into the new relation
\begin{multline}\label{examplerecurrence3}
\vecA^{(II,1)}_{\vecn+\vece_1+\vece_1,\vecm+\vece_1} = (z+
\tilde{c}_{1,1}-c_{1,1}) \vecA^{(II,1)}_{\vecn+\vece_1,\vecm}
- c_{1,2}c_{2,1}\vecA^{(II,1)}_{\vecn+\vece_2,\vecm}\\
- c_{1,3}c_{3,1}\vecA^{(II,1)}_{\vecn,\vecm-\vece_1} -
c_{1,4}c_{4,1}\vecA^{(II,1)}_{\vecn,\vecm-\vece_2},
\end{multline}
where we recall that $c_{i,j}$, $\til{c}_{i,j}$ denote the $(i,j)$th entry of
$(Y_1)_{\vecn,\vecm}$, $(Y_1)_{\vecn+\vece_1,\vecm+\vece_1}$, respectively.
Note in particular that the factors $-2\pi i$ in \eqref{examplerecurrence2}
have disappeared in \eqref{examplerecurrence3}

The relation \eqref{examplerecurrence3} was derived for the special case where
$p=q=2$ and $k=l=1$ in the recurrence relations. In a similar way, one can find
the recurrence relation for general $p,q,k,l$.

\begin{proposition}\label{propforwardrecurrence} (Forward recurrence relations)
For any $k\in\{1,\ldots,p\}$ and $l\in\{1,\ldots,q\}$, we have the $p+q+1$ term
recurrence relation
\begin{multline}\label{forwardrecurrence}
\vecA^{(II,k)}_{\vecn+\vece_k+\vece_k,\vecm+\vece_l} = (z+
\tilde{c}_{k,k}-c_{k,k}) \vecA^{(II,k)}_{\vecn+\vece_k,\vecm}
\\-\sum_{\til{k}\neq k}
c_{k,\til{k}}c_{\til{k},k}\vecA^{(II,k)}_{\vecn+\vece_{\til{k}},\vecm}
-\sum_{\til{l}}
c_{k,p+\til{l}}c_{p+\til{l},k}\vecA^{(II,k)}_{\vecn,\vecm-\vece_{\til{l}}}
\end{multline}
where we use $c_{i,j}$ and $\til{c}_{i,j}$ to denote the entries of
$(Y_1)_{\vecn,\vecm}$ and $(Y_1)_{\vecn+\vece_k,\vecm+\vece_l}$, respectively.
In the above sums, it is assumed that $\til{k}$ runs from $1$ to $p$ while
$\til{l}$ runs from $1$ to $q$.
\end{proposition}

Note that Proposition \ref{propforwardrecurrence} implies the recurrence
relations in Proposition \ref{proprecurrence4x4}, except for the explicit form
of the first term in the right-hand side of each of
\eqref{examplerecurrenceintro1}--\eqref{examplerecurrenceintro4}, which will be
established in Section \ref{subsubsectiondiagrec}. We note that Proposition
\ref{propforwardrecurrence} is valid for general weight functions $w_{1,k}$ and
$w_{2,l}$, not only for Gaussian weights.

\begin{definition}\label{defreccoef} (Off-diagonal and diagonal recurrence coefficients)
We will refer to the coefficients  of the form $c_{k,l}c_{l,k}$ in
\eqref{forwardrecurrence} as the \emph{off-diagonal recurrence coefficients}.
The coefficient $c_{k,k}-\tilde{c}_{k,k}$ in \eqref{forwardrecurrence} will be
called the \emph{diagonal recurrence coefficient}.
\end{definition}

We use the terminology \lq off-diagonal\rq\ and \lq diagonal\rq\ recurrence
coefficients in analogy with the case of classical orthogonal polynomials
$p=q=1$. Indeed, in the latter case these recurrence coefficients correspond
precisely to the off-diagonal and diagonal entries of the Jacobi matrix.

We collect the off-diagonal recurrence coefficients of Definition
\ref{defreccoef} in an upper triangular matrix
\begin{equation}\label{hadamardproduct}
H = \begin{pmatrix}
0 & c_{1,2}c_{2,1} & \ldots  & \ldots & c_{1,p+q}c_{p+q,1}\\
\vdots & 0 & \ddots  & & \vdots\\
\vdots &  & \ddots  & 0 & c_{p+q-1,p+q}c_{p+q,p+q-1}\\
0 & \ldots & \ldots  & \ldots & 0
\end{pmatrix}.
\end{equation}
The matrix $H = (H)_{\vecn,\vecm}$ is the Hadamard product (entry-wise product)
of the strictly upper triangular part of $Y_1 = (Y_1)_{\vecn,\vecm}$ with  the
transpose of the strictly lower triangular part of $Y_1$.
For example, when $p=q=2$ we have
\begin{equation}\label{hadamardproductexample}
H = \begin{pmatrix}
0 & c_{1,2}c_{2,1} & c_{1,3}c_{3,1} & c_{1,4}c_{4,1}\\
0 & 0 & c_{2,3}c_{3,2} & c_{2,4}c_{4,2}\\
0 & 0 & 0 & c_{3,4}c_{4,3}\\
0 & 0 & 0 & 0
\end{pmatrix}.
\end{equation}

The diagonal recurrence coefficient $c_{k,k}-\tilde{c}_{k,k}$ in
\eqref{forwardrecurrence} is equal to a difference of two entries, one from
$(Y_1)_{\vecn+\vece_k,\vecm+\vece_l}$ and the other from $(Y_1)_{\vecn,\vecm}$.
It is possible however to express this quantity completely in terms of the
multi-indices $\vecn,\vecm$. The resulting expression is then somewhat more
complicated, involving entries of both the matrices $Y_1 = (Y_1)_{\vecn,\vecm}$
and $Y_2 = (Y_2)_{\vecn,\vecm}$.

\begin{lemma}\label{lemmadiagonalcoeff}(Diagonal recurrence coefficients)
The diagonal recurrence coefficient $c_{k,k}-\tilde{c}_{k,k}$ in
\eqref{forwardrecurrence} can be expressed as
\begin{align}\nonumber c_{k,k}-\tilde{c}_{k,k} &= (Y_1)_{k,k}+\frac{(Y_2)_{k,p+l}-(Y_1^2)_{k,p+l}}{(Y_1)_{k,p+l}} \\
\label{reccoefdiagonaltype} &= \frac{(Y_2)_{k,p+l}-\sum_{\til{k}\neq
k}c_{k,\til{k}}c_{\til{k},p+l}-\sum_{\til{l}}c_{k,p+\til{l}}c_{p+\til{l},p+l}}{c_{k,p+l}},
\end{align}
where in \eqref{reccoefdiagonaltype} the index $\til{k}$ runs from $1$ to $p$
while $\til{l}$ runs from $1$ to $q$.
\end{lemma}

The proof of \eqref{reccoefdiagonaltype} is similar to the proof in the case of
classical orthogonal polynomials, see \cite[Section 3.2]{Dei}; we omit the
details. A more concise expression in the case of Gaussian weight functions
will be derived in Section \ref{subsubsectiondiagrec}.

In addition to the forward recurrence relations in Propositions
\ref{propforwardmatrixrecurrence} and \ref{propforwardrecurrence}, one can also
run these recurrence relations in backward order. To this end, consider the
matrix function
\begin{equation*}\label{leadstorecurrencenew}
\tilde{U}_{\vecn,\vecm}^{k,l}(z) :=
Y_{\vecn,\vecm}(z)Y_{\vecn+\vece_k,\vecm+\vece_l}^{-1}(z).
\end{equation*}
(Note that this is the inverse of the matrix $U_{\vecn,\vecm}^{k,l}(z)$ in
\eqref{leadstorecurrence1}.) By copying the approach above one finds:

\begin{proposition}\label{propbackwardmatrixrecurrence} (Backward matrix recurrence relations)
Let $k\in\{1,\ldots,p\}$ and $l\in\{1,\ldots,q\}$. Then we have
\begin{equation}\label{backwardmatrixrecurrence}
Y_{\vecn,\vecm}(z) =
\til{U}_{\vecn,\vecm}^{k,l}(z)Y_{\vecn+\vece_k,\vecm+\vece_l}(z),
\end{equation} where
\begin{equation}\label{backwardcoefficientmatrix}
\til{U}_{\vecn,\vecm}^{k,l}(z) = z \vece_{p+l} \vece_{p+l}^T +
\sum_{j\not\in\{k,p+l\}} \vece_j \vece_j^T + (Y_1)_{\vecn,\vecm} \vece_{p+l}
\vece_{p+l}^T - \vece_{p+l} \vece_{p+l}^T
(Y_1)_{\vecn+\vece_{k},\vecm+\vece_l}.
\end{equation} Moreover, the matrix $\til{U}_{\vecn,\vecm}^{k,l}(z)$
is the inverse of the matrix $U_{\vecn,\vecm}^{k,l}(z)$ in
\eqref{forwardcoefficientmatrix}.
\end{proposition}


\begin{proposition}\label{propbackwardrecurrence} (Backward recurrence relations)
For any $k\in\{1,\ldots,p\}$ and $l\in\{1,\ldots,q\}$, we have the $p+q+1$ term
recurrence relation
\begin{multline}\label{backwardrecurrence}
\vecA^{(I,l)}_{\vecn,\vecm-\vece_l} = (z+ c_{p+l,p+l}-\til{c}_{p+l,p+l})
\vecA^{(I,l)}_{\vecn+\vece_k,\vecm}
-\sum_{\til{k}} \til{c}_{p+l,\til{k}}\til{c}_{\til{k},p+l}\vecA^{(I,l)}_{\vecn+\vece_k+\vece_{\til{k}},\vecm+\vece_l}\\
-\sum_{\til{l}\neq l}
\til{c}_{p+l,p+\til{l}}\til{c}_{p+\til{l},p+l}\vecA^{(I,l)}_{\vecn+\vece_k,\vecm+\vece_l-\vece_{\til{l}}}
\end{multline}
where $c_{i,j},\til{c}_{i,j}$ denote the entries of $(Y_1)_{\vecn,\vecm}$ and
$(Y_1)_{\vecn+\vece_k,\vecm+\vece_l}$, respectively.
\end{proposition}

For example, when $p=q=2$, $k=l=1$ the above relations become
\begin{equation*}\label{examplebackwardrecurrence}
Y_{\vecn,\vecm}(z) = \begin{pmatrix} 0 & 0 & c_{1,3} & 0\\
0 & 1 & c_{2,3} & 0 \\ -\tilde{c}_{3,1} & -\tilde{c}_{3,2} & z-\tilde{c}_{3,3}+c_{3,3} & -\tilde{c}_{3,4} \\
0 & 0 & c_{4,3} & 1
\end{pmatrix}Y_{\vecn+\vece_1,\vecm+\vece_1}(z),
\end{equation*}
and
\begin{multline*}\label{backwardrecurrence}
\vecA^{(I,1)}_{\vecn,\vecm-\vece_1} = (z+ c_{3,3}-\til{c}_{3,3})
\vecA^{(I,1)}_{\vecn+\vece_1,\vecm}
- {\til{c}_{3,1}\til{c}_{1,3}}\vecA^{(I,1)}_{\vecn+\vece_1+\vece_{1},\vecm+\vece_1}\\
-
{\til{c}_{3,2}\til{c}_{2,3}}\vecA^{(I,1)}_{\vecn+\vece_1+\vece_{2},\vecm+\vece_1}
-
{\til{c}_{3,4}\til{c}_{4,3}}\vecA^{(I,1)}_{\vecn+\vece_1,\vecm+\vece_1-\vece_{2}}.
\end{multline*}

\subsection{Differential equation for multiple Hermite polynomials}
\label{subsectiondiffeq}

In this subsection we derive a differential equation for multiple Hermite
polynomials, i.e., assuming that the weights are Gaussian as in
\eqref{gaussianw1}--\eqref{gaussianw2}. The differential equation was already
described in \cite{Daems}.

For general $p$ and $q$ we consider the weights
\begin{align*}
w_{1,k} & =  e^{- \frac{N}{2t}(x^2-2a_kx)},\quad k=1,\ldots,p,\\
w_{2,l} & =  e^{- \frac{N}{2(1-t)}(x^2-2b_lx)},\quad l=1,\ldots,q,
\end{align*}
with $N>0$ a certain constant.

We introduce a modification of the RH matrix $Y(z)$.
\begin{definition} \label{defPsifunction}
Define the matrix function
\begin{equation}\label{defPsi}
\Psi(z) := Y(z)\diag(f_1(z),\ldots,f_{p+q}(z)),
\end{equation}
where we used the following functions $f_j$, $j=1,\ldots,p+q$:
\begin{equation} \label{deffs}
\left\{\begin{array}{ll} f_k(z) = e^{-\frac{N}{2t(1-t)}(z^2-2(1-t)a_k z)},\quad & k=1,\ldots,p,\\
f_{p+l}(z) = e^{-\frac{N}{(1-t)}b_l z},\quad & l=1,\ldots,q.
\end{array}\right.
\end{equation}
\end{definition}

The reason for defining the matrix function $\Psi$ in Definition
\ref{defPsifunction} is that it has a constant jump (independent of $x$) along
the real line. More precisely, it satisfies the following Riemann-Hilbert
problem:
\begin{itemize}
\item[(1)] $\Psi$ is analytic on $\mathbb C \setminus\mathbb R $;
\item[(2)] For $x\in\mathbb R $, it holds that \begin{equation*} \Psi_{+}(x) = \Psi_{-}(x)
\begin{pmatrix} I_p & 1_{p\times q} \\ 0 & I_q
\end{pmatrix},\end{equation*}
where $1_{p\times q}$ denotes the $p\times q$ matrix having all entries equal
to one;
\item[(3)] As $z\to\infty$, we have that
\begin{multline*} \Psi(z) =
(I+O(1/z))\times\\
\diag(z^{n_1}f_1(z),\ldots,z^{n_p}f_p(z),z^{-m_1}f_{p+1}(z),\ldots,z^{-m_q}f_{p+q}(z)).
\end{multline*}
\end{itemize}

Since the function $\Psi$ is defined from $Y$ by the multiplication on the
right with an $\vecn,\vecm$-independent matrix (recall that $N$ is assumed to
be constant), it satisfies exactly the same recurrence relations as the
original RH matrix $Y$, cf.\ Section \ref{subsectionrecmop}. Moreover, from the
fact that $\Psi$ has a constant jump matrix we obtain that the matrix valued
function
\begin{equation}\label{leadstodiffeq1} V_{\vecn,\vecm}(z) :=
\Psi'_{\vecn,\vecm}(z)\Psi_{\vecn,\vecm}(z)^{-1},
\end{equation}
where the prime denotes the derivative, is analytic in the full complex plane.
By using \eqref{asymptoticconditionY2}, \eqref{defPsi}, \eqref{deffs}, and
\eqref{leadstodiffeq1}, we obtain that as $z\to\infty$,
\begin{multline}\label{leadstodiffeq2}
V_{\vecn,\vecm}(z) = \left(I+\frac{(Y_1)_{\vecn,\vecm}}{z}+\cdots\right)\times
\\
\frac{N}{t(1-t)}\diag(-z+(1-t)a_1,\ldots,-z+(1-t)a_p,-tb_1,\ldots,-tb_q)\times
\\
\left(I+\frac{(Y_1)_{\vecn,\vecm}}{z}+\cdots\right)^{-1}+O(1/z).
\end{multline}
By Liouville's theorem, it then follows that $V_{\vecn,\vecm}$ is a polynomial
and hence we obtain from \eqref{leadstodiffeq2} that
\begin{align}
\label{diffeqcoeffmatrix} V_{\vecn,\vecm}(z)  =  -
\frac{N}{t(1-t)}\begin{pmatrix} zI_p - D_a & -C_{12}
\\ C_{21} & D_b
\end{pmatrix},
\end{align}
where we define the diagonal matrices
\begin{equation}\label{defD1D2}
D_a:= (1-t) \, \diag \left(a_1,\ldots,a_p \right),\qquad D_b:= t \, \diag
\left(b_1,\ldots,b_q \right),
\end{equation}
and where we define $C_{12}$ and $C_{21}$ by the partitioning of
$(Y_1)_{\vecn,\vecm}$ into blocks:
\begin{equation}\label{partitionY1}
(Y_1)_{\vecn,\vecm} =: \begin{pmatrix} C_{11} & C_{12} \\ C_{21} & C_{22}
\end{pmatrix},
\end{equation}
with the diagonal blocks being of size $p\times p$ and $q\times q$,
respectively. In other words, we have put $C_{11} :=
[c_{i,j}]_{i,j=1,\ldots,p}$, $C_{12} :=
[c_{i,j}]_{i=1,\ldots,p,j=p+1,\ldots,p+q}$, $C_{21} :=
[c_{i,j}]_{i=p+1,\ldots,p+q,j=1,\ldots,p}$, and $C_{22} :=
[c_{i,j}]_{i,j=p+1,\ldots,p+q}$.

Multiplying both sides of \eqref{leadstodiffeq1} on the right with
$\Psi_{\vecn,\vecm}(z)$, we conclude:

\begin{proposition}\label{propdiffeqHermite} (Differential equation for multiple Hermite polynomials)
The multiple Hermite polynomials satisfy the matrix differential equation
\begin{equation}\label{diffeqHermite}
\Psi'_{\vecn,\vecm}(z) = V_{\vecn,\vecm}(z)\Psi_{\vecn,\vecm}(z),
\end{equation}
with $V_{\vecn,\vecm}(z)$ given by
\eqref{diffeqcoeffmatrix}--\eqref{partitionY1}.
\end{proposition}

\subsection{Lax pair and compatibility conditions}
\label{subsectionlaxpair}

In this subsection we investigate the compatibility condition between the
differential equation \eqref{diffeqHermite} in Section \ref{subsectiondiffeq}
and the recurrence relations \eqref{forwardmatrixrecurrence},
\eqref{backwardmatrixrecurrence} in Section \ref{subsectionrecmop}. These
relations together constitute the Lax pair for multiple Hermite polynomials and
the compatibility yields nonlinear difference equations the recurrence
coefficients.

Throughout this subsection we will use the partitioning of
$(Y_1)_{\vecn,\vecm}$ in \eqref{partitionY1} and we will also use the similar
partitioning
\begin{equation}\label{partitionY1til} (Y_1)_{\vecn+\vece_k,\vecm+\vece_l} =
\begin{pmatrix} \til{C}_{11} & \til{C}_{12} \\ \til{C}_{21} & \til{C}_{22}
\end{pmatrix},
\end{equation}
where again $\til{C}_{11}$ has size $p\times p$, and so on.

The derivation of the compatibility conditions below will be merely a
straightforward calculation. The reader who wishes to avoid this kind of
calculations could take these results for granted and move directly to Section
\ref{subsectionscalarproducts} where we discuss the (elegant) scalar product
relations induced by these compatibility conditions.

\subsubsection{Forward relations}
\label{subsubsectionlaxpairforward}

Fix $k\in\{1,\ldots,p\}$ and $l\in\{1,\ldots,q\}$. The compatibility condition
between the differential equation \eqref{diffeqHermite} and the forward
recurrence relation \eqref{forwardmatrixrecurrence} is
\begin{equation}\label{laxpairforward}
(U^{k,l}_{\vecn,\vecm})'(z) =
V_{\vecn+\vece_k,\vecm+\vece_l}(z)U^{k,l}_{\vecn,\vecm}(z)-U^{k,l}_{\vecn,\vecm}(z)V_{\vecn,\vecm}(z),
\end{equation}
where the matrix $V_{\vecn,\vecm}(z)$ is given by \eqref{diffeqcoeffmatrix},
$V_{\vecn+ \vece_k, \vecm + \vece_l}(z)$ is given analogously by
\begin{align}
\label{diffeqcoeffmatrix2} V_{\vecn + \vece_k,\vecm + \vece_l}(z)  =  -
\frac{N}{t(1-t)}\begin{pmatrix} zI_p - D_a & -\tilde{C}_{12}
\\ \tilde{C}_{21} & D_b
\end{pmatrix},
\end{align}
 and $U^{k,l}_{\vecn,\vecm}(z)$ is
given by \eqref{forwardcoefficientmatrix}. Using the partitionings
\eqref{partitionY1} and \eqref{partitionY1til} we rewrite the latter in block
form as
\begin{equation} \label{Uklrewrite}
    U^{k,l}_{\vecn,\vecm}(z) = \begin{pmatrix}I+(z-1)\vece_k\vece_k^T +\til{C}_{11}
    \vece_k\vece_k^T -\vece_k\vece_k^T C_{11} & -\vece_k\vece_k^T C_{12} \\
    \til{C}_{21}\vece_k\vece_k^T & I-\vece_l\vece_l^T \end{pmatrix},
\end{equation} where the top left block is of size $p\times p$, and so
on. In the right-hand side of \eqref{Uklrewrite} we use our convention that
$\vece_k$ and $\vece_l$ denote column vectors, recall Convention
\ref{conventionek}.

Inserting \eqref{diffeqcoeffmatrix}, \eqref{diffeqcoeffmatrix2}, and
\eqref{Uklrewrite} into  \eqref{laxpairforward} we find
\begin{multline*}
\begin{pmatrix} \vece_k\vece_k^T & 0 \\ 0 & 0
\end{pmatrix} =  -\frac{N}{t(1-t)}\times\\ \left[\begin{pmatrix} zI-D_a & -\til{C}_{12} \\ \til{C}_{21} & D_b
\end{pmatrix}
\begin{pmatrix}I+(z-1)\vece_k\vece_k^T +\til{C}_{11}\vece_k\vece_k^T -\vece_k\vece_k^T C_{11} & -\vece_k\vece_k^T C_{12} \\
\til{C}_{21}\vece_k\vece_k^T & I-\vece_l\vece_l^T \end{pmatrix}\right.\\
\left. -
\begin{pmatrix}I+(z-1)\vece_k\vece_k^T +\til{C}_{11}\vece_k\vece_k^T -\vece_k\vece_k^T C_{11} & -\vece_k\vece_k^T C_{12} \\
\til{C}_{21}\vece_k\vece_k^T & I-\vece_l\vece_l^T \end{pmatrix}
\begin{pmatrix} zI-D_a & -C_{12} \\ C_{21} & D_b  \end{pmatrix}\right].
\end{multline*}

The left-hand side is independent of $z$, and therefore all terms on the
right-hand side involving $z$ or $z^2$ cancel out (this can also be checked by
direct calculation). Hence the equation reduces to
\begin{multline}\label{laxpair1234}
\frac{t(1-t)}{N}\begin{pmatrix} \vece_k\vece_k^T & 0 \\ 0 & 0
\end{pmatrix} = \\ -\begin{pmatrix} -D_a & -\til{C}_{12} \\ \til{C}_{21} & D_b  \end{pmatrix}
\begin{pmatrix}I-\vece_k\vece_k^T +\til{C}_{11}\vece_k\vece_k^T -\vece_k\vece_k^T C_{11} & -\vece_k\vece_k^T C_{12} \\
\til{C}_{21}\vece_k\vece_k^T & I-\vece_l\vece_l^T \end{pmatrix}\\
+
\begin{pmatrix}I-\vece_k\vece_k^T +\til{C}_{11}\vece_k\vece_k^T -\vece_k\vece_k^T C_{11} & -\vece_k\vece_k^T C_{12} \\
\til{C}_{21}\vece_k\vece_k^T & I-\vece_l\vece_l^T \end{pmatrix}
\begin{pmatrix} -D_a & -C_{12} \\ C_{21} & D_b  \end{pmatrix}.
\end{multline}
From the identity \eqref{laxpair1234} we obtain the following for the
respective blocks.
\begin{itemize}
\item $(1,1)$ block of \eqref{laxpair1234}:
\begin{equation*}
\frac{t(1-t)}{N}\vece_k\vece_k^T = [D_a,I-\vece_k\vece_k^T
+\til{C}_{11}\vece_k\vece_k^T -\vece_k\vece_k^T
C_{11}]+\til{C}_{12}\til{C}_{21}\vece_k\vece_k^T-\vece_k\vece_k^T C_{12}C_{21},
\end{equation*}
where $[ \cdot, \cdot]$ denotes the usual commutator of square matrices. Since
diagonal matrices commute with each other, we get
\begin{equation*}
\frac{t(1-t)}{N}\vece_k\vece_k^T = [D_a,\til{C}_{11}\vece_k\vece_k^T
-\vece_k\vece_k^T
C_{11}]+\til{C}_{12}\til{C}_{21}\vece_k\vece_k^T-\vece_k\vece_k^T C_{12}C_{21}.
\end{equation*}
This implies for the $(k,k)$ diagonal entry
\begin{equation}\label{laxpair11a}
    \left( \til{C}_{12} \til{C}_{21} \right)_{k,k} -
    \left( C_{12} C_{21} \right)_{k,k} = \frac{t(1-t)}{N},
\end{equation}
for the $(\til{k},k)$ entry, where we also use \eqref{defD1D2},
\begin{equation}\label{laxpair11b}
    \left( \til{C}_{12}  \til{C}_{21} \right)_{\til{k},k} =
    (1-t) (a_k - a_{\til{k}})  \left(\til{C}_{11}\right)_{\til{k},k},
\end{equation}
and similarly for the $(k,\til{k})$ entry,
\begin{equation}\label{laxpair11c}
    \left( C_{12} C_{21} \right)_{k,\til{k}}
        = -(1-t)(a_k - a_{\til{k}}) \left(C_{11}\right)_{k,\til{k}},
\end{equation}
for any $\til{k}\in\{1,\ldots,p\}$ with $\til{k}\neq k$.
\item $(1,2)$ block of \eqref{laxpair1234}:
\begin{multline*}
0= - (1-t)a_k \vece_k\vece_k^T C_{12} +\til{C}_{12}(I-\vece_l\vece_l^T)\\-
(I-\vece_k\vece_k^T +\til{C}_{11}\vece_k\vece_k^T -\vece_k\vece_k^T
C_{11})C_{12} - \vece_k\vece_k^T C_{12}D_b.
\end{multline*}
Evaluating the $l$th column of this equation, one obtains
\begin{multline*}
0= -(1-t) a_k \vece_k\vece_k^T C_{12}\vece_l+0 \\- (I-\vece_k\vece_k^T
+\til{C}_{11}\vece_k\vece_k^T -\vece_k\vece_k^T C_{11})C_{12}\vece_l -
tb_{l}\vece_k\vece_k^T C_{12}\vece_l,
\end{multline*}
which can be rewritten as
\begin{equation}\label{laxpair12a}
\left(I+(t b_{l}+(1-t) a_k-1)\vece_k\vece_k^T +\til{C}_{11}\vece_k\vece_k^T
-\vece_k\vece_k^T C_{11}\right) C_{12} \vece_l = 0.
\end{equation}
Similarly one can evaluate the $\til{l}$th column to obtain
\begin{equation}\label{laxpair12b}
\left(I+(tb_{\til{l}}+(1-t) a_k-1)\vece_k\vece_k^T
+\til{C}_{11}\vece_k\vece_k^T -\vece_k\vece_k^T C_{11}\right) C_{12}
\vece_{\til{l}} =
    \til{C}_{12} \vece_{\til{l}},
\end{equation}
for any $\til{l}\in\{1,\ldots,q\}$ with $\til{l}\neq l$.
\item $(2,1)$ block  of \eqref{laxpair1234}:
\begin{multline*}
0 = -\til{C}_{21}(I-\vece_k\vece_k^T +\til{C}_{11}\vece_k\vece_k^T
-\vece_k\vece_k^T C_{11})\\-D_b\til{C}_{21}\vece_k\vece_k^T - (1-t)
a_k\til{C}_{21}\vece_k\vece_k^T + (I-\vece_l\vece_l^T)C_{21}.
\end{multline*}
Similarly to the case of the $(1,2)$ block entry, one can now evaluate the
different rows of this equation. This yields
\begin{equation}\label{laxpair21a}
 \vece_l^T \til{C}_{21}  \left(I+((1-t)a_k+tb_l -1)\vece_k\vece_k^T
+ \til{C}_{11}\vece_k\vece_k^T -\vece_k\vece_k^T C_{11}\right) = 0,
\end{equation}
and
\begin{equation}\label{laxpair21b}
 \vece_{\til{l}}^T \til{C}_{21} \left(I+((1+t)a_k + tb_{\til{l}} -1)\vece_k\vece_k^T
+\til{C}_{11}\vece_k\vece_k^T -\vece_k\vece_k^T C_{11}\right) =
\vece_{\til{l}}^T  C_{21},
\end{equation}
for any $\til{l}\in\{1,\ldots,q\}$ with $\til{l}\neq l$.
\item $(2,2)$ block of \eqref{laxpair1234}: This is just the trivial relation $0=0$.
\end{itemize}

Note that by comparing \eqref{laxpair12b} and \eqref{laxpair21b}, one finds the
relations
\begin{equation}\label{laxpairsc}
 \left(\til{C}_{21} \til{C}_{12} \right)_{\til{l},\til{l}} =
 \left(C_{21} C_{12} \right)_{\til{l},\til{l}}
\end{equation}
for any $\til{l}\in\{1,\ldots,q\}$ with $\til{l}\neq l$.

\subsubsection{Backward relations}
\label{subsubsectionlaxpairbackward}

The relations in Section \ref{subsubsectionlaxpairforward} constitute in fact
only half of the available compatibility relations. The second half is obtained
from the differential equation \eqref{diffeqHermite} and the \emph{backward}
recurrence relation \eqref{backwardmatrixrecurrence}. Alternatively, these
relations may be obtained from the forward relations by using the involution
symmetry to be described in Corollary \ref{lemmainvolution2}.

We will not derive all these relations in detail. Instead, we list only the
equations that we will need in the sequel. These are the following analogue of
\eqref{laxpair11a}:
\begin{equation}\label{laxpairbackward22a}
 \left( \til{C}_{21}  \til{C}_{12}\right)_{l,l}-
 \left(C_{21}  C_{12}\right)_{l,l} = \frac{t(1-t)}{N},
\end{equation}
and the following analogues of \eqref{laxpair11c} and \eqref{laxpairsc}:
\begin{equation}\label{laxpairbackward22c}
    \left(C_{21} C_{12}\right)_{\til{l},l} =
    t(b_{l}-b_{\til{l}})  \left(C_{22} \right)_{\til{l},l},
\end{equation}
for any $\til{l}\in\{1,\ldots,q\}$ with $\til{l}\neq l$, and
\begin{equation}\label{laxpairbackwardsc}
    \left(\til{C}_{12} \til{C}_{21} \right)_{\til{k}, \til{k}} =
    \left(C_{12} C_{21} \right)_{\til{k},\til{k}}
\end{equation}
for any $\til{k}\in\{1,\ldots,p\}$ with $\til{k}\neq k$.

\subsubsection{Involution symmetry}
\label{subsectioninvolution}

As an aside, this might be a good place to write down explicitly the involution
symmetry which we referred to in Section \ref{subsubsectionlaxpairbackward}, as
well as in Remarks \ref{remarkeliminationothercase} and
\ref{remarkuseinvolutionsymmetry}.

We consider the involution \begin{equation*} \begin{array}{c} p\leftrightarrow
q,\qquad (w_{1,i},n_i)\leftrightarrow (w_{2,j},m_j).\end{array}\end{equation*}
We will temporarily denote the RH matrices corresponding to the original and
reversed RH problem with $Y_{\mathbf{w}_1,\mathbf{n};\mathbf{w}_2,\mathbf{m}}$
and $Y_{\mathbf{w}_2,\mathbf{m};\mathbf{w}_1,\mathbf{n}}$, respectively. Here
we use the vectorial notations $\mathbf{w}_1=(w_{1,1},\ldots,w_{1,p})$ and
$\mathbf{w}_2=(w_{2,1},\ldots,w_{2,q})$.

It is interesting to see what this involution means in terms of
non-intersecting Brownian motions. As already discussed before, in the multiple
Hermite case \eqref{gaussianw1}--\eqref{gaussianw2} the MOP are closely related
to non-intersecting Brownian motions in the sense of Section
\ref{subsectionintrointro}. Then from the non-intersecting Brownian motions
point of
view, the above involution reduces to (cf.\ \eqref{gaussianw1}--\eqref{gaussianw2})
\begin{equation*} \begin{array}{c} t\leftrightarrow 1-t,\qquad
p\leftrightarrow q,\qquad (a_i,n_i)\leftrightarrow
(b_j,m_j).\end{array}\end{equation*} This means that we interchange the role of
starting and ending points and reverse the direction of time. In the $tx$-plane
this corresponds to a reflection of the Brownian motions with respect to the
vertical line $t=1/2$.

In what follows we will denote
\begin{equation}\label{defJ}
J = \begin{pmatrix} 0 & I_q \\ -I_p & 0 \end{pmatrix}.
\end{equation}

\begin{lemma}\label{lemmainvolution1} (See also \cite{AvMV,DK2}.) We have
\begin{equation}\label{vglinvolution1}
Y_{\mathbf{w}_2,\mathbf{m};\mathbf{w}_1,\mathbf{n}} = J
Y^{-T}_{\mathbf{w}_1,\mathbf{n};\mathbf{w}_2,\mathbf{m}} J^{-1}
\end{equation}
where the superscript ${}^{-T}$ denotes the transposed inverse.
\end{lemma}

The proof of Lemma \ref{lemmainvolution1} can be performed in a straightforward
way. It suffices to check that both sides of \eqref{vglinvolution1} satisfy the
same RH problem, and next invoke the fact that the solution to this RH problem
is unique; we omit the details.

As an immediate corollary of Lemma \ref{lemmainvolution1}, we can check what
the involution does with the matrix $Y_1$ in \eqref{asymptoticconditionY2}. The
$Y_1$-matrices corresponding to the original and reversed Brownian motions turn
out to be related as follows.

\begin{corollary}\label{lemmainvolution2}
We have \begin{equation}\label{vglinvolution2}
(Y_1)_{\mathbf{w}_2,\mathbf{m};\mathbf{w}_1,\mathbf{n}} = -J
(Y_1)^{T}_{\mathbf{w}_1,\mathbf{n};\mathbf{w}_2,\mathbf{m}} J^{-1}.
\end{equation}
Put differently, if
\begin{equation}(Y_1)_{\mathbf{w}_1,\mathbf{n};\mathbf{w}_2,\mathbf{m}} =
\begin{pmatrix} C_{1,1} & C_{1,2} \\ C_{2,1} & C_{2,2} \end{pmatrix}\end{equation} denotes a partitioning with diagonal
blocks of size $p\times p$ and $q\times q$, respectively, then
\begin{equation}
(Y_1)_{\mathbf{w}_2,\mathbf{m};\mathbf{w}_1,\mathbf{n}} =
\begin{pmatrix} -C_{2,2}^T & C_{1,2}^T \\ C_{2,1}^T & -C_{1,1}^T \end{pmatrix}.
\end{equation}
\end{corollary}
Using Corollary \ref{lemmainvolution2}, the compatibility relations in Section
\ref{subsubsectionlaxpairbackward} may be immediately deduced from those in
Section \ref{subsubsectionlaxpairforward}.

\subsection{Scalar product relations}
\label{subsectionscalarproducts}

Now we discuss the scalar product relations induced by the compatibility
conditions in Section \ref{subsectionlaxpair}. These relations are given by the
following proposition.

\begin{proposition} \label{stscalarproducts} (Scalar product relations)
Let $\vecn$ and $\vecm$ with $|\vecn|=|\vecm|$ be fixed and consider the
partitioning \eqref{partitionY1} for the matrix $Y_1 = (Y_1)_{\vecn,\vecm}$.
Then we have the relations
\begin{equation}\label{scalarproductsn}
    \left(C_{12} C_{21} \right)_{k,k} = t(1-t)\frac{n_{k}}{N},
\end{equation}
for all $k\in\{1,\ldots,p\}$,
\begin{equation}\label{scalarproductsm}
    \left(C_{21} C_{12} \right)_{l,l} = t(1-t)\frac{m_{l}}{N},
\end{equation}
for all $l\in\{1,\ldots,q\}$,
\begin{equation}\label{scalarproductsskewtop}
    \left(C_{12} C_{21}\right)_{k,\til{k}} = - (1-t)(a_k-a_{\til{k}})
        \left(C_{11} \right)_{k,\til{k}},
\end{equation}
for all $k,\til{k}\in\{1,\ldots,p\}$ with $k\neq \til{k}$, and
\begin{equation}\label{scalarproductsskewbottom}
    \left(C_{21} C_{12} \right)_{l, \til{l}} = -t(b_{l}-b_{\til{l}})
        \left(C_{22} \right)_{l,\til{l}},
\end{equation}
for all $l,\til{l}\in\{1,\ldots,q\}$ with $l\neq \til{l}$.
\end{proposition}

\bewijs. If $\vecn = \mathbf{0}$ and $\vecm = \mathbf{0}$ then the solution
$Y_{\mathbf{0},\mathbf{0}}$ of the RH problem
\eqref{defjumpmatrix0}--\eqref{asymptoticconditionY0} is upper triangular, so
that $C_{21}$ is the zero-matrix and the relations \eqref{scalarproductsn} and
\eqref{scalarproductsm} hold in that case. For arbitrary $\vecn$ and $\vecm$,
these relations then follow from an induction argument based on
\eqref{laxpair11a}, \eqref{laxpairsc}, \eqref{laxpairbackward22a}, and
\eqref{laxpairbackwardsc}.

The  two relations \eqref{scalarproductsskewtop} and
\eqref{scalarproductsskewbottom}  are simply \eqref{laxpair11c} and
\eqref{laxpairbackward22c}. $\bol$

We call \eqref{scalarproductsn}--\eqref{scalarproductsskewbottom} scalar
product relations, since these relations can be interpreted as giving the
scalar products between  row and column vectors of $C_{12}$ and $C_{21}$. The
relations \eqref{scalarproductsn} and \eqref{scalarproductsm} then involve rows
and columns with the same index, while the relations
\eqref{scalarproductsskewtop} and \eqref{scalarproductsskewbottom} involve rows
and columns with different indices. Moreover, note that these last two
relations determine the off-diagonal entries of the diagonal blocks $C_{11}$
and $C_{22}$ in \eqref{partitionY1} in terms of scalar products of the entries
in the off-diagonal blocks $C_{12}$ and $C_{21}$ in \eqref{partitionY1}.

Let us illustrate the above scalar product relations for the $4\times 4$ case
(i.e., $p=q=2$). Then the partitioning \eqref{partitionY1} is, when written out
in full,
\begin{equation*} (Y_1)_{\vecn,\vecm} =
\begin{pmatrix} C_{11} & C_{12} \\ C_{21} & C_{22}
\end{pmatrix} =
\begin{pmatrix} c_{1,1} & c_{1,2} & c_{1,3} & c_{1,4} \\
    c_{2,1} & c_{2,2} & c_{2,3} & c_{2,4} \\
    c_{3,1} & c_{3,2} & c_{3,3} & c_{3,4} \\
    c_{4,1} & c_{4,2} & c_{4,3} & c_{4,4}
    \end{pmatrix}.
\end{equation*}
Then \eqref{scalarproductsn} and \eqref{scalarproductsm} give us the four
relations
\begin{align}
\label{scalarproducts1}\begin{pmatrix} c_{1,3} & c_{1,4}
\end{pmatrix}
\begin{pmatrix} c_{3,1} \\ c_{4,1}
\end{pmatrix} & = t(1-t)\frac{n_1}{N},
\\
\label{scalarproducts2}\begin{pmatrix} c_{2,3} & c_{2,4}
\end{pmatrix}
\begin{pmatrix} c_{3,2} \\ c_{4,2}
\end{pmatrix} &= t(1-t)\frac{n_2}{N},
\\
\label{scalarproducts3}\begin{pmatrix} c_{3,1} & c_{3,2}
\end{pmatrix}
\begin{pmatrix} c_{1,3} \\ c_{2,3}
\end{pmatrix} &= t(1-t)\frac{m_1}{N},
\\
\label{scalarproducts4}\begin{pmatrix} c_{4,1} & c_{4,2}
\end{pmatrix}
\begin{pmatrix} c_{1,4} \\ c_{2,4}
\end{pmatrix} &= t(1-t)\frac{m_2}{N}.
\end{align}
In other words, the number of Brownian particles $n_1$, $n_2$ starting from the
first and second starting point $a_1$ and $a_2$, respectively, is (up to a
common factor $t(1-t)/N$) directly expressed by the scalar products
\eqref{scalarproducts1}, \eqref{scalarproducts2}, while the number of Brownian
particles $m_1$, $m_2$ arriving at the first and second endpoints $b_1$ and
$b_2$, respectively, is expressed by the scalar products
\eqref{scalarproducts3}, \eqref{scalarproducts4} (up to the same common
factor).

The scalar product relations \eqref{scalarproducts1}--\eqref{scalarproducts4}
can be nicely expressed in terms of the Hadamard product matrix $H$ of
\eqref{hadamardproductexample}. Indeed, they express that the top rightmost
$2\times 2$ block of this matrix, i.e.,
\begin{equation}\label{hadamardproductexamplesubmatrix}
    H_{12} := \begin{pmatrix}c_{1,3}c_{3,1} & c_{1,4}c_{4,1}\\
    c_{2,3}c_{3,2} & c_{2,4}c_{4,2}\end{pmatrix},
\end{equation}
has fixed row sums equal to $t(1-t)\frac{n_1}{N}$, $t(1-t)\frac{n_2}{N}$, and
fixed column sums equal to $t(1-t)\frac{m_1}{N}$, $t(1-t)\frac{m_2}{N}$. It
follows from this that the matrix \eqref{hadamardproductexamplesubmatrix} is
fully determined by only one of its entries, since the other  entries then
follow from these row and column sum relations. This implies in particular the
first three relations in Proposition~\ref{propconnectionsrecurrencecoeff4x4}.

In the case of general $p$ and $q$, we have the following result.

\begin{proposition} \label{strowcolumnsums} (Row and column sum relations)
Let $\vecn$ and $\vecm$ with $|\vecn|=|\vecm|$ be fixed and consider the
Hadamard product matrix $H$ in \eqref{hadamardproduct}. Then the top rightmost
$p\times q$ block of the matrix $H$, i.e., the submatrix
\begin{equation}\label{hadamardproductsubmatrix}
H_{12} = \begin{pmatrix}c_{1,p+1}c_{p+1,1} & \ldots & c_{1,p+q}c_{p+q,1}\\
\vdots & & \vdots \\ c_{p,p+1}c_{p+1,p} & \ldots &
c_{p,p+q}c_{p+q,p}\end{pmatrix}
\end{equation}
has fixed row sums equal to $t(1-t)\frac{n_k}{N}$, $k=1,\ldots,p$, and fixed
column sums equal to $t(1-t)\frac{m_l}{N}$, $l=1,\ldots,q$.
\end{proposition}

\bewijs. These relations are equivalent with the scalar product relations
\eqref{scalarproductsn} and \eqref{scalarproductsm}. $\bol$

\subsection{Spectral curve interpretation}
\label{subsectionspectralcurve}

As an aside, let us discuss an equivalent formulation of Proposition
\ref{strowcolumnsums} in terms of the so-called spectral curve. Let
\begin{equation}\label{defspectralcurve}
P_{\vecn,\vecm}(\xi,z) :=  \det \left(\xi I+\frac{1}{n} V_{\vecn,\vecm}(z)
\right),
\end{equation}
where $n = |\vecn| = |\vecm|$. The zero set of the polynomial
$P_{\vecn,\vecm}(\xi,z)$ in \eqref{defspectralcurve}  defines an algebraic
curve (the spectral curve).

For example, when $p=q=2$ we can use
\eqref{diffeqcoeffmatrix}--\eqref{partitionY1} and \eqref{defspectralcurve} to
see that
\begin{multline}
  P_{\vecn,\vecm}(\xi,z)  =   \det \left(\xi I_4 - \frac{N}{n t(1-t)}\begin{pmatrix} zI_2 - D_a & -C_{12}
\\ C_{21} & D_b \end{pmatrix} \right)
\\
\label{vglspectralcurve2}
    =  \left(\frac{N}{nt(1-t)}\right)^4
\det\begin{pmatrix}
\til{\xi}-z+(1-t)a_1 & 0 & c_{1,3} & c_{1,4} \\
0 & \til{\xi}-z+(1-t)a_2 & c_{2,3} & c_{2,4} \\
-c_{3,1} & -c_{3,2} & \til{\xi}-tb_1 & 0 \\
-c_{4,1} & -c_{4,2} & 0 & \til{\xi}-tb_2
\end{pmatrix},
\end{multline}
where $\til{\xi}:=  t(1-t) \frac{n}{N} \xi$.

The spectral curve $P_{\vecn,\vecm}(\xi,z) = 0$ in \eqref{defspectralcurve}
defines an algebraic curve which is of degree $p+q$ in the variable $\xi$.
After resolution of singularities, this curve defines a Riemann surface which
can be used as a tool in the steepest descent analysis of the RH problem. This
is in fact a $(p+q)$-sheeted Riemann surface where the $i$th sheet corresponds
to the $i$th solution function $\xi(z) = \xi_i(z)$, $i=1,\ldots,p+q$.

Now it turns out \cite{ABK,BK2,BK3,DKV,O} that an important role in the
steepest descent analysis is played by the asymptotic expansions for
$z\to\infty$ of the functions $\xi_i(z)$. It was observed for some special
cases that the following asymptotic expansions hold.

\begin{proposition} \label{stasymptoticsspectralcurve} (Asymptotic expansions of the spectral curve)
Let $p$ and $q$ be arbitrary. Then when appropriately labeled, the $p+q$
branches $\xi(z) = \xi_i(z)$, $i=1,\ldots,p+q$ of the algebraic curve
$P_{\vecn,\vecm}(\xi,z) = 0$ in \eqref{defspectralcurve} behave as follows:
\begin{align}
\label{sheetexpansiona} \xi_k(z) & = \frac{N}{nt(1-t)} z - \frac{Na_k}{n t}-
    \frac{n_k}{n} \frac{1}{z} +O\left(\frac{1}{z^2}\right),\\
\label{sheetexpansionb} \xi_{p+l}(z) & =  \frac{N b_l}{n t}+
    \frac{m_l}{n} \frac{1}{z} +O\left(\frac{1}{z^2}\right),
\end{align}
as $z\to\infty$, for any $k\in\{1,\ldots,p\}$ and $l\in\{1,\ldots,q\}$.
\end{proposition}

Note in particular that the $1/z$ terms in the $z\to\infty$ expansion, up to a
factor $\pm 1$, are precisely the fractions of Brownian particles starting from
the $k$th starting point $a_k$ or arriving in the $l$th ending point $b_l$,
respectively.

\bewijs. We give the proof for the case $p=q=2$; the proof for general $p$ and
$q$ is similar. From \eqref{vglspectralcurve2}, we see that the only way to
have $P_{\vecn,\vecm}(\xi,z)=0$ when $z\to\infty$ is that one of the diagonal
entries in \eqref{vglspectralcurve2} vanishes as $z \to \infty$. This yields
the required linear and constant terms in \eqref{sheetexpansiona} and
\eqref{sheetexpansionb}.

To find the $1/z$ terms in \eqref{sheetexpansiona} and \eqref{sheetexpansionb},
one can propose a term of the form $c/z$ for some yet unknown $c\in\mathbb C $, and
subsequently require the coefficient of the $z$ term in the Laurent expansion
of \eqref{vglspectralcurve2} to vanish. For $k=1,2$, we then get for
$\til{\xi}_k =  t(1-t) \frac{n}{N} \xi_k$
\[ \til{\xi}_k(z)
    = z - (1-t) a_k - \frac{c_{k,3}c_{3,k}+c_{k,4}c_{4,k}}{z} + O\left(\frac{1}{z^2}\right) \]
which by \eqref{scalarproducts1} and \eqref{scalarproducts2} leads to
\eqref{sheetexpansiona}. The $1/z$ term in \eqref{sheetexpansionb} follows in a
similar way from \eqref{scalarproducts3} and \eqref{scalarproducts4}. $\bol$

We end this subsection with a brief discussion. We described in
Proposition~\ref{stscalarproducts} the scalar product relations induced by the
Lax pair for multiple Hermite polynomials. In case where $p=1$ these relations
fully describe the off-diagonal entries of $(Y_1)_{\vecn,\vecm}$. Indeed, the
block entries $C_{12}$ and $C_{21}$ in \eqref{partitionY1} then have only one
row and one column, respectively. The row and column sum relations
\eqref{scalarproductsn} and \eqref{scalarproductsm} then give explicit
expressions for the products $c_{i,j}c_{j,i}$, $i<j$. With a little bit more
work one can also obtain explicit expressions for the individual entries
$c_{i,j}$; this leads to the expressions in \cite{BK2}. Similar remarks hold
when $q=1$.

In the case of general $p$ and $q$, however, we are not able to find explicit
expressions for the entries of $(Y_1)_{\vecn,\vecm}$ anymore. In this case, the
best that we could obtain is showing that the compatibility relations in
Section \ref{subsectionlaxpair} can be used to determine the off-diagonal
entries of $(Y_1)_{\vecn+\vece_k,\vecm+\vece_l}$ in terms of the off-diagonal
entries of $(Y_1)_{\vecn,\vecm}$. This yields a recursive scheme for computing
the off-diagonal entries of $(Y_1)_{\vecn,\vecm}$ by induction on $|\vecn| =
|\vecm|$. Unfortunately, we were not able to derive any useful asymptotic
information from this recursive scheme; therefore we will not state these
relations here.

\subsection{Completion of the proofs of Propositions \ref{proprecurrence4x4} and \ref{propconnectionsrecurrencecoeff4x4}}
\label{subsectionlaxpairfurtherconsequences}

\subsubsection{Proof of Proposition \ref{proprecurrence4x4}}
\label{subsubsectiondiagrec}

As already noted above, Proposition \ref{proprecurrence4x4} follows as a
special case of Proposition \ref{propforwardrecurrence}, except for the
explicit form of the diagonal recurrence coefficients in the first term in the
right hand sides of each of
\eqref{examplerecurrenceintro1}--\eqref{examplerecurrenceintro4}. The latter
expressions are only valid in the multiple Hermite case.

To derive these expressions for the diagonal recurrence coefficients, we will
again use the compatibility conditions in Section \ref{subsectionlaxpair}. We
proceed as follows. By evaluating the $k$th row of \eqref{laxpair12a}, one
obtains
\begin{equation*}
((1-t) a_k + tb_{l}\vece_k^T +\vece_k^T\til{C}_{11}\vece_k\vece_k^T -\vece_k^T
C_{11}) C_{12} \vece_{l} = 0,
\end{equation*}
which is
\begin{equation*}
    \left((1-t) a_k+ tb_{l}+ \left(\til{C}_{11}\right)_{k,k}\right) \left(C_{12}\right)_{k,l}
    -\left(C_{11} C_{12}\right)_{k,l} = 0,
\end{equation*}
or equivalently
\begin{equation*}
((1-t) a_k + tb_{l}+ \left(\til{C}_{1,1}\right)_{k,k}-
\left(C_{11}\right)_{k,k})
    \left(C_{12}\right)_{k,l} = \sum_{\til{k}\neq k}
    \left(C_{11}\right)_{k,\til{k}} \left(C_{12}\right)_{\til{k},l}.
\end{equation*}
From this expression, we obtain the relation
\begin{equation*}
    \left(C_{11}\right)_{k,k}-\left(\til{C}_{11}\right)_{k,k}
    = (1-t) a_k+ tb_{l} -\frac{\sum_{\til{k}\neq k}
    \left(C_{11}\right)_{k,\til{k}} \left(C_{12}\right)_{\til{k},l}}{\left(C_{12}\right)_{k,l}}.
\end{equation*}
Rewriting this in terms of the entries $c_{i,j}$, $\til{c}_{i,j}$ of
$(Y_1)_{\vecn,\vecm}$ and $(Y_1)_{\vecn+\vece_k,\vecm+\vece_l}$, we have the
equivalent expression
\begin{align}
c_{k,k}-\til{c}_{k,k} \label{reccoefdiagonaltypelax}&=
(1-t)a_k+tb_l-\frac{\sum_{\til{k}\neq k}
    c_{k,\til{k}}c_{\til{k},p+l}}{c_{k,p+l}},
\end{align}
for the diagonal recurrence coefficients $c_{k,k}-\til{c}_{k,k}$. Note that the
expressions in the first term in the right hand sides of each of
\eqref{examplerecurrenceintro1}--\eqref{examplerecurrenceintro4} follows as a
special case of \eqref{reccoefdiagonaltypelax}. This ends the proof of
Proposition \ref{proprecurrence4x4}.

By the way, note that the right-hand side of \eqref{reccoefdiagonaltypelax}
contains entries of $(Y_1)_{\vecn,\vecm}$ only, and so it is more convenient
than \eqref{reccoefdiagonaltype} which also contains an entry of
$(Y_2)_{\vecn,\vecm}$. In addition, \eqref{reccoefdiagonaltypelax} has only
about half as many terms as \eqref{reccoefdiagonaltype}.

Recall also that the left-hand side of \eqref{reccoefdiagonaltypelax} depends
on $l$ since $\til{c}_{k,k}$ is an entry of the matrix
$(Y_1)_{\vecn+\vece_k,\vecm+\vece_l}$. $\bol$

\subsubsection{Proof of Proposition \ref{propconnectionsrecurrencecoeff4x4}}
\label{subsubsectiontwodeterminants}

As already noted above, Proposition \ref{propconnectionsrecurrencecoeff4x4}
follows from the more general row and column sum relations in Proposition
\ref{strowcolumnsums}. The only thing that remains is to prove
\eqref{fourthrelation}.

Assume again that $p=q=2$. We will compute the determinant
\begin{equation}\label{determinanttwoways}
    \det C_{12} \det C_{21}
    = \det \begin{pmatrix} c_{1,3} & c_{1,4}\\ c_{2,3} & c_{2,4} \end{pmatrix}
    \det \begin{pmatrix} c_{3,1} & c_{3,2}\\ c_{4,1} & c_{4,2} \end{pmatrix}
\end{equation}
first via
\begin{align}
    \det C_{12} \det C_{21} & = \nonumber
    \det (C_{12} C_{21}) \\
    & = \nonumber
     \det \left(\begin{pmatrix} c_{1,3} & c_{1,4}\\ c_{2,3} &  c_{2,4}
    \end{pmatrix}\begin{pmatrix} c_{3,1} & c_{3,2}\\ c_{4,1} &
    c_{4,2} \end{pmatrix} \right) \\
    \nonumber
    &= \det \begin{pmatrix} t(1-t)\frac{n_1}{N} & -(1-t)(a_{1}-a_{2})c_{1,2} \\
    (1-t)(a_1-a_{2})c_{2,1} & t(1-t)\frac{n_2}{N} \end{pmatrix}\\
    \label{thirdunderbraced}
    & = t^2(1-t)^2\frac{n_1n_2}{N^2}+ (1-t)^2(a_1-a_{2})^2 c_{1,2}c_{2,1},
\end{align}
where we made use of the scalar product relations \eqref{scalarproductsn} and
\eqref{scalarproductsskewtop}.

We similarly compute
\begin{align}
    \det C_{12} \det C_{21} & = \nonumber
    \det (C_{21} C_{12}) \\
    & = \nonumber
    \det \left(\begin{pmatrix} c_{3,1} & c_{3,2}\\ c_{4,1} & c_{4,2} \end{pmatrix}
    \begin{pmatrix} c_{1,3} & c_{1,4}\\ c_{2,3} & c_{2,4} \end{pmatrix}\right) \\
    \nonumber
    &= \det \begin{pmatrix} t(1-t)\frac{m_1}{N} & t(b_{1}-b_{2})c_{3,4} \\
    -t(b_{1}-b_{2})c_{4,3} & t(1-t)\frac{m_2}{N}
    \end{pmatrix}\\
    \label{thirdunderbracedbis} & = t^2(1-t)^2\frac{m_1m_2}{N^2}+
    t^2(b_{1}-b_{2})^2c_{3,4}c_{4,3},
\end{align}
where we now used  \eqref{scalarproductsm} and
\eqref{scalarproductsskewbottom}. By equating \eqref{thirdunderbraced} and
\eqref{thirdunderbracedbis} we obtain a relation between the off-diagonal
recurrence coefficients $c_{1,2}c_{2,1}$ and $c_{3,4}c_{4,3}$
\[ (1-t)^2(a_1-a_{2})^2 c_{1,2}c_{2,1}
    - t^2(b_{1}-b_{2})^2c_{3,4}c_{4,3} =
    t^2(1-t)^2\frac{n_1n_2-m_1m_2}{N^2} \]
which for the special case $n_1 = m_1$ and $n_2 = m_2$ reduces to
\eqref{fourthrelation}. This ends the proof of Proposition
\ref{propconnectionsrecurrencecoeff4x4}. $\bol$

\end{document}